\documentclass[preprint,12pt]{elsarticle}

\usepackage{amssymb}
\usepackage{amsmath}
\usepackage{anyfontsize}
\usepackage[pagewise]{lineno}
\usepackage{epsfig,latexsym,amssymb, amsmath,amscd,amsxtra}
\usepackage{subfigure}
\usepackage{graphicx}
\usepackage{algorithm}  
\usepackage{algpseudocode}    
\usepackage{bm}
\usepackage{multirow}
\usepackage[final]{changes} %禁用修订，输出最终修订完成的版本
\usepackage{subfigure}
\usepackage{color}
\usepackage[numbers]{natbib}
\usepackage{indentfirst}
\usepackage{comment}
\usepackage[misc]{ifsym}
\usepackage{makecell} % 用于单元格内换行
\newcommand{\ignore}[1]{}
\newcommand{\RNum}[1]{\uppercase\expandafter{\romannumeral #1\relax}}
\newtheorem{remark}{Remark}
\usepackage{enumitem}
\setdeletedmarkup{}
\usepackage{etoolbox}
\makeatletter
\def\ps@pprintTitle{%
  \let\@oddhead\@empty
  \let\@evenhead\@empty
  \def\@oddfoot{\centerline{\thepage}}%
  \let\@evenfoot\@oddfoot}
\makeatother

\begin{document}

\begin{frontmatter}



\title{Adaptive feature capture method for solving  partial differential
equations with near singular solutions} %% Article title


\author[1]{Yangtao Deng}
\ead{ytdeng1998@foxmail.com}

\author[1]{Qiaolin He \corref{cor}}
\ead{qlhejenny@scu.edu.cn}

\author[2,3]{Xiaoping Wang \corref{cor}}
\ead{wangxiaoping@cuhk.edu.cn}

\cortext[cor]{Corresponding author}

\address[1]{School of Mathematics, Sichuan University, Chengdu, 610065, China}

\address[2]{School of Science and Engineering, The Chinese University of Hong Kong, Shenzhen, Guangdong, 518172, China}
\address[3]{Shenzhen International Center for Industrial and Applied Mathematics, Shenzhen Research Institute of Big Data,
Guangdong, 518172, China}

%% Abstract
\begin{abstract}
Partial differential equations (PDEs) with near singular solutions pose significant challenges for traditional numerical methods, particularly in complex geometries where mesh generation and adaptive refinement become computationally expensive. Although deep-learning-based approaches, such as Physics-Informed Neural Networks (PINNs) and the Random Feature Method (RFM), offer mesh-free alternatives, they often lack adaptive resolution in critical regions, limiting their accuracy for solutions with steep gradients or singularities.
In this work, we propose the Adaptive Feature Capture Method (AFCM), a novel machine learning framework that adaptively redistributes neurons and collocation points in high-gradient regions to enhance local expressive power. Inspired by adaptive moving mesh techniques, AFCM uses the gradient norm of an approximate solution as a monitor function to guide the reinitialization of feature function parameters. This ensures that partition hyperplanes and collocation points cluster where they are most needed, achieving higher resolution without increasing computational overhead.
The AFCM extends the capabilities of RFM to handle PDEs with near-singular solutions while preserving its mesh-free efficiency. Numerical experiments demonstrate the method's effectiveness in accurately resolving near-singular problems with a performance that is better than that of the traditional finite element method in terms of accuracy and efficiency. 
AFCM offers a robust and scalable approach to solving challenging PDEs in scientific and engineering applications.

\end{abstract}

%% Keywords

\end{frontmatter}



\section{Introduction}\label{sec01}
\setcounter{equation}{0}

Partial differential equations (PDEs) are widely applied in diverse fields such as physics, engineering, economics, and biology \cite{debnath2005nonlinear,achdou2014partial,leung2013systems}.  Traditional numerical methods, including finite difference \cite{leveque2007finite}, finite volume\cite{moukalled2016finite}, and finite element methods \cite{thomee2007galerkin}, have made significant theoretical and practical contributions to solving PDEs. However, these methods face notable challenges. For instance, complex geometries often lead to distorted mesh elements, which degrade computational accuracy and efficiency \cite{leveque2007finite,zienkiewicz2005finite,rajendran2010technique,blazek2015computational}. 
% Additionally, the discretization process can introduce biases that diverge from the intrinsic mathematical nature of PDEs, resulting in approximate models that fail to fully capture the underlying physical behaviors [9]. Extending solutions beyond the computational domain further complicates matters, as extrapolation tends to be inaccurate, while grid regeneration is computationally expensive [4].

In contrast, the success of deep learning in computer vision and natural language processing \cite{bengio2017deep}  has spurred interest in its application to scientific computing. Neural networks, with their universal approximation capabilities \cite{cybenko1989approximation},  have been explored for solving ordinary and partial differential equations (ODEs and PDEs) \cite{weinan2017deep,han2018solving,weinan2018deep,sirignano2018dgm,zang2020weak,raissi2019physics,weinan2021algorithms,lin2024adaptive}. Various deep-learning-based approaches have emerged, such as the Deep Ritz Method (DRM) \cite{weinan2018deep}, Deep Galerkin Method (DGM) \cite{sirignano2018dgm}, Physics-Informed Neural Networks (PINNs) \cite{raissi2019physics}, and Weak Adversarial Networks (WAN) \cite{zang2020weak}. These methods offer mesh-free alternatives, circumventing the need for computationally intensive mesh generation. However, a critical limitation of these approaches is the lack of reliable error estimation. Without knowledge of the exact solution, numerical approximations often fail to exhibit clear convergence trends, even as network parameters increase \cite{zhang2024physics}, raising concerns about their reliability in scientific and engineering applications.

Recent studies highlight the potential of randomized neural networks, such as the Extreme Learning Machine (ELM) \cite{huang2006extreme} or Random Feature Method (RFM) \cite{neal2012bayesian,rahimi2007random}, for solving ODEs and PDEs. ELM, a single-hidden-layer feedforward network, randomly initializes hidden-layer weights and biases while analytically optimizing output-layer weights via least squares \cite{huang2015trends}. This architecture eliminates the need for iterative training of hidden-layer parameters, offering significant computational efficiency over conventional deep networks. As a mesh-free universal approximator \cite{ huang2006extreme, chen2022bridging, chen2023random, huang2006universal} , ELM excels in handling PDEs in complex geometries. Extensions like the Physics-Informed ELM (PIELM) \cite{dwivedi2020physics} have further demonstrated its utility in solving differential equations \cite{chen2022bridging,chen2023random,calabro2021extreme,dong2021local,dwivedi2020physics,fabiani2021numerical,yang2018novel,wang2024extreme,dang2024adaptive,SUN2024115830,huang2024adaptive}.

Building on these advances, the Random Feature Method (RFM) \cite{chen2022bridging} combines partition of unity (PoU) with random feature functions to solve linear PDEs in complex geometries, achieving high accuracy in both space and time. However, for PDEs with near-singular solutions, RFM may struggle due to insufficient local expressive power in high-gradient regions.  In classic adaptive numerical methods, the mesh as well as the domain may be refined or decomposed, in order to improve the accuracy. Adaptive mesh refinement techniques, such as the moving mesh method \cite{renwang2000, huang2010},  address this issue by dynamically clustering mesh points in critical regions using monitor functions (e.g., solution gradients or error estimates). Therefore, it is desirable to transfer such important and successful strategies to the field of neural-network-based solutions.  

In this paper,  we propose the Adaptive Feature Capture Method (AFCM), an extension of RFM that enhances resolution in high-gradient regions without increasing computational overhead.
The AFCM leverages the gradient norm of an approximate solution to redistribute feature function hyperplanes and collocation points, concentrating them in regions of steep gradients. This adaptive refinement is iteratively applied until convergence, ensuring accurate approximations even for near-singular solutions. Crucially, AFCM preserves the mesh-free nature of RFM, making it suitable for complex geometries while maintaining computational efficiency.  Numerical experiments  demonstrate significant accuracy improvement and adaptive capability across various near-singular PDE problems. By gradient-driven redistribution of partition hyperplanes, resampling of collocation points, and adjustment of shape parameters, the method enhances local expressive power in high-gradient regions without increasing overall computational complexity. AFCM performs exceptionally well in problems with peaks, line-type singularities, and time-dependent singularities, achieving error reductions of up to 10 orders of magnitude. In geometric near-singularity problems (e.g., corners), the method still improves accuracy, though the extent of improvement is relatively limited, suggesting that further optimization of local adaptive strategies and feature function construction may be needed. Overall, AFCM provides an efficient, adaptive, and easily implementable numerical framework for mesh-free solving of near-singular PDEs.

The remainder of this paper is organized as follows: Section 2 introduces the RFM and TransNet initialization. Section 3 details the AFCM algorithm. Section 4 presents numerical experiments, and Section 5 concludes with remarks and future directions.

\section{Random Feature Method and the General Neural Feature Space}
\label{sec02}

\subsection{Random Feature Method}
\label{sec02:sec01}

Consider the following  linear  boundary-value problem
\begin{equation}
	\begin{aligned}
\begin{cases}\mathcal{L} \boldsymbol{\phi}(\boldsymbol{x})=\boldsymbol{f}(\boldsymbol{x}), & \boldsymbol{x} \in \Omega, \\ \mathcal{B} \boldsymbol{\phi}(\boldsymbol{x})=\boldsymbol{g}(\boldsymbol{x}), & \boldsymbol{x} \in \partial \Omega,\end{cases}
	\end{aligned}\label{eq:spde}
\end{equation}
where $\Omega$ is a bounded  spatial domain with the boundary  $\partial \Omega$. The $\mathcal{L}$ and $\mathcal{B}$ are linear differential and  boundary operators, respectively. We use $d$ and $d_{\phi}$ to denote the dimensions of $\boldsymbol{x}=(x_1,x_2,...,x_{d})$ and $\boldsymbol{\phi} = (\phi_1, \phi_2,...,\phi_{d_\phi})$, respectively. 

%The RFM \cite{chen2022bridging} primarily includes the following steps. 
In RFM, the domain $\Omega$ is partitioned into $M_p$ non-overlapping subdomains ${\Omega_n}$, each centered at $\boldsymbol{x}_n$, such that $\Omega =\cup_{n=1}^{M_p} \Omega_n$. For each $\Omega_n$,  \added{$\boldsymbol{r}_n\in \mathbb{R}^d$ is defined as the half-width vector of the minimum axis-aligned bounding box of $\Omega_n$.} RFM applies a linear transformation
\begin{equation}
	\begin{aligned}
\tilde{\boldsymbol{x}}=\frac{1}{\boldsymbol{r}_n}\left(\boldsymbol{x}-\boldsymbol{x}_n\right), \quad n=1, \cdots, M_p,
\end{aligned}\label{eq:normalize}
\end{equation}
\deleted{to  map $\Omega_n$ into $[-1,1]^d$, where $\boldsymbol{r}_n\in \mathbb{R}^d$ represents the radius of $\Omega_n$.}\added{to maps $\Omega_n$ into the normalized coordinate space. The transformed subdomain $\Omega_n$ is strictly enclosed within the standard hypercube $[-1,1]^d$, and the use of the minimum bounding box enables $\Omega_n$ to fit inside $[-1,1]^d$ as tightly as possible.} 
The PoU function $\psi_n$ is constructed in the normalized coordinate $\tilde{\boldsymbol{x}}$. For $d=1$, two commonly used PoU functions are
\begin{eqnarray}
%	\begin{aligned}
& &\psi_n(\boldsymbol{x})  = \mathbb{I}_{[-1,1]}(\tilde{\boldsymbol{x}}), \label{eq:psi1} %\\ 
  %& & \nonumber \\
\end{eqnarray}
and 
\begin{eqnarray}
& &  \psi_n(\boldsymbol{x})  = \mathbb{I}_{\left[-\frac{5}{4},-\frac{3}{4}\right]}(\tilde{\boldsymbol{x}}) \frac{1+\sin (2 \pi \tilde{\boldsymbol{x}})}{2}+\mathbb{I}_{\left[-\frac{3}{4}, \frac{3}{4}\right]}(\tilde{\boldsymbol{x}})\label{eq:psi2}\\
& & +\mathbb{I}_{\left[\frac{3}{4}, \frac{5}{4}\right]}(\tilde{\boldsymbol{x}}) \frac{1-\sin (2 \pi \tilde{\boldsymbol{x}})}{2}, \nonumber
%\end{aligned}
\end{eqnarray}
where $\mathbb{I}_{[a,b]}(\boldsymbol{x})=1, \boldsymbol{x}\in[a,b]$ and $a\le b$. For $d>1$, the PoU function $\psi_n(\boldsymbol{x})$ is defined as $\psi_n(\boldsymbol{x})=\Pi_{i=1}^{d} \psi_n\left(x_i\right)$.

Next, a random feature function $\varphi_{n j}$ on $\Omega_n$ is constructed using a two-layer neural network
\begin{equation}
	\begin{aligned}
\varphi_{n j}(\boldsymbol{x})=\sigma\left(\boldsymbol{W}_{n j}\cdot \tilde{\boldsymbol{x}}+b_{n j}\right), \quad j=1,2, \cdots, J_n,
\end{aligned}\label{eq:basicfunction}
\end{equation}
where $\sigma$ is the nonlinear activation function. The $\boldsymbol{W}_{n j}$ and $b_{n j}$ are randomly generated from the uniform distribution $\mathbb{U}\left(-R_m, R_m\right)$ and then fixed. The $R_m$ controls the magnitude of the parameters, and $J_n$ is the number of random feature functions.
The approximate solution in RFM is formed by a linear combination of the random feature functions and the PoU functions as follows
\begin{equation}
	\begin{aligned}
\tilde{\boldsymbol{\phi}}(\boldsymbol{x})=\left(\sum_{n=1}^{M_p} \psi_n(\boldsymbol{x}) \sum_{j=1}^{J_n} u_{n j}^{1} \varphi_{n j}^{1}(\boldsymbol{x}),\cdots,\sum_{n=1}^{M_p} \psi_n(\boldsymbol{x}) \sum_{j=1}^{J_n} u_{n j}^{d_{\phi}} \varphi_{n j}^{d_{\phi}}(\boldsymbol{x})\right)^T,
\end{aligned}\label{eq:approximatesolution}
\end{equation}
where $u_{n j}^{i} \in \mathbb{R}$ are unknowns to be determined, and $M=\sum_{n=1}^{M_p}J_n$ denotes the degree of freedom.  

Then, the linear least-squares method is utilized to minimize the loss function defined by
\begin{equation}
	\begin{aligned}
\operatorname{Loss}\left(\left\{u_{n j}^{i}\right\}\right)=&\sum_{n=1}^{M_p}\left(\sum_{q=1}^{Q_n}\left\|\boldsymbol{\lambda}_{n, q}\left(\mathcal{L} \tilde{\boldsymbol{\phi}}\left(\boldsymbol{x}_{q}^{n}\right)-\boldsymbol{f}\left(\boldsymbol{x}_{q}^{n}\right)\right)\right\|_2^2\right)\\
&+\sum_{n=1}^{M_p}\left(\sum_{\boldsymbol{x}_{q}^{n} \in \partial \Omega}\left\|\boldsymbol{\lambda}_{n, b}\left(\mathcal{B}\tilde{\boldsymbol{\phi}}\left(\boldsymbol{x}_{q}^{n}\right)-\boldsymbol{g}\left(\boldsymbol{x}_{q}^{n}\right)\right)\right\|_2^2\right).
\end{aligned}\label{eq:loss}
\end{equation}
 When employing the PoU function $\psi_n$ defined in \eqref{eq:psi1}, the regularization terms \added{
%  \begin{equation}
% 	\begin{aligned}
% \sum_{\substack{\Gamma \in \partial \Omega_n \cap \partial \Omega_m \\ n \neq m}} \sum_{\boldsymbol{x}_q^{\Gamma} \in \Gamma}\left(\left\|\tilde{\boldsymbol{\phi}}_n\left(\boldsymbol{x}_q^{\Gamma}\right)-\tilde{\boldsymbol{\phi}}_m\left(\boldsymbol{x}_q^{\Gamma}\right)\right\|_2^2+\left\|\nabla \tilde{\boldsymbol{\phi}}_n\left(\boldsymbol{x}_q^{\Gamma}\right)-\nabla \tilde{\boldsymbol{\phi}}_m\left(\boldsymbol{x}_q^{\Gamma}\right)\right\|_2^2\right)
% \end{aligned}\label{eq:regterms}
% \end{equation}
\begin{equation}
\begin{aligned}
\sum_{\substack{\Gamma \in \partial \Omega_n \cap \partial \Omega_m \\ n \neq m}}
\sum_{\boldsymbol{x}_q^{\Gamma} \in \Gamma}
\bigg(
&\left\|\boldsymbol{\lambda}_{n,m,\Gamma,0}
\left(\tilde{\boldsymbol{\phi}}_n\left(\boldsymbol{x}_q^{\Gamma}\right)
-\tilde{\boldsymbol{\phi}}_m\left(\boldsymbol{x}_q^{\Gamma}\right)\right)\right\|_2^2 \\
+&\left\|\boldsymbol{\lambda}_{n,m,\Gamma,1}
\left(\left(\nabla \tilde{\boldsymbol{\phi}}_n\left(\boldsymbol{x}_q^{\Gamma}\right)-\nabla \tilde{\boldsymbol{\phi}}_m\left(\boldsymbol{x}_q^{\Gamma}\right)\right)\cdot\boldsymbol{n}_{\Gamma}\right)\right\|_2^2
\bigg)
\end{aligned}\label{eq:regterms}
\end{equation}
}must be added to the loss function \eqref{eq:loss} to ensure continuity between neighboring subdomains, \added{where $\Gamma$ denotes the interface between two adjacent subdomains $\Omega_n$ and $\Omega_m \ (n\neq m)$, $\tilde{\boldsymbol{\phi}}_n$, $\tilde{\boldsymbol{\phi}}_m$ are the approximate solutions on $\Omega_n$, $\Omega_m$, $\boldsymbol{x}_q^{\Gamma}$ are collocation points sampled on $\Gamma$, $\boldsymbol{\lambda}_{n,m,\Gamma,0}$ and $\boldsymbol{\lambda}_{n,m,\Gamma,1}$ are the rescaling parameters defined in the same manner as $\boldsymbol{\lambda}_{n,q}$ and $\boldsymbol{\lambda}_{n,b}$, and $\boldsymbol{n}_{\Gamma}$ denotes the unit normal vector to $\Gamma$}. In contrast, when utilizing the $\psi_n$ defined in \eqref{eq:psi2}, the regularization terms are not required. The loss function \eqref{eq:loss} can be written in matrix form 
\begin{equation}
\mathcal{A}\boldsymbol{U}=\boldsymbol{G},
\label{eq:loss1}
\end{equation}
where $\mathcal{A}$ is the coefficient matrix related to both $\mathcal{L} \tilde{\boldsymbol{\phi}}\left(\boldsymbol{x}_{q}^{n}\right)$ and $\mathcal{B} \tilde{\boldsymbol{\phi}}\left(\boldsymbol{x}_{q}^{n}\right)$, $\boldsymbol{G}$ is the right-hand side term associated with $\boldsymbol{f}\left(\boldsymbol{x}_{q}^{n}\right)$ and $\boldsymbol{g}\left(\boldsymbol{x}_{q}^{n}\right)$.
To find the optimal parameter set $\boldsymbol{U} =\left(u_{n j}^{i}\right)^T$, the RFM samples $Q_n$ collocation points $\left\{\boldsymbol{x}_{q}^{n}\right\}_{q=1}^{Q_n}$ within each subdomain $\Omega_n$. It then calculates the rescaling parameters  $\boldsymbol{\lambda}_{n, q}=diag(\lambda_{n, q}^{1},\cdots,\lambda_{n, q}^{d_{\phi}})$ and $\boldsymbol{\lambda}_{n, b}=diag(\lambda_{n, b}^{1},\cdots,\lambda_{n, b}^{d_{\phi}})$. Specifically, the rescaling parameters are determined through the following formulas
\begin{equation}
\begin{aligned} 
& \lambda_{n, q}^{i}=\frac{c}{\mathop{\max}\limits_{1\leq j \leq J_n}\left|\mathcal{L}\left(\psi_n(\boldsymbol{x}_{q}^{n})\varphi_{n j}^{i}(\boldsymbol{x}_{q}^{n})  \right)\right|}, \\
&\quad q=1, \cdots, Q_n, \ n=1, \cdots, M_p,  \ i=1, \cdots, d_{\phi}, \\
&\lambda_{n, b}^{i}=\frac{c}{\mathop{\max}\limits_{1\leq j \leq J_n}\left|\mathcal{B}\left(\psi_n(\boldsymbol{x}_{q}^{n})\varphi_{n j}^{i}(\boldsymbol{x}_{q}^{n})  \right)\right|}, \\
&\quad \boldsymbol{x}_{q}^{n}\in{\partial\Omega}, \ n=1, \cdots, M_p, \ i=1, \cdots, d_{\phi},
\end{aligned}\label{eq:RFMlambda}
\end{equation}
where $c>0$ is a constant. Finally, the numerical result is obtained using equation \eqref{eq:approximatesolution}.  

\begin{remark}(Nonlinear PDEs)
When either the operator $\mathcal{L}$, the operator $\mathcal{B}$, or both are nonlinear, we opt to embed the least squares problem \eqref{eq:loss1} into a  nonlinear iterative solver, such as Picard's iterative method \cite{ramos2009picard}, for solving the PDE. In each iteration, the PDE is linearized, allowing the coefficient $\boldsymbol{U}$ to be updated by solving  \eqref{eq:loss1} via the linear least-squares method.
\end{remark}

\subsection{The Construction of the  Neural Feature Space}\label{sec02.2}

We follow the approach in \cite{zhang2024transferable}.
For simplicity, let's assume $\Omega_n$ to be the unit ball $B_1(0)$. For regions of other shapes and sizes, the target region can be embedded within the unit ball $B_1(0)$ via simple translations and scaling operations. The random feature function $\varphi_{n j}$ on $\Omega_n$ is redefined as follows
\begin{equation}
	\begin{aligned}
\varphi_{n j}(\boldsymbol{x}) = \sigma\left(\boldsymbol{W}_{n j}\cdot \tilde{\boldsymbol{x}}+b_{n j}\right) = \sigma\left(\gamma_{n j}\left(\boldsymbol{a}_{n j}\cdot \tilde{\boldsymbol{x}}+r_{n j}\right)\right), \quad j=1,2, \cdots, J_n,
\end{aligned}\label{eq:reparameterizedbasicfunction}
\end{equation}
where $\tilde{\boldsymbol{x}}$ represents $\boldsymbol{x}$ after the linear transformation \eqref{eq:normalize},  $\boldsymbol{a}_{n j} = \frac{\boldsymbol{W}_{n j}}{\vert\boldsymbol{W}_{n j}\vert}$, $r_{n j} = \frac{b_{n j}}{\vert\boldsymbol{W}_{n j}\vert}$ and $\gamma_{n j} = \vert\boldsymbol{W}_{n j}\vert$.   %The location parameter $\left(\boldsymbol{a}_{n j}, r_{n j}\right)$ determines 
The position of the partition hyperplane is determined by the location parameter $\left(\boldsymbol{a}_{n j}, r_{n j}\right)$ as follows
\begin{equation}
	\begin{aligned}
\boldsymbol{a}_{n j}\cdot \tilde{\boldsymbol{x}}+r_{n j}=0, \quad j=1,2, \cdots, J_n,
\end{aligned}\label{eq:partitionhyperplane}
\end{equation} 
where the unit vector $\boldsymbol{a}_{n j}$ represents the normal direction of the partition hyperplane and $\left|r_{n j}\right|$ indicates its distance from the origin. The shape parameter $\gamma_{n j}$ modulates‌ the steepness of the pre-activation value $\gamma_{n j}\left(\boldsymbol{a}_{n j}\cdot \tilde{\boldsymbol{x}}+r_{n j}\right)$ in the normal direction $\boldsymbol{a}_{n j}$.

We define the distance from a point $\boldsymbol{x}$ to the partition hyperplane \eqref{eq:partitionhyperplane} of  the random feature function $\varphi_{n j}$ on $\Omega_n$ as
\begin{equation}
	\begin{aligned}
\operatorname{dist}_j(\boldsymbol{x})=\left|\boldsymbol{a}_{n j}\cdot \tilde{\boldsymbol{x}}+r_{n j}\right|, \quad j = 1,2,\cdots,J_n, 
\end{aligned}\label{eq:dist}
\end{equation} 
and the partition hyperplane density as
\begin{equation}
	\begin{aligned}
D_{J_n}^{\tau}(\boldsymbol{x})=\frac{1}{J_n} \sum_{j=1}^{J_n} 1_{\operatorname{dist}_j(\boldsymbol{x})<\tau}(\boldsymbol{x}),
\end{aligned}\label{eq:density}
\end{equation} 
where $\tau>0$ denotes the bandwidth for density estimation, and 
$1_{\operatorname{dist}_j(\boldsymbol{x})<\tau}(\boldsymbol{x})$ represents the indicator function ‌evaluating whether the distance between $\boldsymbol{x}$ and the partition hyperplane of  the random feature function $\varphi_{n j}$ on $\Omega_n$  satisfies $\operatorname{dist}_j(\boldsymbol{x}) < \tau$.
 The $D_{J_n}^{\tau}(\boldsymbol{x})$ quantifies the proportion of the random feature functions whose partition hyperplanes intersect the ball $B_\tau(\boldsymbol{x})$, centered at $\boldsymbol{x}$  with radius $\tau$. 
% Theoretical proofs show that the neurons generated in this way are uniformly distributed within $\Omega$. .

In solving PDEs, the distribution‌ of partition hyperplanes \eqref{eq:partitionhyperplane} within $\Omega_n$ can be controlled. For example, in  \cite{zhang2024transferable}, uniform distribution of partition hyperplanes is preferred.  The location parameters $\left(\boldsymbol{a}_{n j}, r_{n j}\right)$ are randomly chosen as
%it is necessary to guarantee that the partition hyperplanes \eqref{eq:partitionhyperplane} are uniformly distributed within
%$\Omega_n$, the location parameter $\left(\boldsymbol{a}_{n j}, r_{n j}\right)$ is randomly chosen as
\begin{equation}
	\begin{aligned}
\boldsymbol{a}_{n j}=\frac{\boldsymbol{X}_{n j}}{\vert\boldsymbol{X}_{n j}\vert} \  \text{and}  \ r_{n j}=D_{n j}, \quad j=1,2, \cdots, J_n,
\end{aligned}\label{eq:chooseparameter}
\end{equation}  
where $\boldsymbol{X}_{n j}$ is distributed as a $d$-dimensional standard Gaussian distribution, and $D_{n j}$ follows a uniform distribution over $[0,1]$. This sampling method generates a set of uniformly distributed partition hyperplanes \eqref{eq:partitionhyperplane} in $\Omega_n$. For a fixed $\tau \in(0,1)$, the expectation has $\mathbb{E}\left[D_{J_n}^{\tau}(\boldsymbol{x})\right]=\tau$ when $\vert\tilde{\boldsymbol{x}}\vert \leq 1-\tau$.
%when $\|\boldsymbol{x}\|_2 \leq 1-\tau$, there is $\mathbb{E}\left[D_{J_n}^{\tau}(\boldsymbol{x})\right]=\tau$.
The  same shape parameter $\gamma_{nj} = \gamma_n$ are also adopted  for all random feature functions $\varphi_{n j}$ defined on $\Omega_n$  . To compute the optimal shape parameter $\gamma_n$,  $L$ realizations of the Gaussian random fields (GRFs) $G\left(\boldsymbol{x} \mid \omega_l, \eta\right)_{l=1}^L$ are simulated, where $\omega_l$ denotes abstract randomness and $\eta$ represents a fixed correlation length. The fitting error for each realization of GRFs and the approximate solution on $\Omega_n$ is defined as 
\begin{equation}
\begin{aligned}
\operatorname{Loss}^{\gamma_n}_{l}=\min\limits_{\left\{u_{n j}\right\}}\left(\sum_{q=1}^{Q_n}\left\| \sum_{j=1}^{J_n} u_{n j}\sigma\left(\gamma_n\left(\boldsymbol{a}_{n j}\cdot \tilde{\boldsymbol{x}}_{q}^{n}+r_{n j}\right)\right)-G\left(\boldsymbol{x}_{q}^{n} \mid \omega_l, \eta\right)\right\|_2^2\right).
\end{aligned}\label{eq:lossgamma}
\end{equation}
The optimal shape parameter $\gamma_n = \gamma_n^{opt}$ is obtained by 
\begin{equation}
\begin{aligned}
\gamma_n^{o p t}=\arg \min _{\gamma_n} \frac{1}{L} \sum_{l=1}^L\operatorname{Loss}^{\gamma_n}_{l},
\end{aligned}\label{eq:optgamma}
\end{equation}
which can be solved via grid search.

\section{The Adaptive Feature Capture Method}

We now present the Adaptive Feature Capture Method (AFCM). For the boundary value problem in Equation (1), we first partition the domain $\Omega$ into $M_p$ non-overlapping subdomains $\{\Omega_n\}_{n=1}^{M_p}$. Within each $\Omega_n$, we sample $Q_n$ collocation points $\{\bm{x}_q^n\}_{q=1}^{Q_n}$ and construct $J_n$ feature functions $\{\varphi_{nj}\}_{j=1}^{J_n}$ using the methodology described in Section \ref{sec02.2}. We employ the partition of unity (PoU) function $\psi_n$ defined in Equation \eqref{eq:psi1}. The total number of feature functions is denoted as 
$ J=\sum_{n=1}^{M_p} J_n$. 
To ensure continuity between subdomains, regularization terms enforcing $C^1$ continuity conditions at interfaces are added to the loss function \eqref{eq:loss} following the approach in \cite{dong2021local}. The loss matrix \eqref{eq:loss1} is then assembled and solved via linear least squares to obtain an initial approximate solution $\tilde{\bm{\phi}}(\bm{x})$.

The accuracy of this initial solution $\tilde{\bm{\phi}}(\bm{x})$ depends on the solution's regularity. While satisfactory for smooth solutions, its performance degrades significantly for near-singular problems with steep gradients. This limitation arises because uniformly distributed partition hyperplanes and collocation points in the Random Feature Method (RFM) lack sufficient local expressive power in high-gradient regions. To address this, the AFCM adaptively increases feature function density and collocation point concentration in critical areas while preserving computational efficiency. The adaptation process consists of two key components:

\textbf{Step 1:  Feature Function Adaptation} 

We enhance the density of partition hyperplanes and the steepness of feature function pre-activations in high-gradient regions by adjusting the shape parameter $\gamma_{nj}$ and location parameters $(\bm{a}_{nj}, r_{nj})$, using the gradient norm $|\nabla\tilde{\bm{\phi}}(\bm{x})|$ as an indicator. %We first uniformly sample $m$ points $S$ across $\Omega$ and define a probability density function (PDF) on $S$:
We first uniformly sample a set $S$ of $m$ points across $\Omega$
and define a probability density function (PDF) on $S$:
\begin{equation}\label{eq:PDF}
    p(\bm{x}) = \frac{(|\nabla\tilde{\bm{\phi}}(\bm{x})| + c_1)}{\sum_{\bm{x}\in S}(|\nabla\tilde{\bm{\phi}}(\bm{x})| + c_1)}, \quad \bm{x} \in S,
\end{equation}
where $c_1 > 0$ prevents excessive concentration. \added{It is required that $m\gg\sum_{n=1}^{M_p}J_n=J$ and $m\gg \sum_{n=1}^{M_p}Q_n$  to ensure adequate coverage and accurate estimation of the solution gradient distribution.} From this PDF, we sample $J$  points via weighted random sampling (WRS) without replacement from the set $S$. \added{This sampling scheme draws a fixed number of distinct points from the finite set $S$ according to the PDF without repetition: each point $\bm{x} \in S$ is selected with probability proportional to the PDF, and once chosen, it will not be sampled again in the same iteration.} With $\{\bm{x}_j^n\}_{j=1}^{J'_n}$ denoting points in $\Omega_n$ (note that $J'_n$ is usually  different from $J_n$), \deleted{where the weights for sampling are determined by the PDF.} the feature functions on $\Omega_n$ are reconstructed via:
\begin{align}
    \gamma_{nj} &= \deleted{\gamma_n} \added{\beta_n} \cdot \frac{(|\nabla\tilde{\bm{\phi}}(\bm{x}_j^n)| + c_2)}{\min_{\bm{x}\in\{\bm{x}_j^n\}_{j=1}^{J'_n}}(|\nabla\tilde{\bm{\phi}}(\bm{x})| + c_2)}, \quad j=1,\ldots,J'_n, \label{eq:regenerategamma}\\
    \bm{a}_{nj} &= \frac{\bm{X}_{nj}}{|\bm{X}_{nj}|}, \quad r_{nj} = -\bm{a}_{nj} \cdot \tilde{\bm{x}}_j^n, \quad j=1,\ldots,J'_n, \label{eq:regenerateb}
\end{align}
where \deleted{$\gamma_n$ is the initial shape parameter from Equations \eqref{eq:lossgamma}--\eqref{eq:optgamma},} \added{$\beta_n=\gamma_n$ denotes the fixed initial baseline shape parameter of subdomain $\Omega_n$, computed once in the initial iteration via Equations \eqref{eq:lossgamma}--\eqref{eq:optgamma}, and remains constant throughout all subsequent adaptation iterations.} $\bm{X}_{nj}$ follows a $d$-dimensional standard Gaussian distribution, and $c_2 > 0$ controls shape parameter variation to avoid numerical instability.

During the adaptation process, Equation \eqref{eq:regenerategamma}  determines the shape parameter $\gamma_{nj}$ at each sampling point $\bm{x}_j^n$ by scaling the \deleted{initial constant $\gamma_n$} \added{fixed baseline $\beta_n$} proportionally to the local gradient magnitude $|\nabla\tilde{\bm{\phi}}(\bm{x}_j^n)|$. \added{Notably, $\gamma_{nj}$ is always adjusted from the invariant initial baseline $\beta_n$ at each adaptation step, rather than being updated recursively from the $\gamma_{nj}$ obtained in the previous iteration.} This proportional scaling ensures feature functions become steeper in high-gradient regions while maintaining moderate variations elsewhere. Specifically, locations with larger solution gradients receive amplified $\gamma_{nj}$ values, enhancing local approximation capability without global parameter modifications.

Equation \eqref{eq:regenerateb}  geometrically enforces that each feature function's partition hyperplane intersects its corresponding sampling point $\bm{x}_j^n$. This critical constraint ensures the adapted basis functions remain anchored to regions requiring enhanced resolution, creating a dynamic alignment between network architecture and solution characteristics.

% While the adaptation is taking place, the equation \eqref{eq:regenerategamma} determines the shape parameter $\gamma_{nj}$ at each sampling point $\boldsymbol{x}_{j}^{n}$ by multiplying the initial shape parameter $\gamma_n$ ($\gamma_n$ is a constant) by a proportional factor associated with the approximate solution gradient $\nabla\tilde{\boldsymbol{\phi}}({\boldsymbol{x}_{j}^{n}})$.  Intuitively, for points where the magnitude of the gradient is relatively large, the corresponding  $\gamma_{nj}$ is relatively large. 
 % and it is always the case that $\gamma_{nj}\geq \gamma_n$. 
 %Meanwhile, the equation \eqref{eq:regenerateb} ensures that the partition hyperplane of the feature function passes through the sampling point $\boldsymbol{x}_{j}^{n}$.

\textbf{Step 2:  Collocation Point Adaptation} 

Interior collocation points are redistributed to high-gradient regions using the same PDF $p(\bm{x})$. We sample $\sum_{n=1}^{M_p} I_n$ points via WRS without replacement from the set $S$, where $I_n$ denotes the number of initial interior collocation points in 
$\Omega_n$.  The set $\{\bm{x}_i^n\}_{i=1}^{I'_n}$  represents the interior points in $\Omega_n$, and the weights for sampling are determined by the PDF. The updated collocation points on $\Omega_n$ become:
\begin{equation}\label{eq:regeneratepoint}
    \left\{\bm{x}_q^n \mid \bm{x}_q^n \in \partial\Omega_n\right\} \cup \left\{\bm{x}_i^n\right\}_{i=1}^{I'_n}.
\end{equation}

\begin{figure}[htp]
	\center{\includegraphics[scale=0.17]  {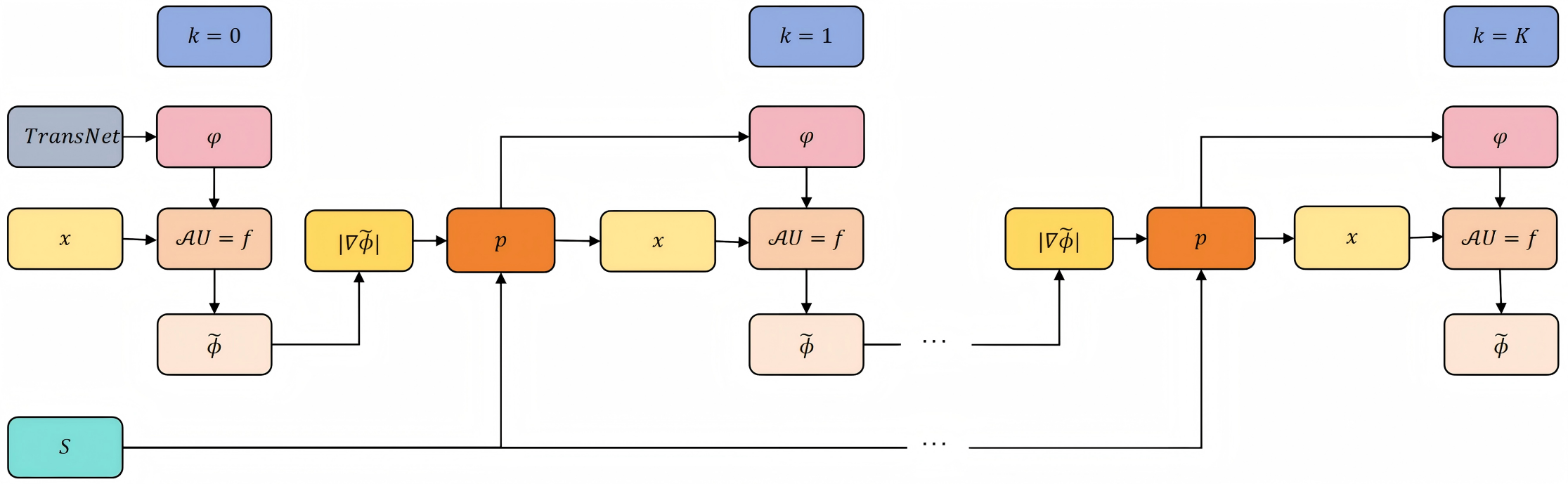}} 
	\caption{The computational framework of of Algorithm \ref{alg:Algorithm 1}.}
	\label{fig:archi}
\end{figure}

After these adaptations, the loss matrix \eqref{eq:loss1} is reassembled and solved to obtain an improved approximation $\tilde{\bm{\phi}}(\bm{x})$. The AFCM iteratively applies this process until convergence, as formalized in Algorithm \ref{alg:Algorithm 1}. Figure \ref{fig:archi} illustrates the computational framework.

\begin{remark} \added{(Reuse of Previous‑round Solution)}
\added{In the proposed AFCM, the linear system is reconstructed and resolved from scratch following adaptive updates of collocation points and feature functions. The prior approximate solution guides the adaptive strategy rather than serving as an initial guess: its gradient norm constructs a probability density that regulates weighted sampling, clustering key components in high-gradient regions.}

\added{This ensures stable convergence: high-gradient regions where resolution is insufficient are targeted to enhance local expressiveness. Numerical experiments demonstrate monotonic reduction of relative error, with up to 10 orders of magnitude decrease for near-singular solutions, validating the adaptive mechanism.}

\added{Computationally, extra overhead such as gradient computation and sampling is negligible. With fixed degrees of freedom and collocation points (only their positions adjusted), each linear solve incurs a cost similar to the original RFM. Merely 4–5 iterations attain near-machine precision, leading to a slight increase in total cost, while avoiding the high mesh-related costs of traditional methods and retaining the mesh-free efficiency of RFM.}
\end{remark}

\begin{algorithm}
	\caption{ The Adaptive Feature Capture Method}
	\label{alg:Algorithm 1}
	\begin{algorithmic}[0]
		\Require
        ~$\bullet$\hspace{0.5em}Number of subdomains $M_p$
        \begin{itemize}[leftmargin=3.96em, labelsep=0.5em]
            % \item Number of subdomains $M_p$
            \item Number of feature functions per subdomain $J_n$
            \item Number of collocation points per subdomain $Q_n$
            \item Number of interior collocation points per subdomain $I_n$
            \item Rescaling parameter $c$, coefficients $c_1$, $c_2$
            \item Number of iterations $K$
            \item Number of sampling points $m$
            \item Number of GRF realizations $L$
            \item Correlation length $\eta$
        \end{itemize}
		\Ensure
		The $K$-th approximate solution $\tilde{\boldsymbol{\phi}}(\boldsymbol{x})$;
		\State Divide $\Omega$ into $M_p$ non-overlapping subdomains ${\Omega_n}$;
        \For{$k = 0$}
            \State \textbf{1.}  Sample $Q_n$ collocation points $\left\{ \boldsymbol{x}_{q}^{n}\right\}_{q=1}^{Q_n}$  in each $\Omega_n$;
            \State \textbf{2.} Generate $J_n $ location parameters $\{\left(\boldsymbol{a}_{n j}, r_{n j}\right)\}_{j=1}^{J_{n}}$ on $\Omega_n$ according to \eqref{eq:chooseparameter};
            \State \textbf{3.} Compute the shape parameter $\gamma_n$ on $\Omega_n$ according to \eqref{eq:lossgamma} and \eqref{eq:optgamma};
            \State \textbf{4.} Construct $J_n$ feature functions $\{\varphi_{n j}\}_{j=1}^{J_n}$ on $\Omega_n$ by \eqref{eq:reparameterizedbasicfunction};
            \State \textbf{5.}  Assemble  the loss matrix \eqref{eq:loss1} and solve it by the linear least-squares method to obtain the $k$-th approximate solution $\tilde{\boldsymbol{\phi}}(\boldsymbol{x})$;
        \EndFor
        \State  Sample $m$ points $S$ in the entire $\Omega$  by the uniform distribution;
	\For{$k = 1,2, \cdots, K$}
            \State  \textbf{1.} Define PDF $p(\boldsymbol{x})$ for the points $S$  according to \eqref{eq:PDF};
            \State  \textbf{2.} Sample $\sum_{n=1}^{M_p}J_n$ points from $S$ according to $p(\boldsymbol{x})$ and reconstruct the feature functions $\{\varphi_{n j}\}_{j=1}^{J'_n}$ on $\Omega_n$ by \eqref{eq:regenerategamma} and \eqref{eq:regenerateb}
            \State  \textbf{3.} Sample $\sum_{n=1}^{M_p}I_n$ points from $S$ according to $p(\boldsymbol{x})$ and regenerate the collocation points on $\Omega_n$ according to \eqref{eq:regeneratepoint}.
            \State   \textbf{4.} Assemble  the loss matrix \eqref{eq:loss1} and solve it by the linear least-squares method to obtain the $k$-th approximate solution $\tilde{\boldsymbol{\phi}}(\boldsymbol{x})$;
        \EndFor
        \State Return the $K$-th approximate solution $\tilde{\boldsymbol{\phi}}(\boldsymbol{x})$.
	\end{algorithmic}

\end{algorithm}

This adaptive strategy ensures that partition hyperplanes and collocation points cluster in critical regions, achieving enhanced resolution without increasing computational overhead. The mesh-free nature of RFM is preserved, making AFCM suitable for complex geometries while maintaining efficiency.

\section{Numerical Experiments}
\label{sec04}

\begin{figure}[htp]
	\center
  {\subfigure[The arrangement of subdomains. ]{\includegraphics[scale=0.14]{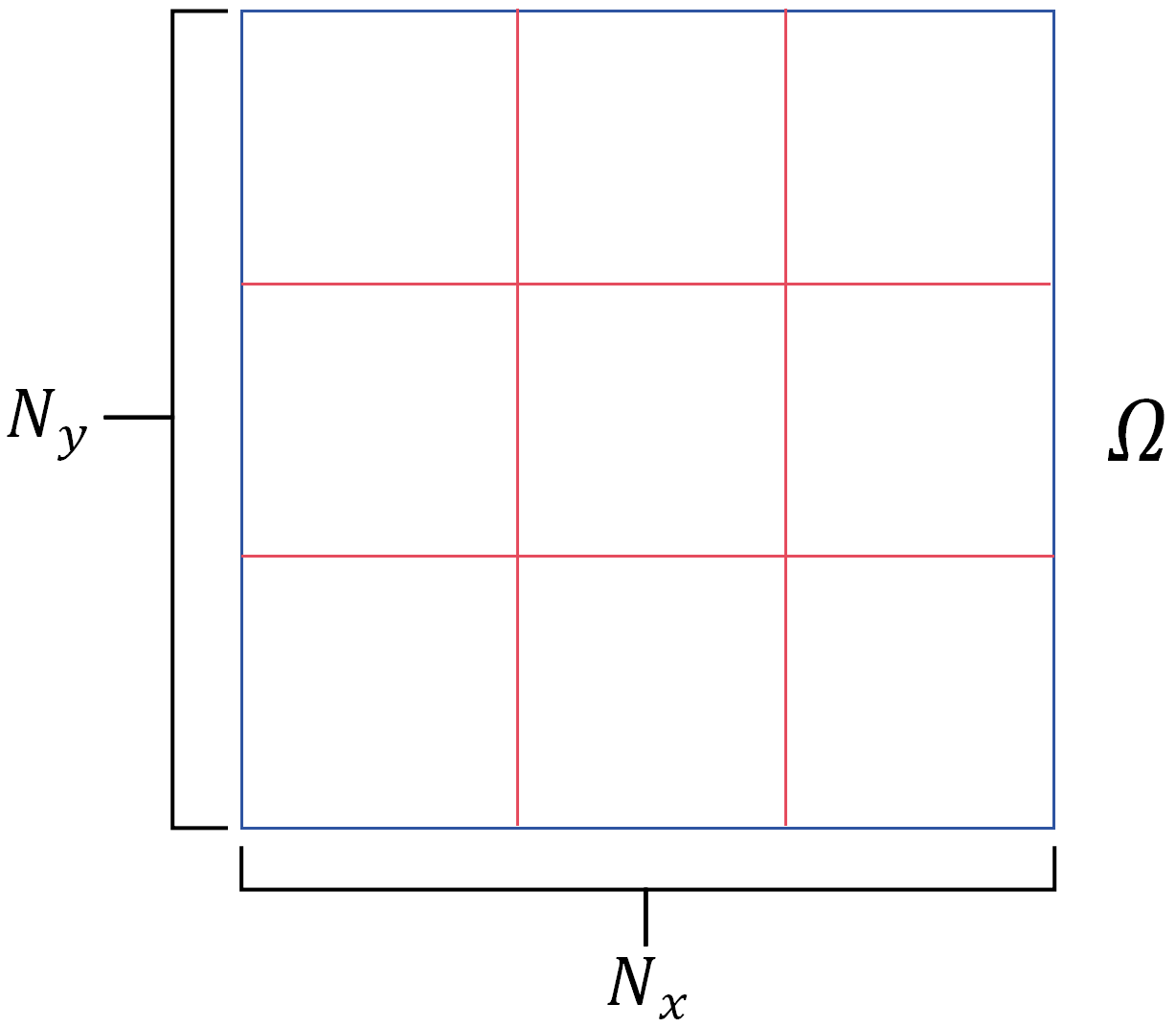}}
   \subfigure[The distribution of collocation points for $k=0$ iteration in subdomain $\Omega_{n}$. ]{
            \includegraphics[scale=0.15]{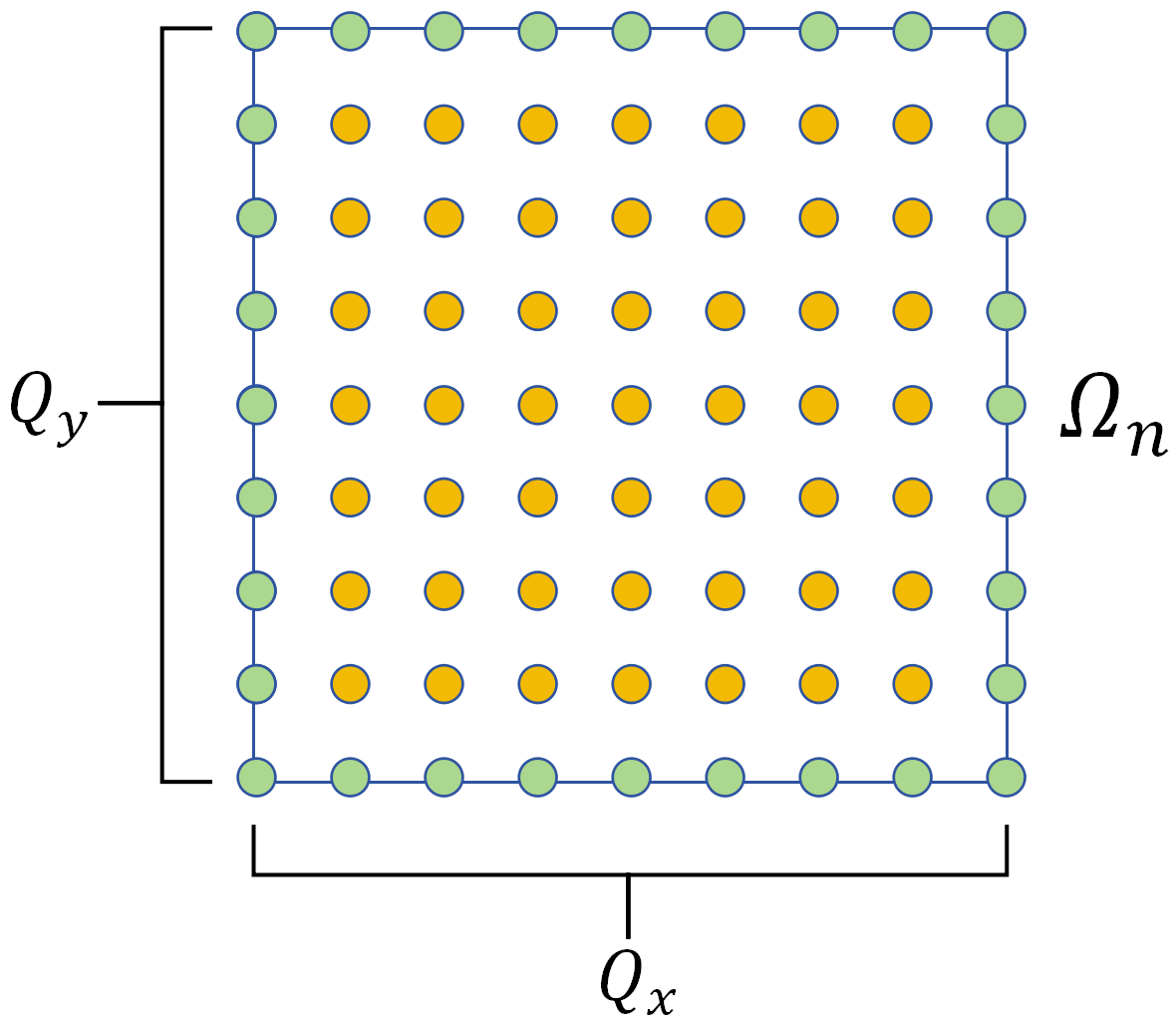}}}
	\caption{The arrangement of subdomains and the distribution of collocation points for $k=0$ iteration in subdomain $\Omega_{n}$. (a) The blue lines and red lines represent the boundaries and the interfaces between subdomains, respectively. (b) The green points and yellow points represent the boundary points and interior points in each subdomain, respectively.}
	\label{fig:MQ}
\end{figure}

We present numerical experiments to validate the proposed method. The relative $L^{\infty}$ and $L^{2}$ error norms between the numerical solution $\tilde{\boldsymbol{\phi}}_k(\boldsymbol{x})$ and exact solution $\boldsymbol{\phi}(\boldsymbol{x})$ are defined as:
\begin{equation}
    \|e_k\|_{L^{\infty}} = \frac{\|\tilde{\boldsymbol{\phi}}_k(\boldsymbol{x}) - \boldsymbol{\phi}(\boldsymbol{x})\|_{L^{\infty}}}{\|\boldsymbol{\phi}(\boldsymbol{x})\|_{L^{\infty}}}, \quad
    \|e_k\|_{L^{2}} = \frac{\|\tilde{\boldsymbol{\phi}}_k(\boldsymbol{x}) - \boldsymbol{\phi}(\boldsymbol{x})\|_{L^{2}}}{\|\boldsymbol{\phi}(\boldsymbol{x})\|_{L^{2}}}
\end{equation}
where $k=0,\ldots,K$ denotes the number of adaptive iterations. For $d=2$, we set $M_p = N_x \times N_y$ and $Q_n = Q_x \times Q_y$, where $N_x,N_y$ represent subdomain counts and $Q_x,Q_y$ collocation points per subdomain. Figure~\ref{fig:MQ} shows the initial subdomain arrangement and collocation point distribution.

All linear least squares computations use PyTorch's \texttt{torch.linalg.lstsq} solver. Key experimental parameters are summarized in Table~\ref{Table:set}

% \begin{table}[htp]
% 	\begin{center} 
%     \caption{
% Parameters and experimental configurations for numerical experiments }\label{Table:set}
% \vspace{10pt}
% \small
% \begin{tabular}{lc}
% \hline 
% \text {Number of Subdomains $M_p=N_x\times N_y$} & $3\times 3$ \\
% \text {Number of Collocation points per subdomain $Q_n=Q_x\times Q_y$} &  $79\times 79$\\
% \text {Number of Interior collocation points per subdomain $I_n$} &  $77\times 77$  \\
% % \text {\small  number of test points in each subdomain $Q_{test}$ } & $158\times 158$   \\
% \text {Constant of Rescaling parameter $c$}  & 1 \\
% \text {Number of GRFs realizations $L$} & 10  \\
% \text {Constant $c_1$} & 0.01 \\
% \text {Constant $c_2$} & 50 \\
% \text {Density bandwidth $\tau$} & 0.2 \\
% \text {Correlation length $\eta$}  & \text { 0.5 } \\
% \text {Activation function $\sigma$} &  $tanh^3$ \\
% \text {CPU} & \text {\small Intel Xeon Platinum 8358} \\
% \text {GPU} & \text {\small NVIDIA A100 (80GB)} \\
% \hline
% \end{tabular}
% 	\end{center}
% \end{table}

\begin{table}[htp]
    \centering
    \caption{Parameters and experimental configurations for numerical experiments}
    \label{Table:set}
    \vspace{10pt}
    \begin{tabular}{lc}
        \hline 
        \makecell[l]{Number of Subdomains \\ $M_p=N_x\times N_y$} & $3\times 3$ \\
        \makecell[l]{Number of Collocation points per subdomain \\ $Q_n=Q_x\times Q_y$} & $79\times 79$\\
        \makecell[l]{Number of Interior Collocation Points \\ per Subdomain $I_n$} & $77\times 77$  \\
        Constant of Rescaling Parameter $c$ & 1 \\
        Number of GRFs Realizations $L$ & 10  \\
        Constant $c_1$ & 0.01 \\
        Constant $c_2$ & 50 \\
        Density Bandwidth $\tau$ & 0.2 \\
        Correlation Length $\eta$ & 0.5 \\
        Activation Function $\sigma$ & $\tanh^3$ \\
        CPU & Intel Xeon Platinum 8358 \\
        GPU & NVIDIA A100 (80GB) \\
        \hline
    \end{tabular}
\end{table}

\subsection{Two-Dimensional Poisson equation with a near singular solution with two peaks}
We consider the Poisson equation on $\Omega = (-1,1)^2$:
\begin{equation} \label{eq:poissoneq}
    \begin{cases}
        -\Delta\boldsymbol{\phi}(x,y) = f(x,y), & (x,y) \in \Omega \\
        \boldsymbol{\phi}(x,y) = g(x,y), & (x,y) \in \partial\Omega
    \end{cases},
\end{equation}
the exact solution with two peaks is given by
\begin{equation}\label{eq:exactsolution1}
\begin{aligned}
\boldsymbol{\phi}(x,y)=e^{-1000\left(x^2+(y-\frac{2}{3})^2\right)}+e^{-1000\left(x^2+(y+\frac{2}{3})^2\right)}.
\end{aligned}
\end{equation}
The source term $f(x,y)$ and boundary condition $g(x,y)$ are derived directly from \eqref{eq:poissoneq}--\eqref{eq:exactsolution1}. 
% Using $m = 1.26 \times 10^6$ sampling points, we obtain adaptive solutions via Algorithm~\ref{alg:Algorithm 1}. 
We implement Algorithm \ref{alg:Algorithm 1} with $m = 1.26 \times 10^6$ sampling points to resolve the steep gradients near $(0, \pm\frac{2}{3})$.
% \begin{table}[htp]
% 	\begin{center} 
%  \caption{\deleted{The initial shape parameter $\gamma_n$ and the relative errors for \eqref{eq:exactsolution1}  with different $J_n$ at  $k=0$ and $k=K=4$ iterations}\added{The initial shape parameter $\gamma_n$, the relative errors at $k=0$ and $k=K=4$ iterations, the computational time, and the peak GPU memory consumption  for \eqref{eq:exactsolution1} with different $J_n$}}
%  \vspace{10pt}
% 		\begin{tabular}{|l|l|l|l|l|l|l|l|l|l|}
% 			\hline      
% 		\multirow{2}{*}{$m$}& $J_n$ & 1500 & 2000  & 3000 & 4000         \\
%         \cline{2-6}
%         &  $\gamma_n$ & 2.0 & 2.6  & 2.8 & 3.4       \\
% 			\hline
% 	\multirow{6}{*}{$1.26\times 10^6$}& $\vert\vert e_{0}\vert\vert_{L^{\infty}} $  & 9.79E-2 & 2.55E-1 & 8.26E-2 & 9.22E-2     \\
%     \cline{2-6}
% 		&  $\vert\vert e_{0}\vert\vert_{L^2} $  & 2.50E-1 & 1.02E-0 & 2.70E-1 & 1.04E-0   \\
%             \cline{2-6}
% 		&  $\vert\vert e_{K=4}\vert\vert_{L^\infty} $  & 2.48E-5 & 1.19E-8 & 1.02E-9 & 6.35E-11   \\
%                 \cline{2-6}
% 		&  $\vert\vert e_{K=4}\vert\vert_{L^2} $  & 3.83E-5 & 2.11E-8 & 2.99E-9 & 1.38E-10   \\

% 			\cline{2-6}
%         &  \added{time (s)}  & \added{56.2} & \added{73.7} & \added{134} & \added{206.6}   \\
%                 \cline{2-6}
% 		&  \added{Memory (GB)}  & \added{75.0} & \added{72.8} & \added{71.5} & \added{75.5}   \\

% 			\hline
% 		\end{tabular}
% 		\label{Table:1}
% 	\end{center}
% \end{table}
\begin{table}[htp]
	\begin{center} 
 \caption{ \deleted{The initial shape parameter $\gamma_n$ and the relative errors for \eqref{eq:exactsolution1}  with different $J_n$ at  $k=0$ and $k=K=4$ iterations}\added{The initial shape parameter $\gamma_n$, the relative errors at $k=0$ and $k=K=4$ iterations,  the computational time, and the peak GPU memory consumption  for \eqref{eq:exactsolution1} with different $J_n$}}
  \vspace{10pt}
		\begin{tabular}{|l|l|l|l|l|l|l|l|l|l|}
			\hline
            
		\multirow{2}{*}{$m$}& $J_n$ & 1500 & 2000  & 3000 & 4000         \\
        \cline{2-6}
        &   $\gamma_n$ & 2.0 & 2.6  & 2.8 & 3.4       \\
			\hline
	\multirow{6}{*}{$1.26\times 10^6$}& $\vert\vert e_{0}\vert\vert_{L^{\infty}} $  & 1.18E-1 & 1.08E-1 & 5.93E-2 & 1.82E-1     \\
    \cline{2-6}
		&  $\vert\vert e_{0}\vert\vert_{L^2} $  & 1.64E-1 & 2.48E-1 & 1.31E-1 & 4.65E-1   \\
            \cline{2-6}
		&  $\vert\vert e_{K=4}\vert\vert_{L^\infty} $  & 4.10E-5 & 9.79E-9 & 1.75E-9 & 1.07E-10   \\
                \cline{2-6}
		&  $\vert\vert e_{K=4}\vert\vert_{L^2} $  & 6.75E-5 & 1.28E-8 & 2.02E-9 & 9.76E-11   \\
				\cline{2-6}
        &  \added{time (s)}  & \added{56.1} & \added{73.9} & \added{134.8} & \added{205.3}   \\
                \cline{2-6}
		& \added{Memory (GB)} & \added{74.6} & \added{46.1} & \added{75.4} & \added{64.7}  \\

			\hline
		\end{tabular}
		\label{Table:1}
	\end{center}
\end{table}

\begin{figure}[htp]
	\center
        {\subfigure[exact solution $\boldsymbol{\phi}$ ]{		 
            \includegraphics[scale=0.23]{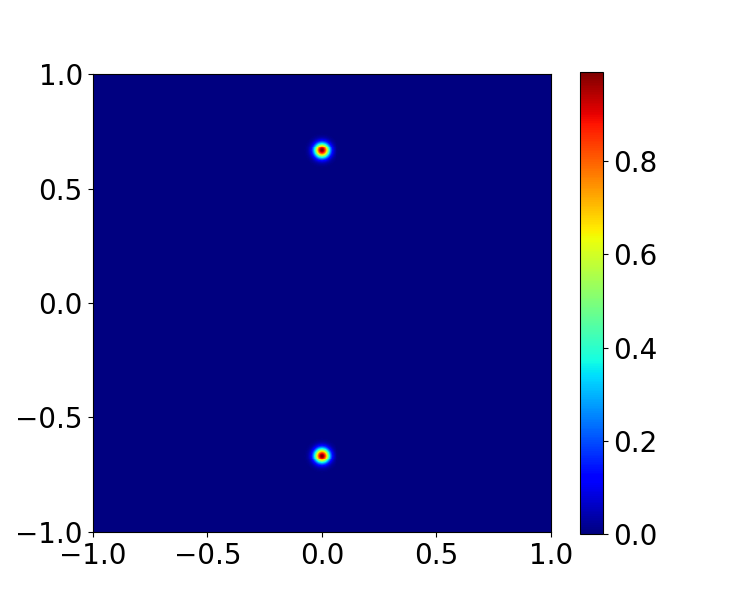}\label{3a}}
        \subfigure[approximte solution $\tilde{\boldsymbol{\phi}}_{0}$]{		 
            \includegraphics[scale=0.23]{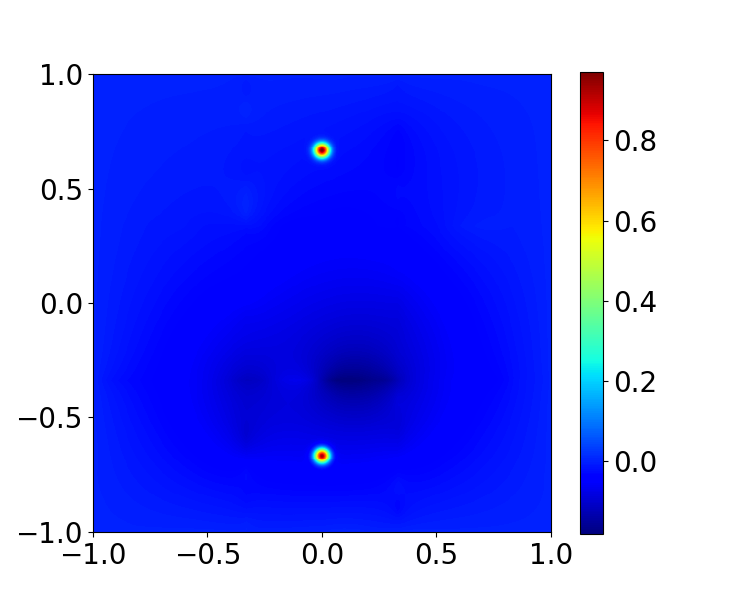}\label{3b}}
            \subfigure[local surface of  $\tilde{\phi}_0$ ]{		 
            \includegraphics[scale=0.32]{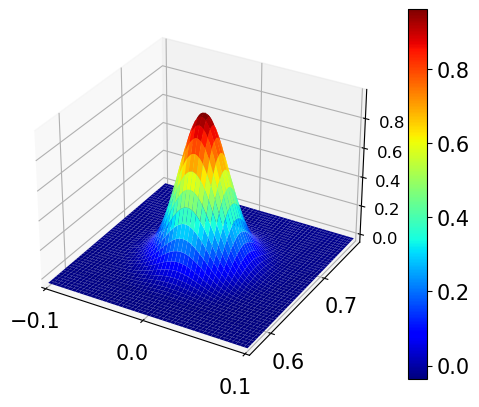}}
            \subfigure[Relative errors and  training losses]{		
            \includegraphics[scale=0.195]{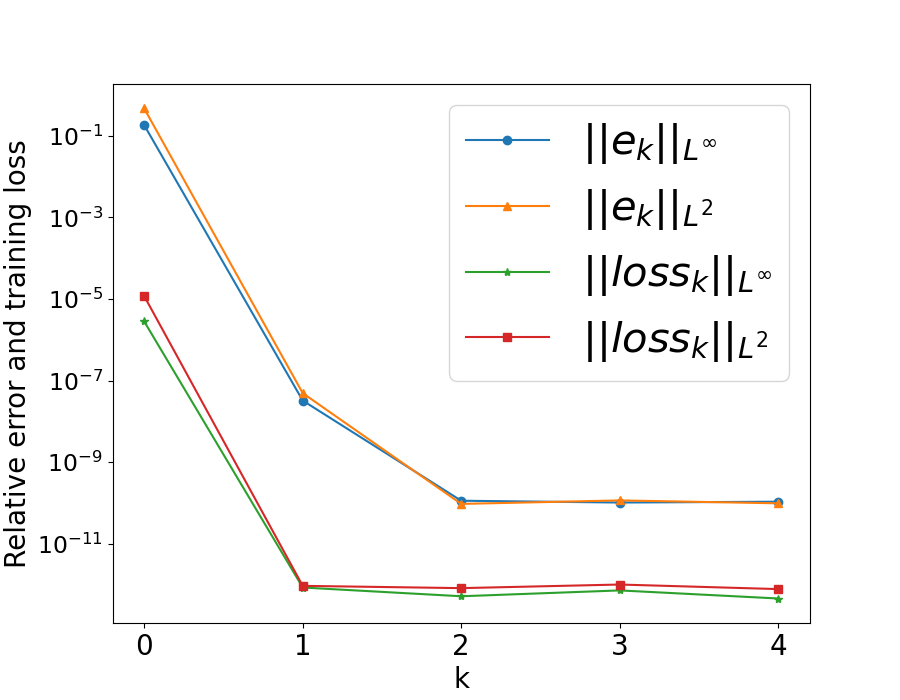}\label{3d}}
            \subfigure[approximte solution $\tilde{\boldsymbol{\phi}}_{K} $  ]{		
            \includegraphics[scale=0.23]{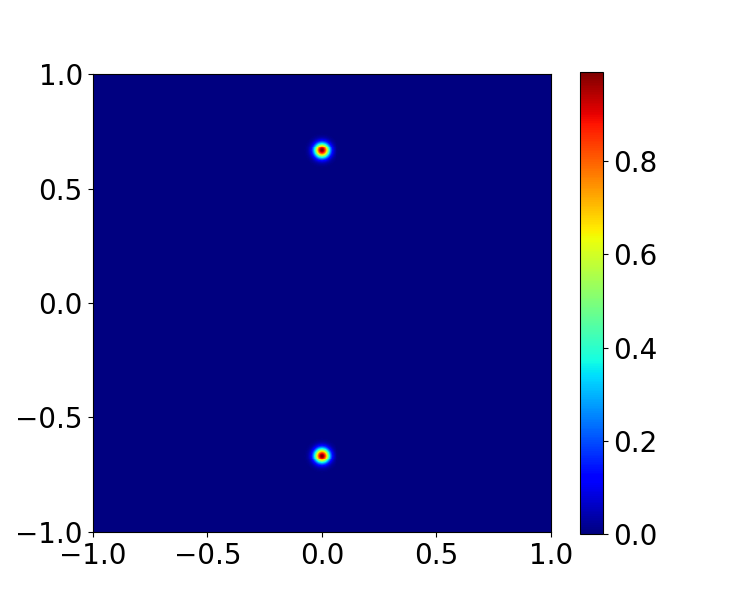}}
            \subfigure[local surface of  $\tilde{\phi}_K$ ]{		 
            \includegraphics[scale=0.32]{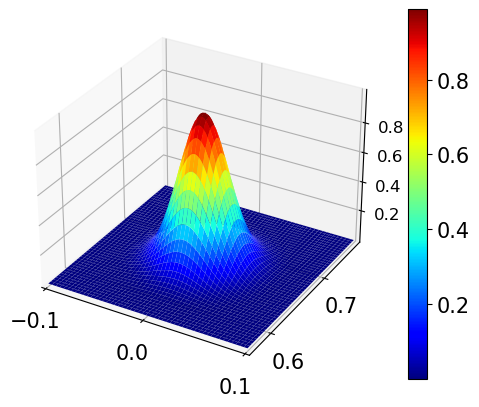}\label{3f}}
            }
	\caption{The exact solution and the numerical results for \eqref{eq:exactsolution1}  with $J_n=4000$ and $K=4$. }
	\label{fig:fig1}
\end{figure}

We show in Table \ref{Table:1} the relative errors  for $\boldsymbol{\phi}(x,y)$  at the initial  iteration ($k = 0$) and the  final iteration  ($k = 4$), \added{ as well as the computational time and peak GPU memory consumption of the full adaptive procedure,} for different  $J_n$.  The results demonstrates significant error reduction through adaptive iterations, with monotonic decay as $J_n$ increases. 
\added{We compare the performance of AFCM with the adaptive finite element methods (AFEM) \cite{babuvvska1978error} using P1 and P2 elements, where AFEM adopts the standard residual-based error indicator for adaptivity on triangular meshes. The AFEM implementation is based on FEniCSx; the initial meshes are generated with $189\times189$ and $94\times94$ subdivisions for AFEM-P1 and AFEM-P2, respectively; and the cells with the largest indicators are marked with a marking fraction of 0.25 for refinement. Both methods are initialized with a nearly identical initial number of degrees of freedom (DoF). The adaptive iterations terminate when both the relative $L^\infty$ and $L^2$ errors are below the prescribed tolerance of 1.0E-8. As shown in Table \ref{Table:AFEM1}, AFCM achieves substantially higher accuracy with a fixed and much smaller DoF of 36000. In contrast, AFEM-P1 fails to attain the prescribed error bound even after intensive adaptive refinement, while AFEM-P2 requires 8 adaptive iterations with significantly increased DoF to meet the stopping criterion. AFCM converges to the target tolerance in only 2 iterations and runs faster with a runtime of merely 95.6s. These results demonstrate that AFCM resolves steep gradients and near-singular solutions with better accuracy and computational efficiency.}
\begin{table}[htp]
	\begin{center} 
 \caption{\added{Comparison of AFCM with the traditional AFEM with P1 and P2 elements. We show the degrees of freedom (at final iteration), the relative errors, the iteration numbers $K$ and the computational time at $k = K$ iterations of AFEM and AFCM  with $J_n = 4000$}}
 \vspace{10pt}
\begin{tabular}{|l|l|l|l|l|l|}
			\hline      
		\added{Method} &  \added{DoF ($k=K$)} & \added{$\vert\vert e_{k=K}\vert\vert_{L^\infty} $} & \added{$\vert\vert e_{k=K}\vert\vert_{L^2} $} & \added{K} & \added{time (s)} \\
            \hline
            \added{AFEM-P1} & \added{6956242} & \added{1.15E-6} & \added{1.10E-6} & \added{9} & \added{1781.4}\\
			\hline
        \added{AFEM-P2} & \added{4409761} & \added{8.44E-9} & \added{3.78E-9} & \added{8} & \added{935.1}\\
            \hline
            % \added{AFEM-P3} & \added{1185142} & \added{2.89E-9} & \added{6.31E-10} & \added{6} &\added{114.2}\\
            % \hline
        \added{AFCM} & \added{36000} & \added{5.74E-11} & \added{1.28E-10} & \added{2} &\added{95.6}\\
            \hline
		\end{tabular}
		\label{Table:AFEM1}
	\end{center}
\end{table}
% \begin{table}[htp]
% 	\begin{center} 
%  \caption{\added{The degrees of freedom, the relative errors and the computational time at $k=K=4$ iterations of AFEM and AFCM for \eqref{eq:exactsolution1} with $J_n=4000$}}
%  \vspace{10pt}
% \begin{tabular}{|l|l|l|l|l|l|l|l|l|l|}
% 			\hline      
% 		\added{Method} &  \added{DoF (K=4)} & \added{$\vert\vert e_{K=4}\vert\vert_{L^\infty} $} & \added{$\vert\vert e_{K=4}\vert\vert_{L^2} $} & \added{time (s)} \\
%         \hline
%         \added{AFEM-P1} & \added{101871} & \added{3.62E-4} & \added{1.94E-4} & \added{9.2} \\
%         \hline
%         \added{AFEM-P2} &  \added{414909} & \added{2.14E-6} & \added{1.10E-6} & \added{32.6} \\
%         \hline
%         \added{AFEM-P3} & \added{918226} & \added{5.51E-8} & \added{1.31E-8} & \added{82.9}\\
% 			\hline
%         \added{AFEM-P4} & \added{1793861} & \added{1.97E-10} & \added{1.13E-10} & \added{279.1}\\
% 			\hline
%         \added{AFCM} & \added{36000} & \added{6.35E-11} & \added{1.38E-10} & \added{206.6}\\
% 			\hline
% 		\end{tabular}
% 		\label{Table:AFEM}
% 	\end{center}
% \end{table}
For $J_n=4000$, we show in Figure~\ref{fig:fig1} that 
\begin{itemize}
    \item The exact solution $\phi$ and approximate solutions $\tilde{\phi}_k$ on $\Omega$ (Fig.\ref{3a} \ref{3b})
     
    \item  Relative errors and  training losses  (Fig.\ref{3d}) 
    % $\|e_k\|_{L^{\infty}}$ and $\|e_k\|_{L^{2}}$
    
    \item local surfaces of approximate solutions $\tilde{\phi}_k$ on $[-0.1,0.1]\times[\frac{2}{3}-0.1,\frac{2}{3}+0.1]$ (Fig. \ref{3f})
\end{itemize}
\color{black}
    It is proved  in Proposition 2.9 of \cite{ming2025spectral}, the error between the exact and estimated solutions is bounded by the training loss. The results in Figure~\ref{fig:fig1} indicate that the relative training loss decreases with adaptive iterations, which implies a tighter error bound and thereby verifies the effectiveness of the proposed method.

\begin{figure}[htp]
	\center
        {\subfigure[partition hyperplane density $(k=0)$ ]{		 
            \includegraphics[scale=0.23]{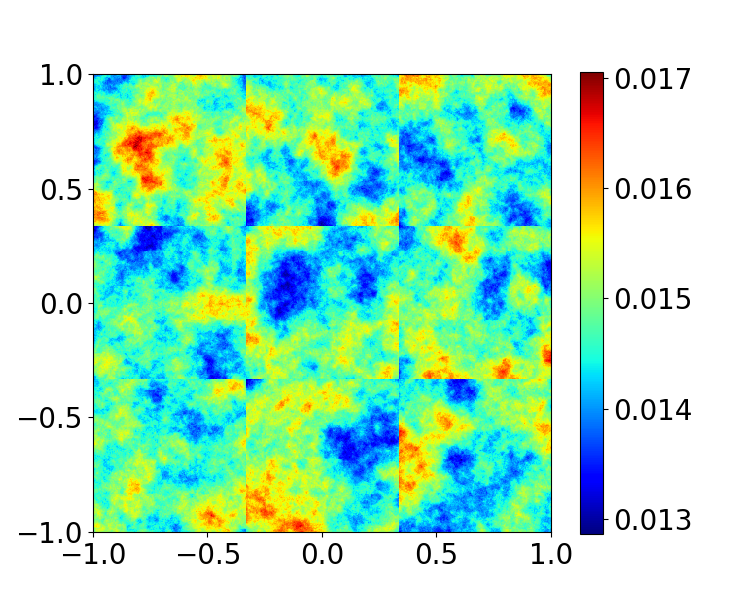}}
        \subfigure[collocation points $(k=0)$ ]{		 
            \includegraphics[scale=0.23]{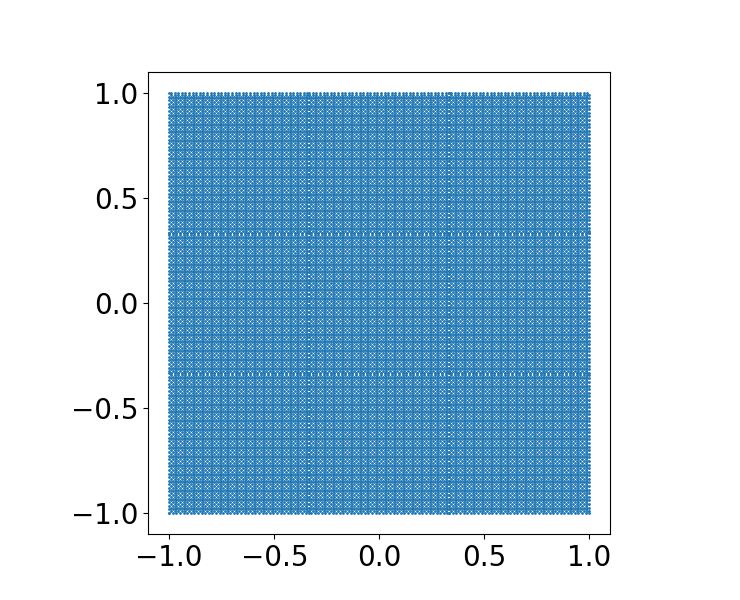}}
        \subfigure[shape parameter $(k=0)$ ]{		 
            \includegraphics[scale=0.23]{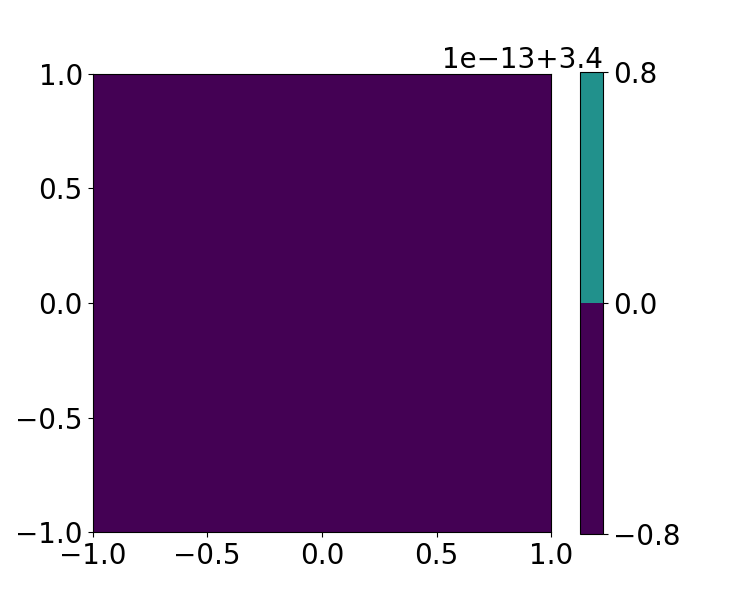}}
            \subfigure[partition hyperplane density  $(k=K)$  ]{		
            \includegraphics[scale=0.23]{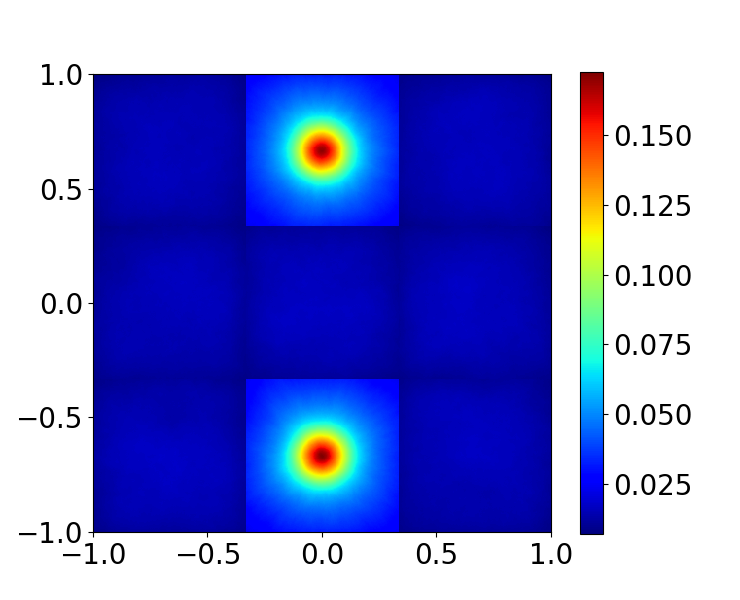}}
            \subfigure[collocation points $(k=K)$  ]{		
            \includegraphics[scale=0.23]{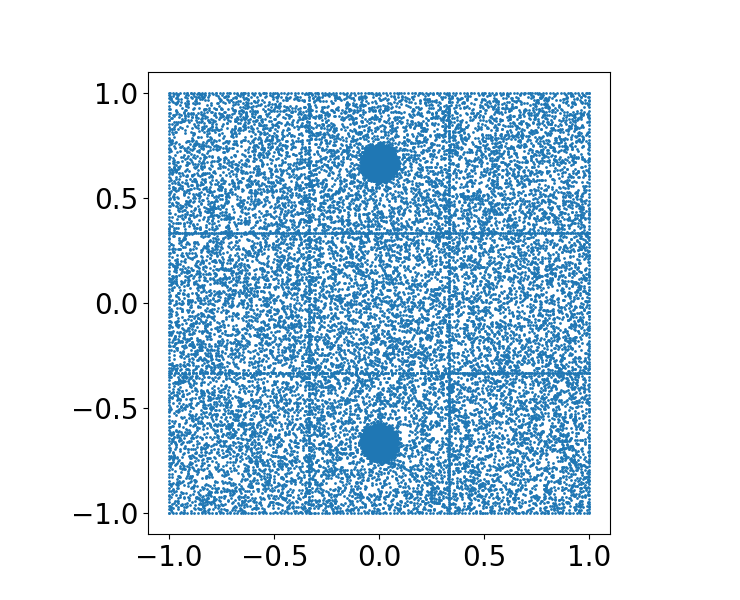}}
            \subfigure[shape parameter $(k=K)$ ]{		 
            \includegraphics[scale=0.23]{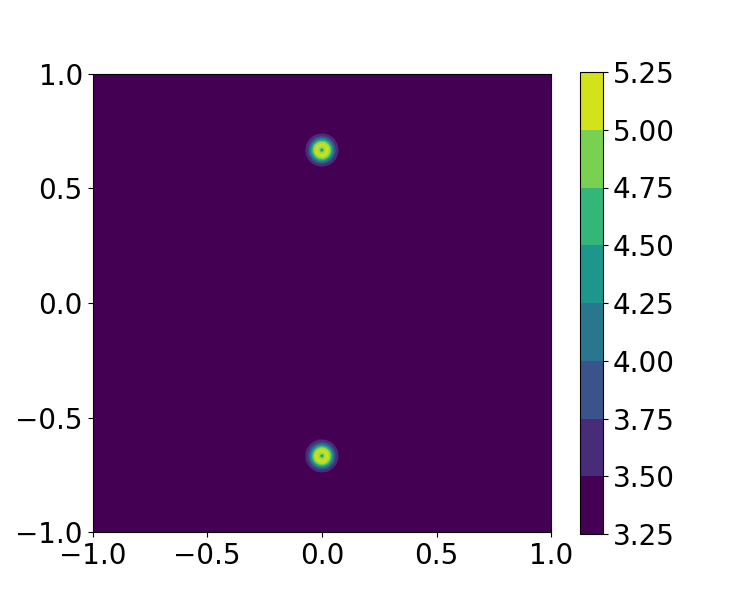}}
            }
	\caption{The partition hyperplane density, the collocation points and the shape parameters for \eqref{eq:exactsolution1} with $J_n=4000$ at  $k=0$ and $k=K=4$ iterations.}
	\label{fig:fig2}
\end{figure}

This case demonstrates the effectiveness of AFCM in handling problems with concentrated steep gradients. The initially uniformly distributed partition hyperplanes and collocation points provide insufficient resolution near the two peaks, leading to relatively large initial errors (e.g., $L^\infty$ error $\approx 1.82\times 10^{-1}$ for $J_n=4000$). After four adaptive iterations, the error decreases significantly to the order of $1.07\times 10^{-10}$, a reduction of over $9$ orders of magnitude. Figure ~\ref{fig:fig2} further shows that the partition hyperplane density, collocation point distribution, and shape parameters all become noticeably concentrated in the  two peaks region. This indicates that AFCM can dynamically enhance local expressive power based on gradient information, achieving higher approximation accuracy in high-gradient regions.
\color{black}

\subsection{Two-Dimensional Poisson equation with line near-singularity} \label{sec04:sec03}
We investigate a challenging test case of the two-dimensional Poisson equation \eqref{eq:poissoneq} defined on the computational domain $\Omega = (-1, 1)^2$, where the exact solution exhibits line-concentrated near-singular behavior with sharp gradients along the straight line $x=y/20$. The exact solution is given by
\begin{equation}\label{eq:exactsolution3}
    \phi(x, y) = e^{-1000 \left( x- \frac{y}{20} \right)^2}.
\end{equation}

The Dirichlet boundary condition $g(x, y)$ on $\partial\Omega$ and the source term $f(x, y)$ are analytically derived directly from the exact solution. We set $m = 1.8 \times 10^5$ sampling points for the test. The proposed Algorithm \ref{alg:Algorithm 1} is adopted to resolve the steep gradient along the singularity line $x=y/20$ via adaptive refinement.

Table \ref{Table:3} quantitatively reports the initial shape parameter $\gamma_n$, the initial ($k=0$) and final ($k=K=4$) $L^\infty$ and $L^2$ relative errors, \added{as well as the computational time and peak GPU memory consumption of the full adaptive procedure} under different $J_n$ for the line near-singularity case. For $J_n=4000$, after four adaptive iterations, the $L^\infty$ error drops from $\mathcal{O}(10^{-2})$ to $\mathcal{O}(10^{-10})$, achieving an accuracy improvement of eight orders of magnitude.

% \added{We compare AFCM with adaptive finite element methods (AFEM) using P2 and P3 elements, where AFEM adopts the standard residual-based error indicator for adaptivity. Both methods start from the same initial number of degrees of freedom (DoF), and the adaptive iterations terminate when both the relative $L^\infty$ and $L^2$ errors are below the prescribed bound 1.0E-6. As shown in Table \ref{Table:AFEM2}, AFCM achieves substantially higher accuracy with a fixed and much smaller DoF of 36000, whereas AFEM-P2 and AFEM-P3 require significant DoF growth to meet the stopping criterion. AFCM converges in only 1 iterations for the error bound, while AFEM-P2/P3 require 7/5 iterations with much larger DoF. Moreover, AFCM runs faster, taking only 58.3s. These results demonstrate that AFCM resolves steep gradients and near-singular solutions with better accuracy and computational efficiency.}

\begin{figure}[htp]
	\center
        {\subfigure[exact solution $\boldsymbol{\phi}$ ]{		 
            \includegraphics[scale=0.23]{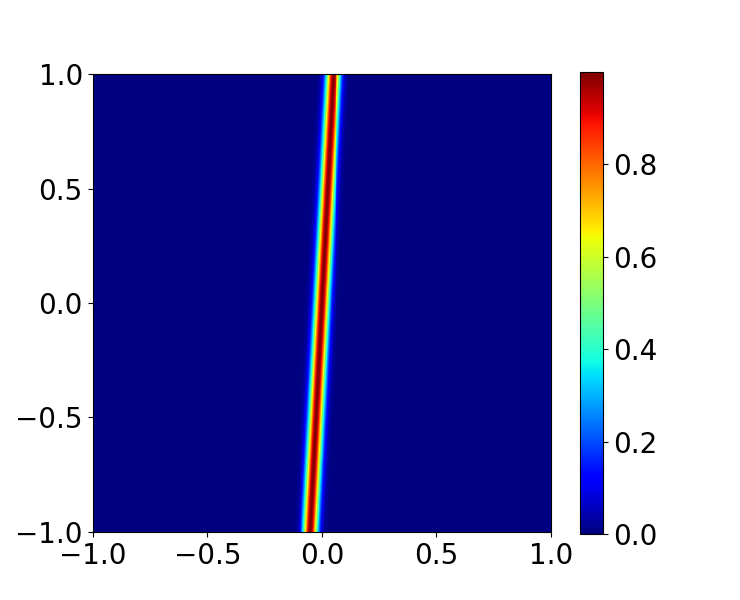}}
        \subfigure[approximte solution $\tilde{\boldsymbol{\phi}}_{0}  $ ]{		 
            \includegraphics[scale=0.23]{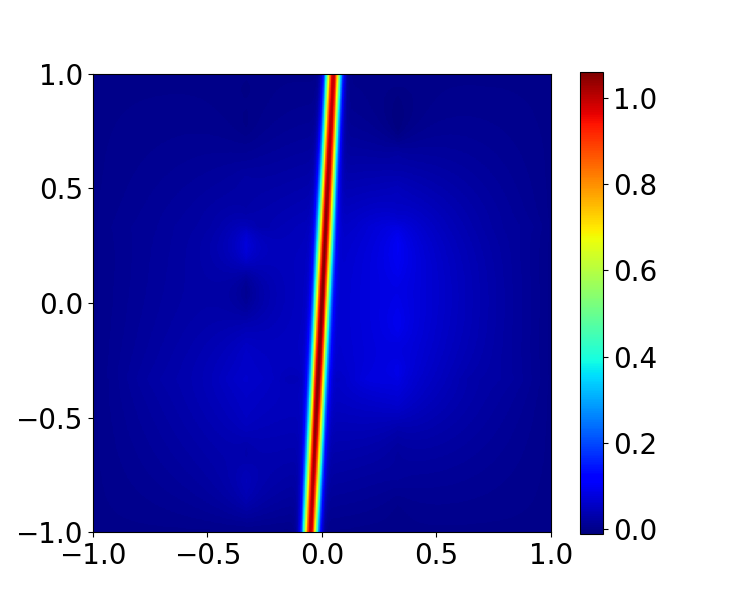}}
            \subfigure[local surface of  $\tilde{\phi}_0$ ]{		 
            \includegraphics[scale=0.32]{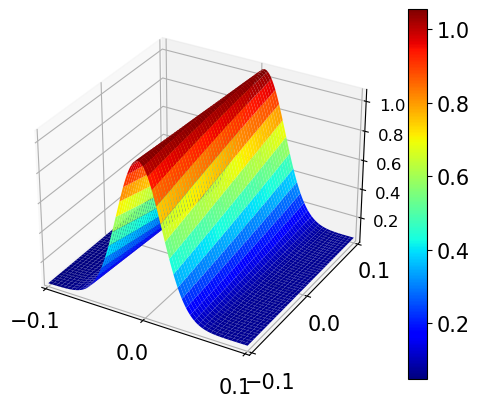}}
            \subfigure[Relative errors and training losses ]{		
            \includegraphics[scale=0.195]{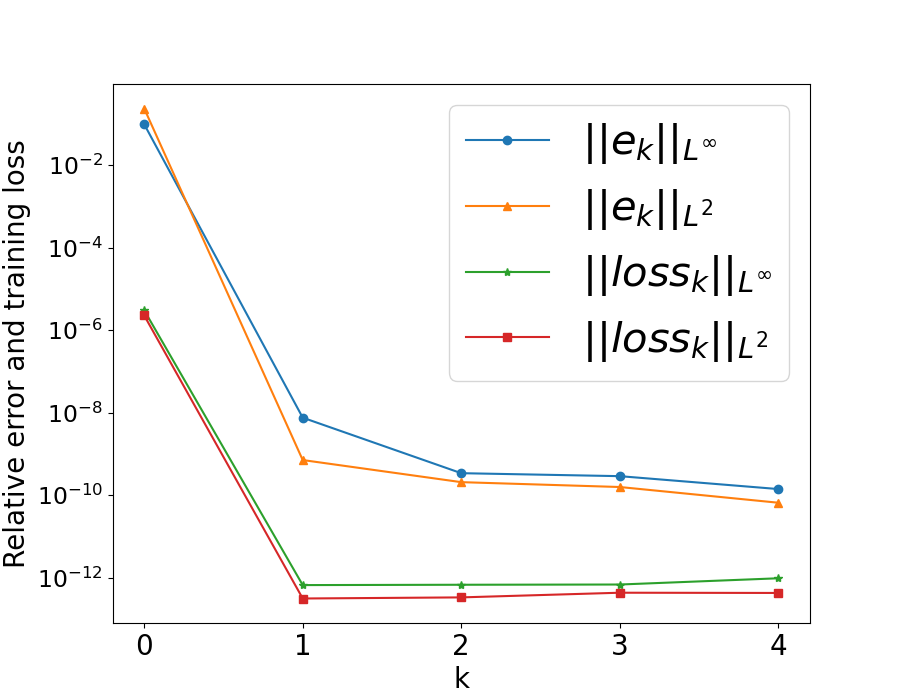}}
            \subfigure[approximte solution $\tilde{\boldsymbol{\phi}}_{K}  $  ]{		
            \includegraphics[scale=0.23]{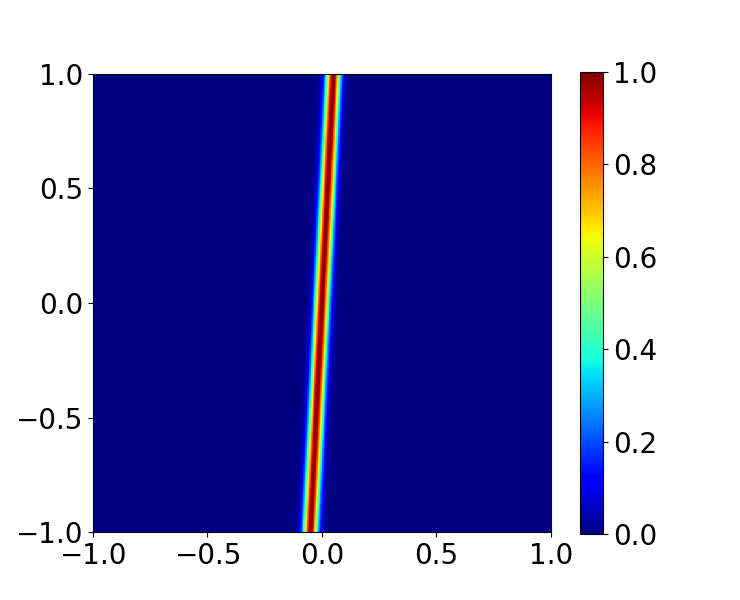}}
            \subfigure[local surface of  $\tilde{\phi}_K$]{		 
            \includegraphics[scale=0.32]{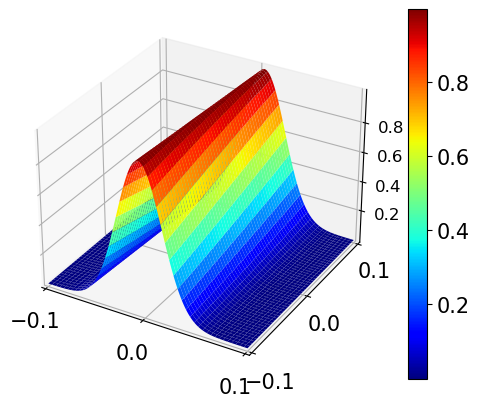}}
            }
	\caption{The exact solution and the numerical results for \eqref{eq:exactsolution3}  with $J_n=4000$ and $K=4$.}
	\label{fig:fig5}
\end{figure}

Figure \ref{fig:fig5} visualizes the exact solution, initial and final numerical approximations, local surface profiles, as well as the evolution of relative errors and training losses for the test case with $J_n=4000$ and $K=4$ iterations. The results clearly show that the solution possesses extremely steep gradients concentrated along the prescribed line.
It is also observed that the relative training loss gradually decreases with each adaptive iteration, which further verifies the effectiveness of the proposed adaptive method.

\begin{figure}[htp]
	\center
        {\subfigure[partition hyperplane density  $(k=0)$ ]{		 
            \includegraphics[scale=0.23]{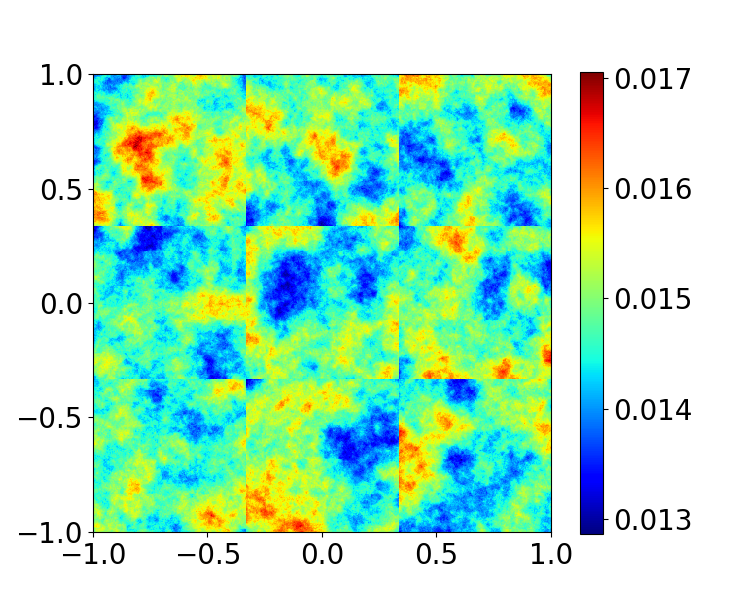}}
        \subfigure[collocation points $(k=0)$ ]{		 
            \includegraphics[scale=0.23]{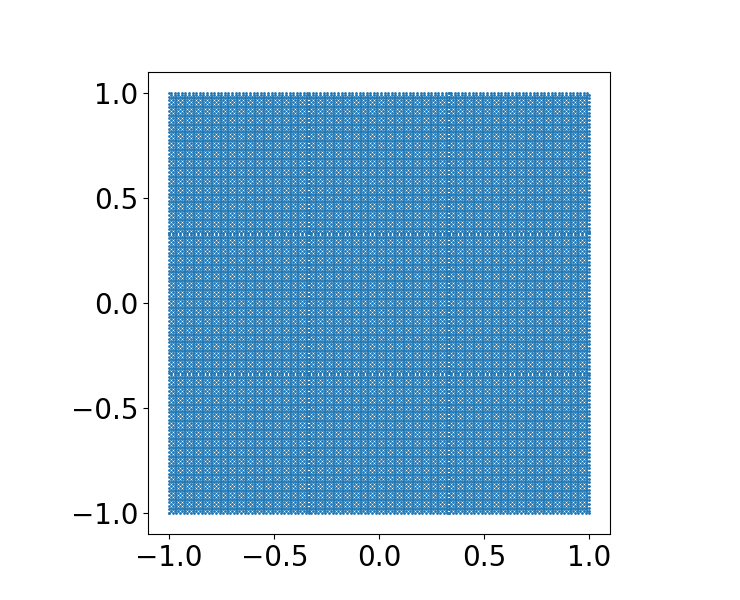}}
        \subfigure[shape parameter $(k=0)$ ]{		 
            \includegraphics[scale=0.23]{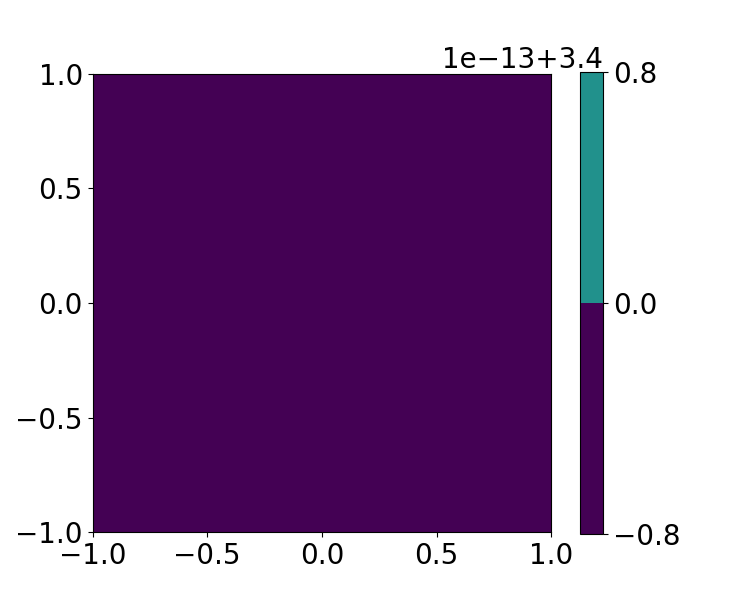}}
            \subfigure[partition hyperplane density  $(k=K)$  ]{		
            \includegraphics[scale=0.23]{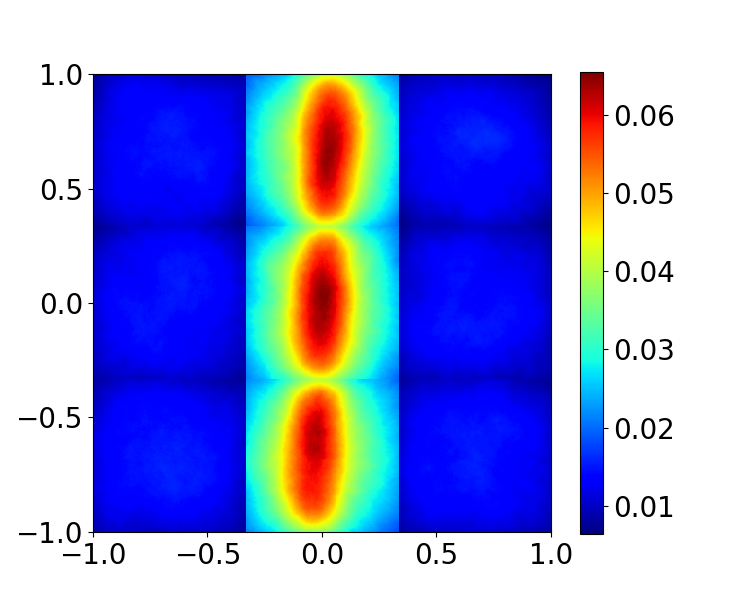}}
            \subfigure[collocation points $(k=K)$  ]{		
            \includegraphics[scale=0.23]{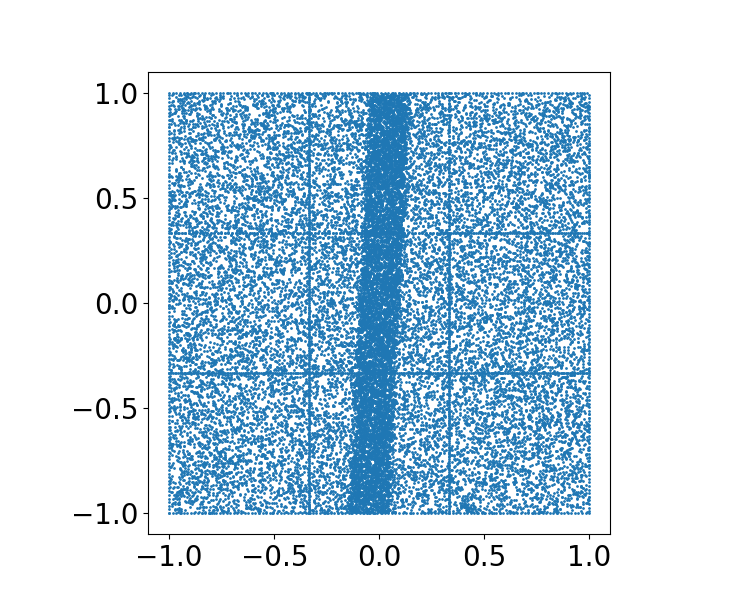}}
            \subfigure[shape parameter $(k=K)$ ]{		 
            \includegraphics[scale=0.23]{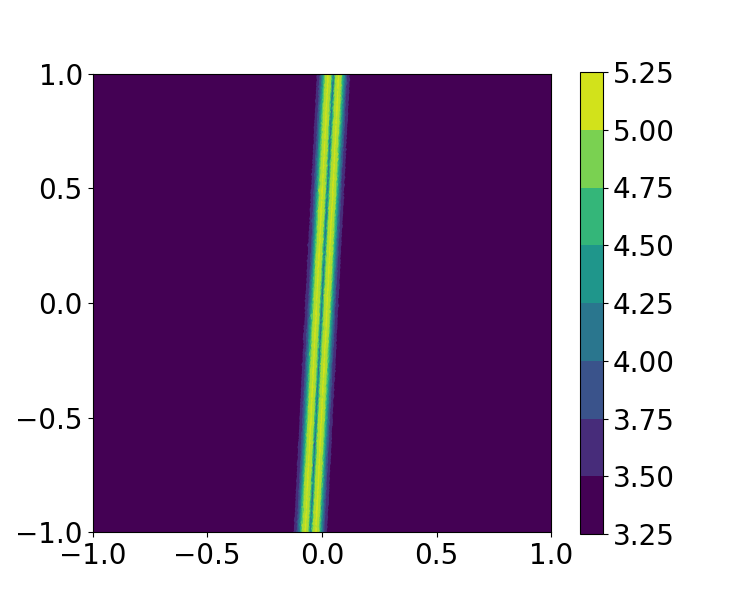}}
            }
	\caption{The partition hyperplane density, the collocation points and the shape parameters for \eqref{eq:exactsolution3} with $J_n=4000$ at $k=0$ and $k=K=4$ iterations}
	\label{fig:fig6}
\end{figure}

The adaptive refinement mechanism is further illustrated in Figure \ref{fig:fig6}. It can be observed that partition hyperplanes are densely clustered along the near-singularity line, collocation points automatically concentrate in the critical gradient region, and the shape parameters are adaptively increased substantially to well capture the sharp solution variation.
Numerical results demonstrate that the AFCM can accurately detect line-type near-singularity regions, automatically refine local numerical resolution, and efficiently resolve PDE problems with strong directional and steep gradients.

\begin{table}[htp]
	\begin{center} 
 \caption{\deleted{The initial shape parameter $\gamma_n$ and the relative errors for \eqref{eq:exactsolution3}  with different $J_n$ at  $k=0$ and $k=K=4$ iterations}\added{The initial shape parameter $\gamma_n$, the relative errors at $k=0$ and $k=K=4$ iterations, the computational time, and the peak GPU memory consumption  for \eqref{eq:exactsolution3} with different $J_n$}}
  \vspace{10pt}
		\begin{tabular}{|l|l|l|l|l|l|l|l|l|l|}
			\hline
            
		\multirow{2}{*}{$m$}& $J_n$ & 1500 & 2000  & 3000 & 4000         \\
        \cline{2-6}
        &  $\gamma_n$ & 2.0 & 2.6  & 2.8 & 3.4       \\
			\hline
	\multirow{6}{*}{$1.8\times 10^5$}& $\vert\vert e_{0}\vert\vert_{L^{\infty}} $  & 2.40E-1 & 3.30E-1 & 2.44E-1 & 9.80E-2     \\
    \cline{2-6}
		&  $\vert\vert e_{0}\vert\vert_{L^2} $  & 2.25E-1 & 3.06E-1 & 3.17E-1 & 2.36E-1   \\
            \cline{2-6}
		&  $\vert\vert e_{K=4}\vert\vert_{L^\infty} $  & 8.58E-6 & 4.54E-9 & 1.15E-9 & 1.41E-10   \\
                \cline{2-6}
		&  $\vert\vert e_{K=4}\vert\vert_{L^2} $  & 8.00E-6 & 1.82E-9 & 1.09E-9 & 6.56E-11   \\
				\cline{2-6}
        &  \added{time (s)}  & \added{48.8} & \added{65.9} & \added{122.3} & \added{195.8}   \\
                \cline{2-6}
		&  \added{Memory (GB)}  & \added{16.1} & \added{21.1} & \added{30.9} & \added{40.82}   \\
			\hline
		\end{tabular}
		\label{Table:3}
	\end{center}
\end{table}

\subsection{One-Dimensional Burgers equation with one line}
\label{sec04:sec05}
%In this section, we test the RK-RFM for multiphase flow problem \eqref{eq:multiphasepde}  by  numerically simulating 
Consider the one-dimensional Burgers equation defined on the computational domain $\Omega\times(0,T]=(0,1)\times(0,1]$,
\begin{equation}\label{eq:burgerseq}
\left\{\begin{aligned}
\boldsymbol{\phi}_t(x, t) +\boldsymbol{\phi}(x, t)\boldsymbol{\phi}_x(x, t)-\epsilon \boldsymbol{\phi}_{xx}(x, t) &=f(x, t), & & (x, t) \in \Omega\times(0,T], \\
\boldsymbol{\phi}(x, t) &=g(x, t), & & (x, t) \in \partial \Omega\times(0,T],\\
\boldsymbol{\phi}(x, 0) &=h(x), & & x \in \Omega,
\end{aligned}\right.
\end{equation}
the exact solution is given by
\begin{equation}\label{eq:exactsolution5}
\begin{aligned}
\boldsymbol{\phi}(x,t) = \frac{1}{1+e^{\frac{x-t}{2\epsilon}}}.
\end{aligned}
\end{equation}

The Dirichlet boundary condition $g(x, t)$ on $\partial \Omega\times(0,T]$, the initial condition $h(x)$ on $\Omega$ and the function $f(x, t)$ are given by the exact solution.  For this case, we treat the time variable $t$ as the spatial variable $y$.  We set $\epsilon=0.006$ and $m = 9.0\times 10^4$.  Algorithm \ref{alg:Algorithm 1} is employed to estimate $\boldsymbol{\phi}(x,t)$. 
We utilize Picard’s iterative methods to handle nonlinearity \deleted{address nonlinearity}, with the number of iterations set to 40. 

\begin{figure}[htp]
	\center
        {\subfigure[exact solution $\boldsymbol{\phi}$ ]{		 
            \includegraphics[scale=0.23]{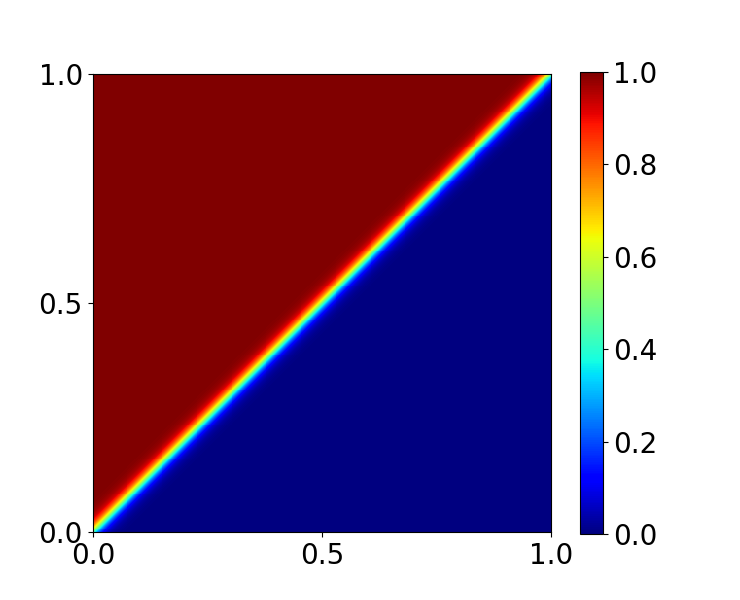}}
        \subfigure[approximte solution $\tilde{\boldsymbol{\phi}}_{0}  $ ]{		 
            \includegraphics[scale=0.23]{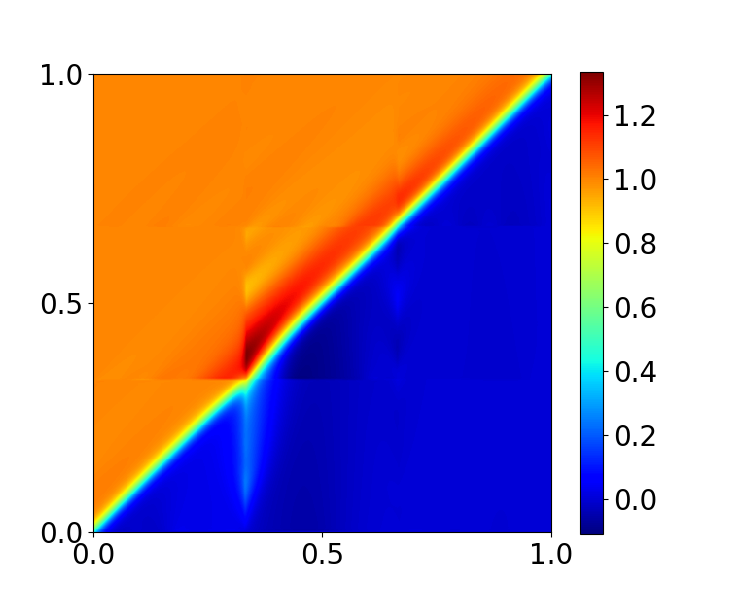}}
            \subfigure[local surface of  $\tilde{\phi}_0$ ]{		 
            \includegraphics[scale=0.32]{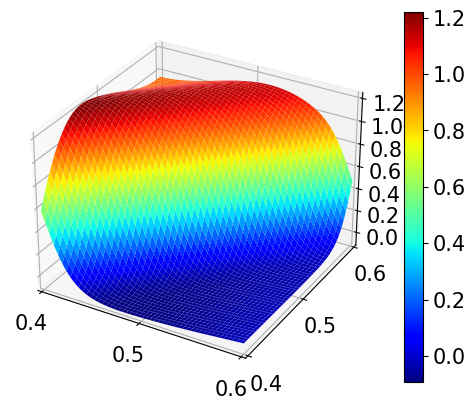}}
            \subfigure[Relative errors $\vert\vert e_{k}\vert\vert_{L^\infty}$ and $\vert\vert e_{k}\vert\vert_{L^2}$ ]{		
            \includegraphics[scale=0.191]{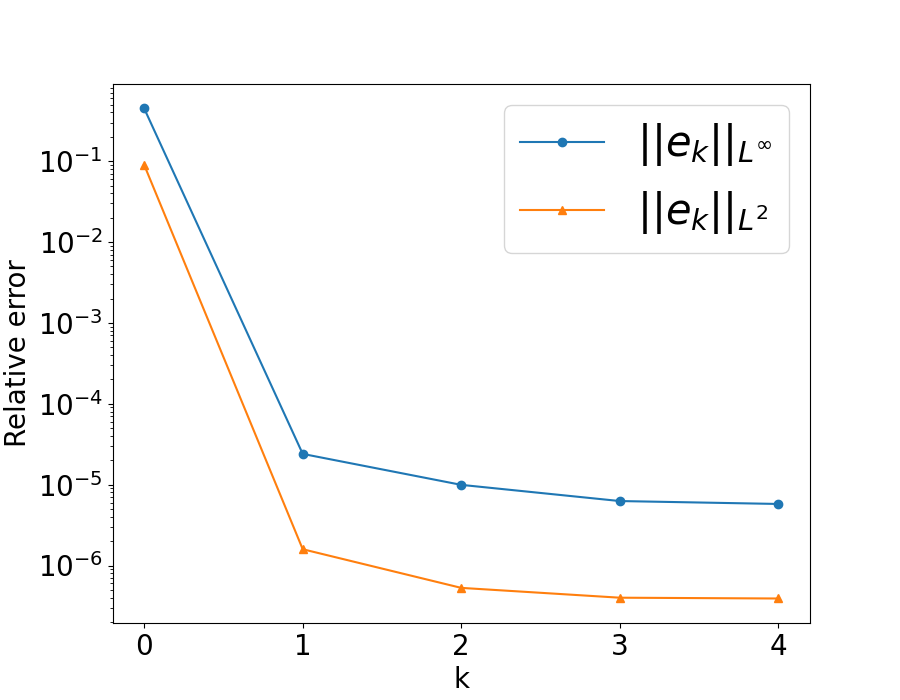}}
            \subfigure[approximte solution $\tilde{\boldsymbol{\phi}}_{K}  $  ]{		
            \includegraphics[scale=0.23]{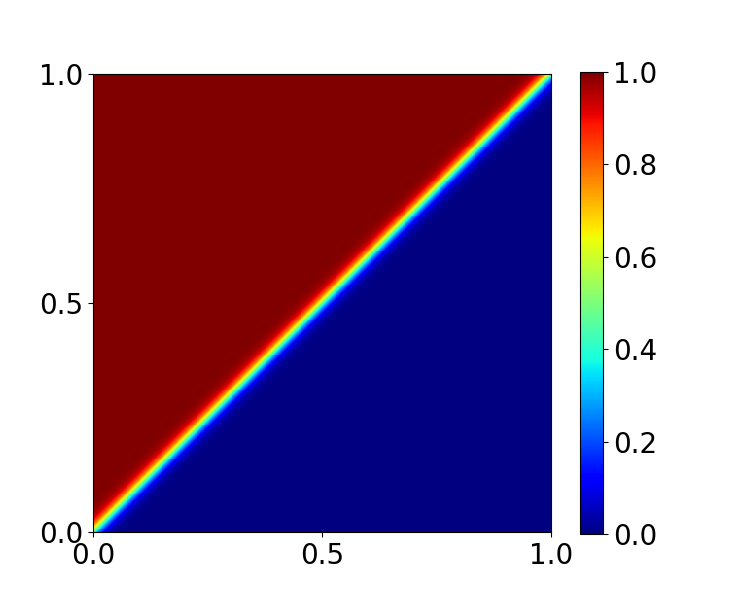}}
            \subfigure[local surface of  $\tilde{\phi}_K$]{		 
            \includegraphics[scale=0.32]{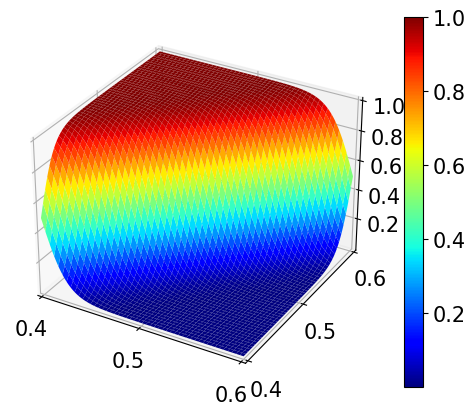}}
            }
	\caption{The exact solution and the numerical results for \eqref{eq:exactsolution5}  with $J_n=4000$ and $K=4$.}
	\label{fig:fig9}
\end{figure}

\begin{figure}[htp]
	\center
        {\subfigure[partition hyperplane density  $(k=0)$ ]{		 
            \includegraphics[scale=0.23]{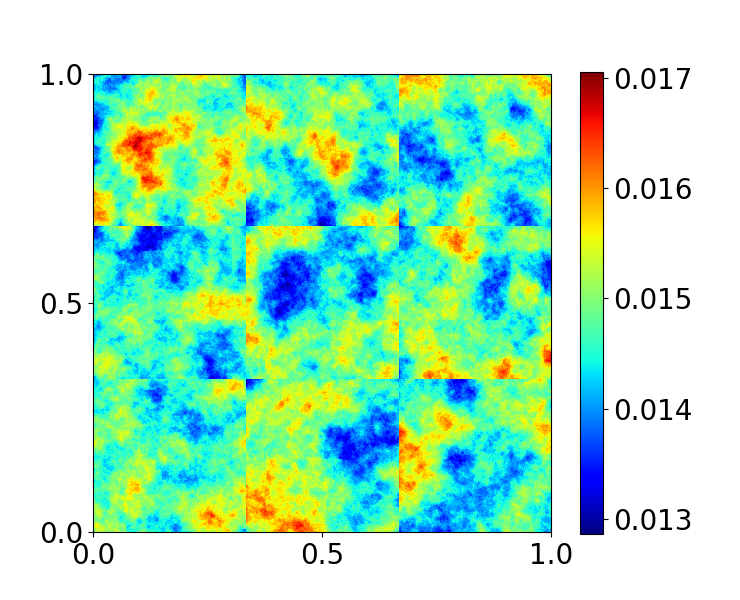}}
        \subfigure[collocation points $(k=0)$ ]{		 
            \includegraphics[scale=0.23]{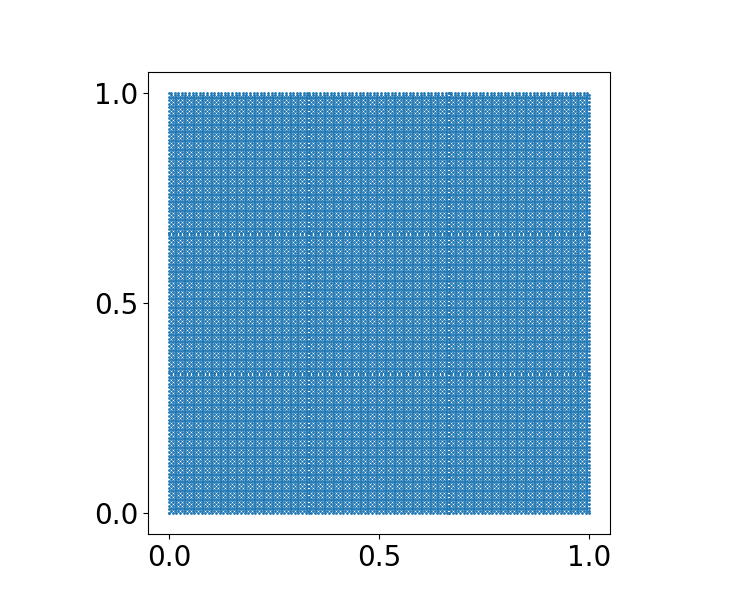}}
        \subfigure[shape parameter $(k=0)$ ]{		 
            \includegraphics[scale=0.23]{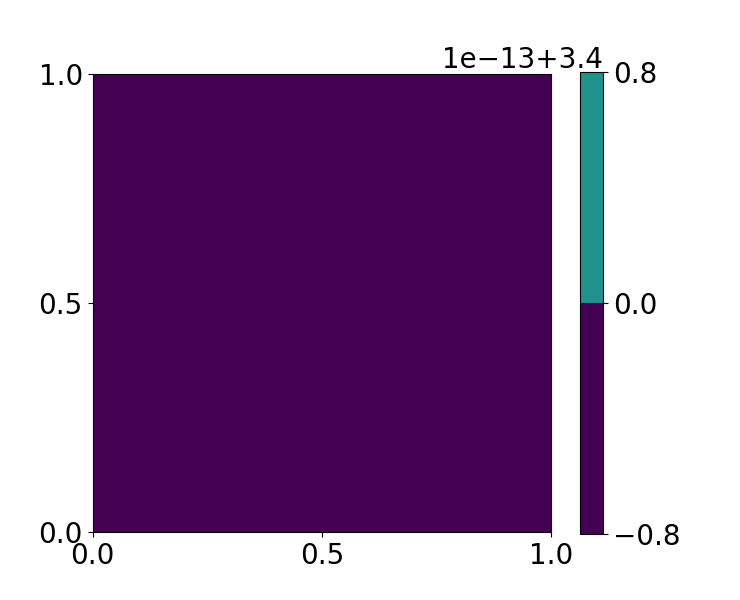}}
            \subfigure[partition hyperplane density  $(k=K)$  ]{		
            \includegraphics[scale=0.23]{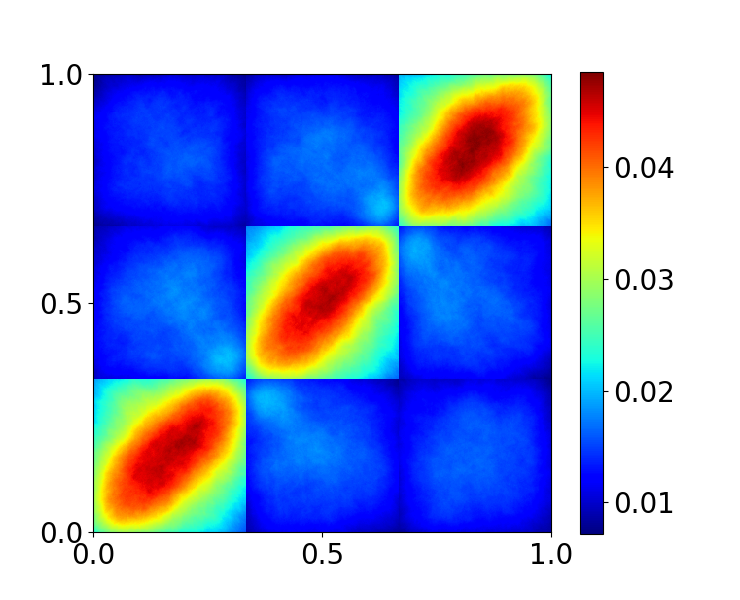}}
            \subfigure[collocation points $(k=K)$  ]{		
            \includegraphics[scale=0.23]{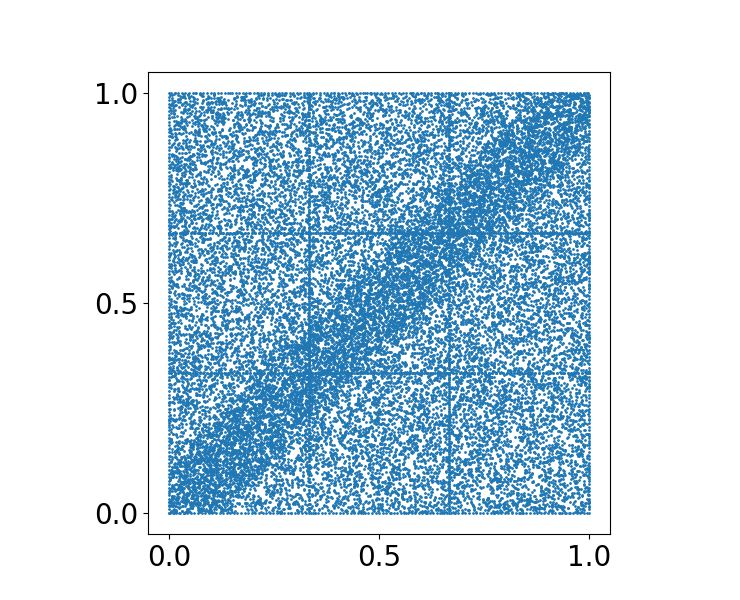}}
            \subfigure[shape parameter $(k=K)$ ]{		 
            \includegraphics[scale=0.23]{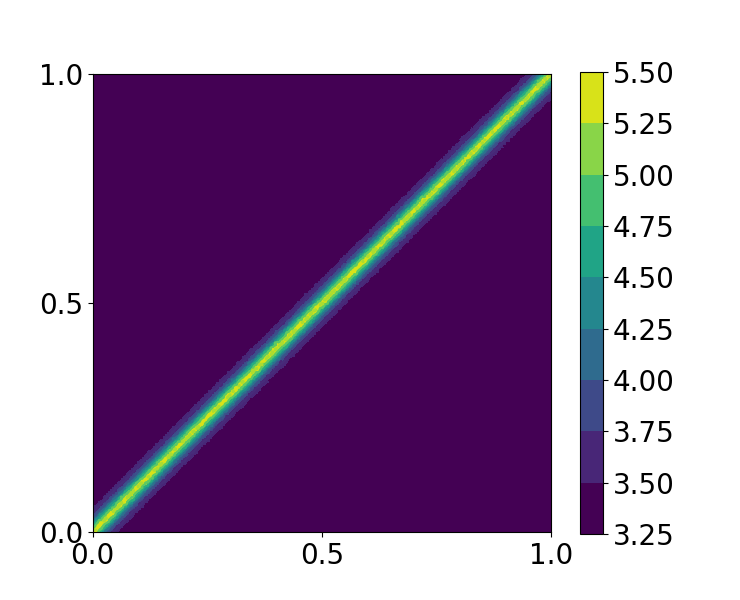}}
            }
	\caption{The partition hyperplane density, the collocation points and the shape parameters for \eqref{eq:exactsolution5} with $J_n=4000$ at $k=0$ and $k=K=4$ iterations}
	\label{fig:fig10}
\end{figure}

The initial shape parameter $\gamma_n$  and the relative errors for $\boldsymbol{\phi}(x,t)$  \deleted{with different $J_n$} at initial ($k=0$) and final ($k=K=4$) iterations, \added{as well as the computational time and peak GPU memory consumption of the full adaptive procedure,} are summarized in Table \ref{Table:5} \added{for different $J_n$}. For $J_n=4000$, Figure \ref{fig:fig9} displays the exact solution $\boldsymbol{\phi}$, the relative errors $\vert\vert e_{k}\vert\vert_{L^\infty}$ and $\vert\vert e_{k}\vert\vert_{L^2}$, the approximate solutions $\tilde{\boldsymbol{\phi}}_{k}$, and the local surfaces of approximate solution $\tilde{\phi}_k$ on $[0.4,0.6]^2$ at the initial ($k=0$) and final ($k=K=4$) iterations. The partition hyperplane density, the collocation points and the shape parameters with $J_n=4000$ at $k=0$ and $k=K=4$ iterations are shown in Figure \ref{fig:fig10}. \color{black} The results demonstrate that AFCM combined with Picard iteration can effectively handle nonlinear terms, and the adaptive mechanism can dynamically track time-dependent singularities, making it suitable for evolutionary PDEs.
\color{black}
\begin{table}[htp]
	\begin{center} 
 \caption{\deleted{The initial shape parameter $\gamma_n$ and the relative errors for \eqref{eq:exactsolution5}  with different $J_n$ at  $k=0$ and $k=K=4$ iterations}\added{The initial shape parameter $\gamma_n$, the relative errors at $k=0$ and $k=K=4$ iterations, the computational time, and the peak GPU memory consumption for \eqref{eq:exactsolution5} with different $J_n$}}
  \vspace{10pt}
		\begin{tabular}{|l|l|l|l|l|l|l|l|l|l|}
			\hline
            
		\multirow{2}{*}{$m$}& $J_n$ & 1500 & 2000  & 3000 & 4000         \\
        \cline{2-6}
        &  $\gamma_n$ & 2.0 & 2.6  & 2.8 & 3.4       \\
			\hline
	\multirow{6}{*}{$9.0\times 10^4$}& $\vert\vert e_{0}\vert\vert_{L^{\infty}} $  & 1.04E-2 & 1.60E-2 & 3.45E-2 & 4.48E-1     \\
    \cline{2-6}
		&  $\vert\vert e_{0}\vert\vert_{L^2} $  & 2.51E-3 & 4.10E-3 & 9.66E-3 & 8.97E-2   \\
            \cline{2-6}
		&  $\vert\vert e_{K=4}\vert\vert_{L^\infty} $  & 1.16E-2 & 3.20E-4 & 1.46E-5 & 5.77E-6   \\
                \cline{2-6}
		&  $\vert\vert e_{K=4}\vert\vert_{L^2} $  & 6.58E-4 & 2.40E-5 & 3.40E-6 & 3.91E-7   \\
				\cline{2-6}
      &  \added{time (s)}  & \added{1678.5} & \added{2250.7} & \added{4498.9} & \added{7196.0}   \\
                \cline{2-6}
        &  \added{Memory (GB)}  & \added{22.1} & \added{29.0} & \added{42.8} & \added{56.7}   \\
    \hline
		\end{tabular}
		\label{Table:5}
	\end{center}
\end{table}

\subsection{Two-Dimensional Heat equation with one peak}
\label{sec04:sec06}
Consider the two-dimensional Heat equation in $\Omega\times(0,T]=(-1,1)^2\times(0,2]$,
\begin{equation}\label{eq:2dheateq}
\left\{\begin{aligned}
\boldsymbol{\phi}_t(x,y,t) - \alpha\Delta\boldsymbol{\phi}(x,y,t) &=f(x,y,t), & & (x,y,t) \in \Omega\times(0,T], \\
\boldsymbol{\phi}(x,y, t) &=g(x,y,t), & & (x,y,t) \in \partial \Omega\times(0,T],\\
\boldsymbol{\phi}(x,y, 0) &=h(x,y), & & (x,y) \in \Omega,
\end{aligned}\right.
\end{equation}
 the exact solution is given by
\begin{equation}\label{eq:exactsolution6}
\begin{aligned}
\boldsymbol{\phi}(x,y,t) = e^{-1000\left((x-\frac{t}{10})^2+(y-\frac{t}{10})^2\right)}.
\end{aligned}
\end{equation}

The Dirichlet boundary condition $g(x,y,t)$ on $\partial \Omega\times(0,T]$, the initial condition $h(x,y)$ on $\Omega$ and the function $f(x,y,t)$ are given by the exact solution. For this example,  we use the Crank-Nicolson scheme for time discretization and apply the AFCM with $K$ adaptive iterations at each time step. The collocation points and feature functions after the adaptation of the previous time step are used as the initial collocation points and feature functions of the next time step. We set $\alpha=1000$, the number of time steps $N = 10$, the time interval $dt = 0.2$, $m = 1.26 \times 10^6$  and Algorithm \ref{alg:Algorithm 1} is employed to estimate $\boldsymbol{\phi}(x,y,t)$. 

\begin{figure}[htp]
	\center
        {\subfigure[exact solution $\boldsymbol{\phi}$ ]{		 
            \includegraphics[scale=0.236]{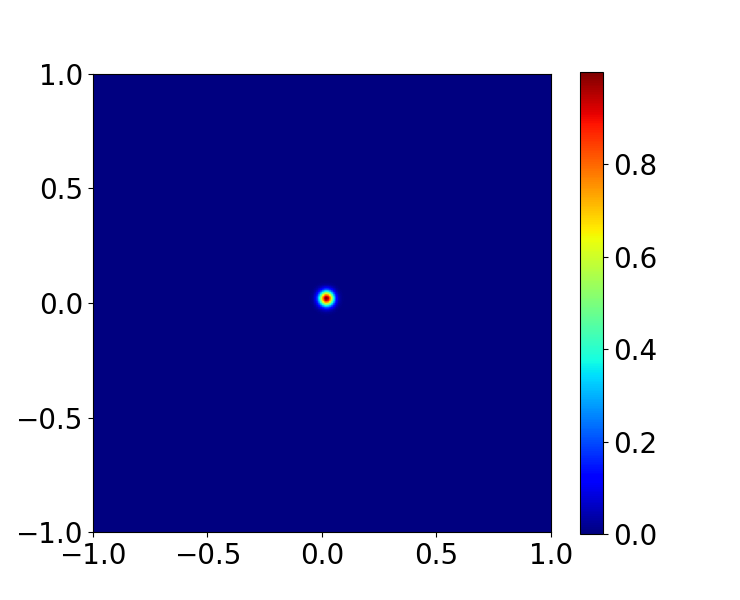}		 
            \includegraphics[scale=0.236]{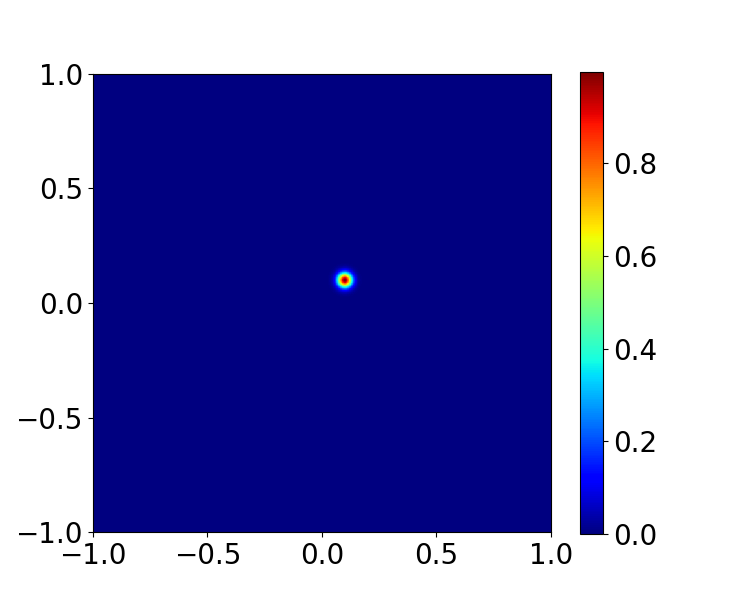}		 
            \includegraphics[scale=0.236]{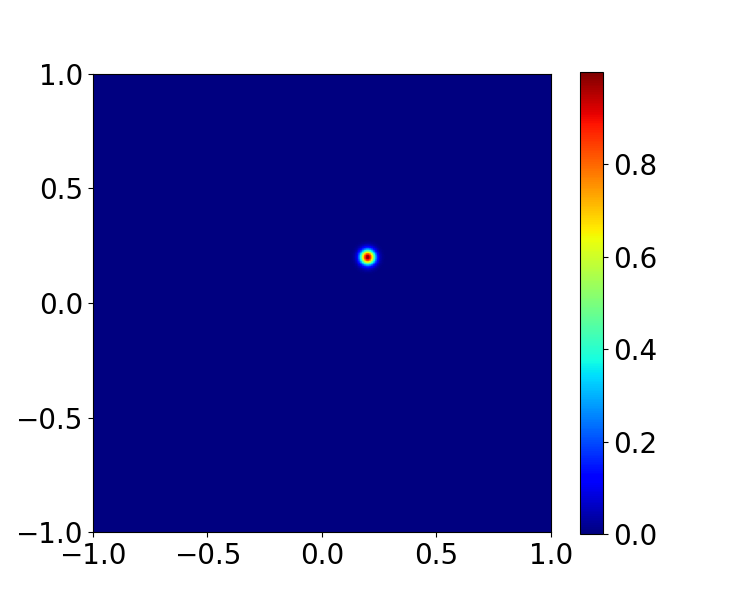}}
            \subfigure[approximte solution $\tilde{\boldsymbol{\phi}}_{K=0}$ ]{		
            \includegraphics[scale=0.236]{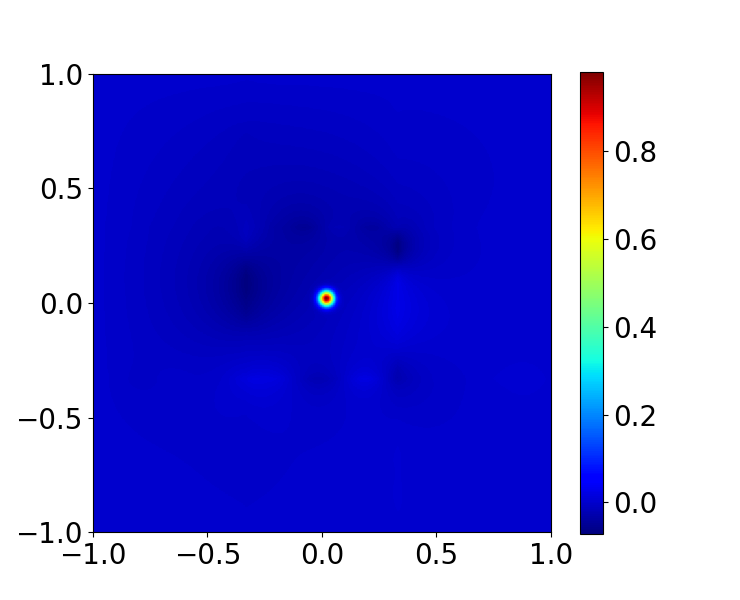}
            \includegraphics[scale=0.236]{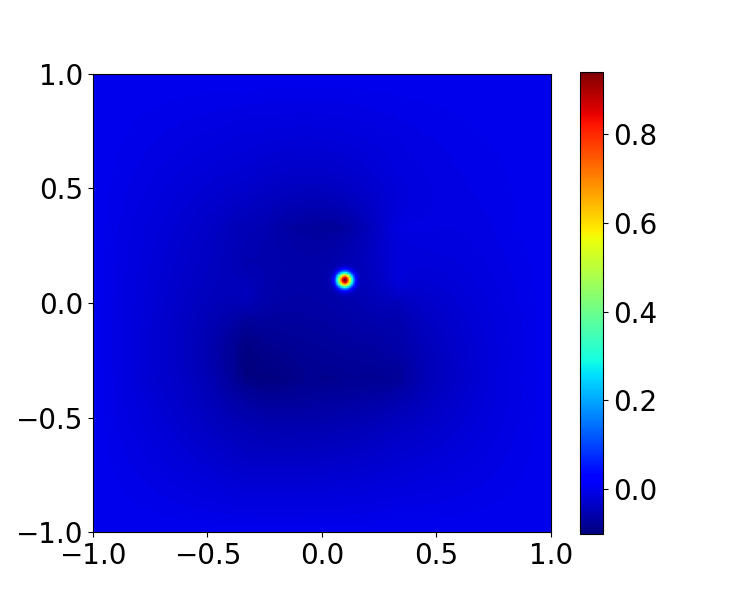}
            \includegraphics[scale=0.236]{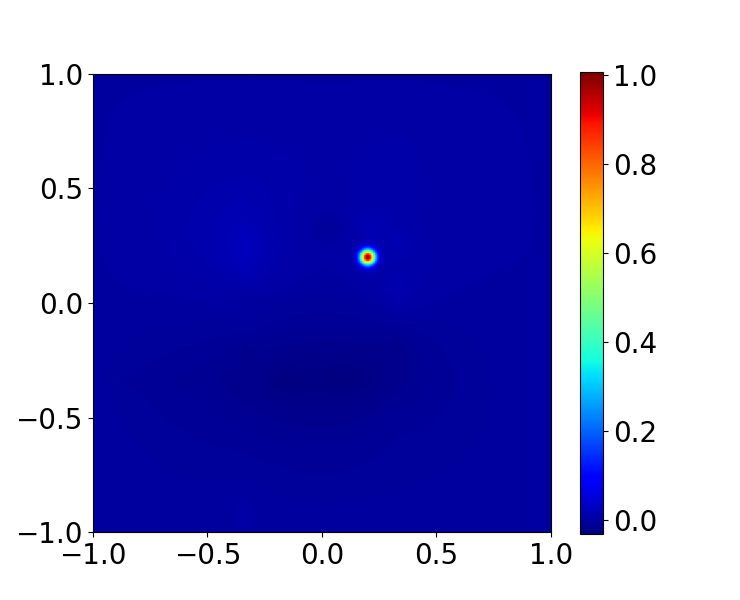}}
            \subfigure[local surface of  $\tilde{\phi}_{K=0}$ ]{		
            \includegraphics[scale=0.355]{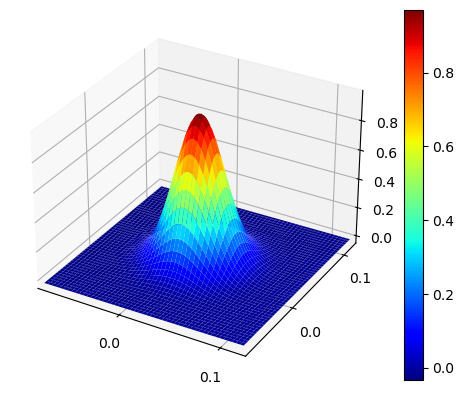}
            \includegraphics[scale=0.355]{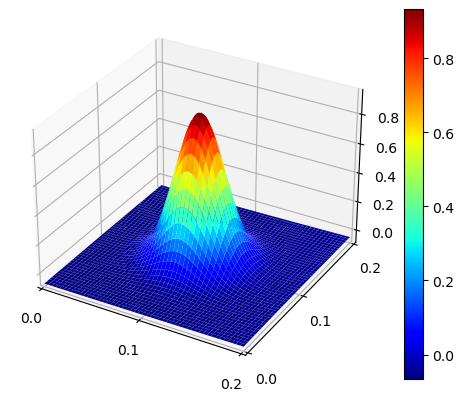}
            \includegraphics[scale=0.355]{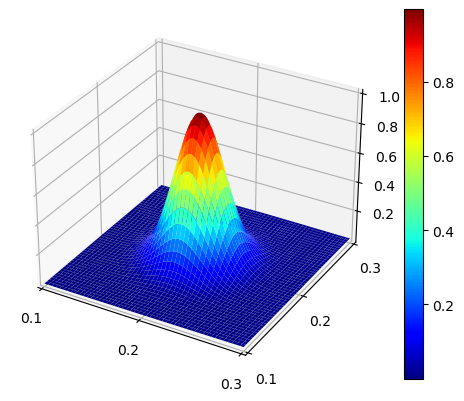}}
            }

	\caption{The exact solution and the numerical results for \eqref{eq:exactsolution6} at the final $(k = K)$ iteration and $t =0.2$ (left), $1.0$ (middle), $2.0$ (left)   with $J_n=4000$ and $K=0$.}
	\label{fig:fig11}
\end{figure}

\begin{figure}[htp]
	\center
        {\subfigure[exact solution $\boldsymbol{\phi}$ ]{		 
            \includegraphics[scale=0.236]{2dheat_as_1t5.png}		 
            \includegraphics[scale=0.236]{2dheat_as_5t5.png}		 
            \includegraphics[scale=0.236]{2dheat_as_10t5.png}}
            \subfigure[approximte solution $\tilde{\boldsymbol{\phi}}_{K=4}$ ]{		
            \includegraphics[scale=0.236]{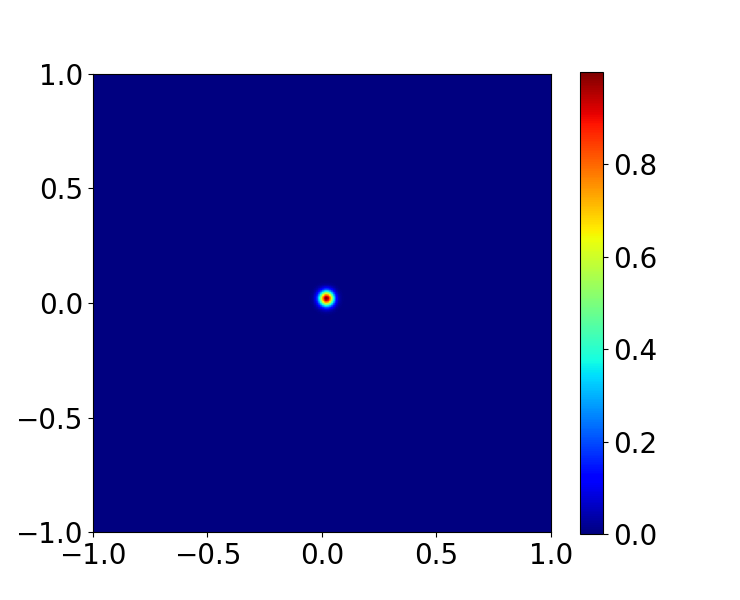}
            \includegraphics[scale=0.236]{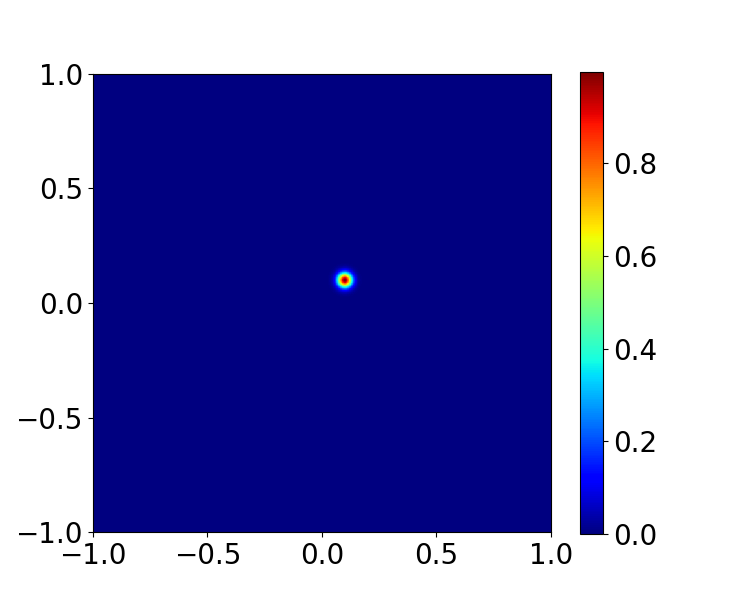}
            \includegraphics[scale=0.236]{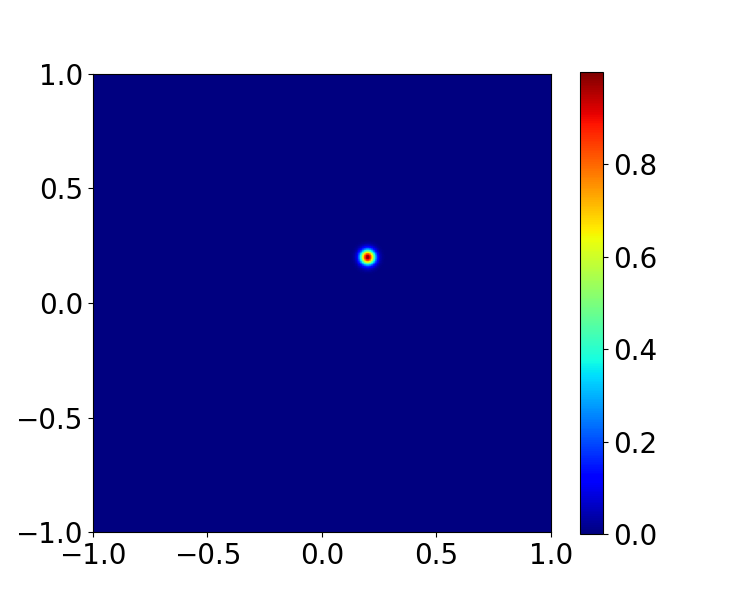}}
            \subfigure[local surface of  $\tilde{\phi}_{K=4}$ ]{		
            \includegraphics[scale=0.355]{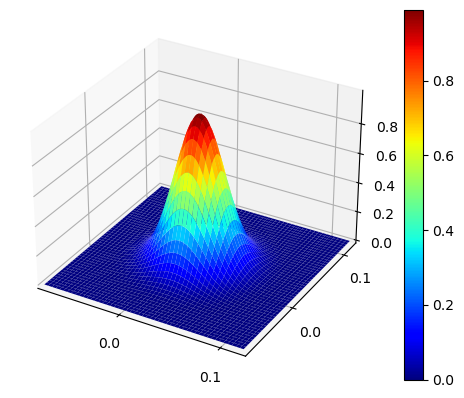}
            \includegraphics[scale=0.355]{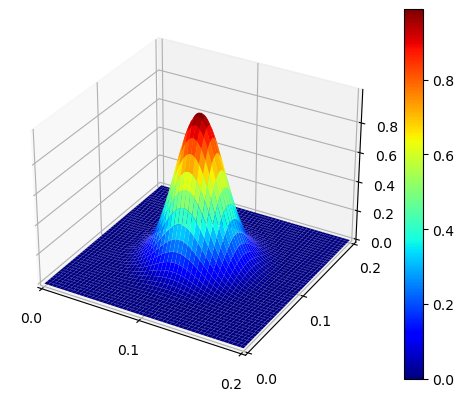}
            \includegraphics[scale=0.355]{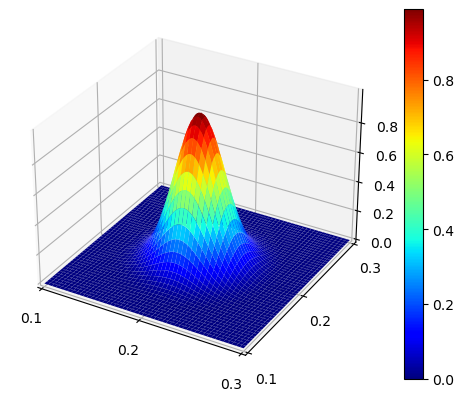}}
            }

	\caption{The exact solution and the numerical results for \eqref{eq:exactsolution6} at the final $(k = K)$ iteration and $t =0.2$ (left), $1.0$ (middle), $2.0$ (left)   with $J_n=4000$ and $K=4$.}
	\label{fig:fig12}
\end{figure}

\begin{figure}[htp]
	\center
        {\subfigure[  $k = 4, \ t=0.2$ ]{	
            \includegraphics[scale=0.23]{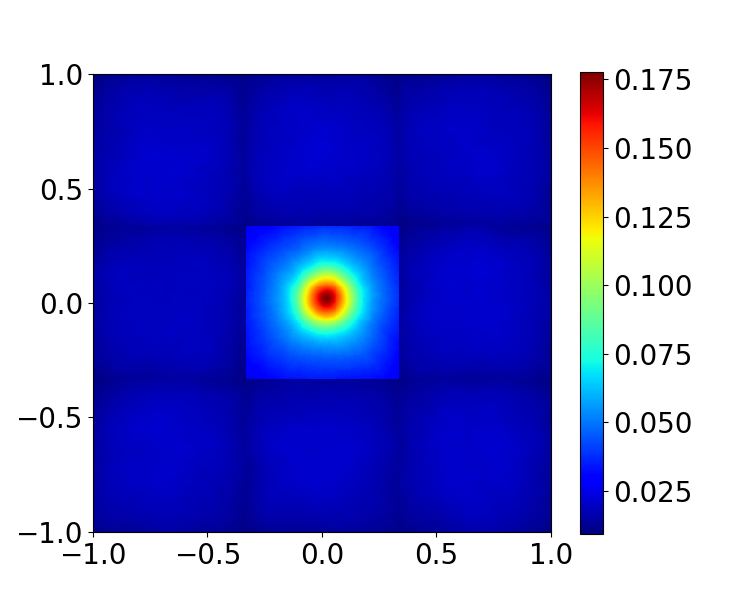}
            \includegraphics[scale=0.23]{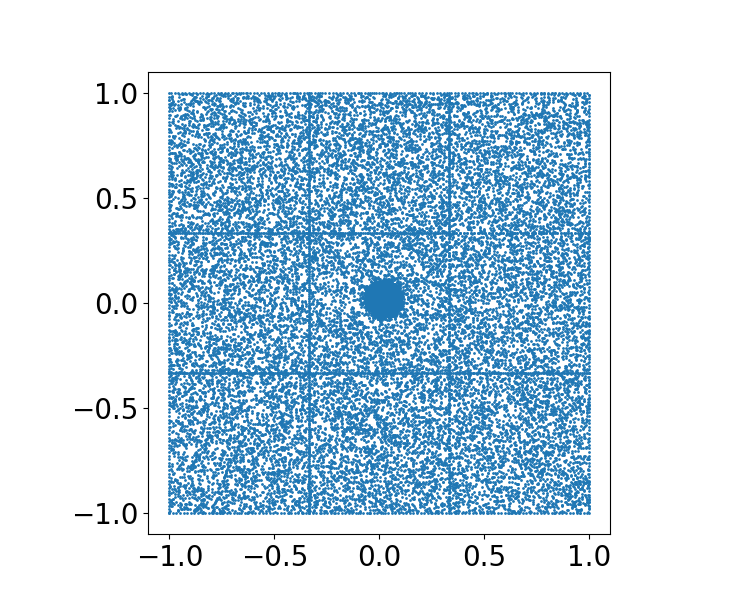}
            \includegraphics[scale=0.23]{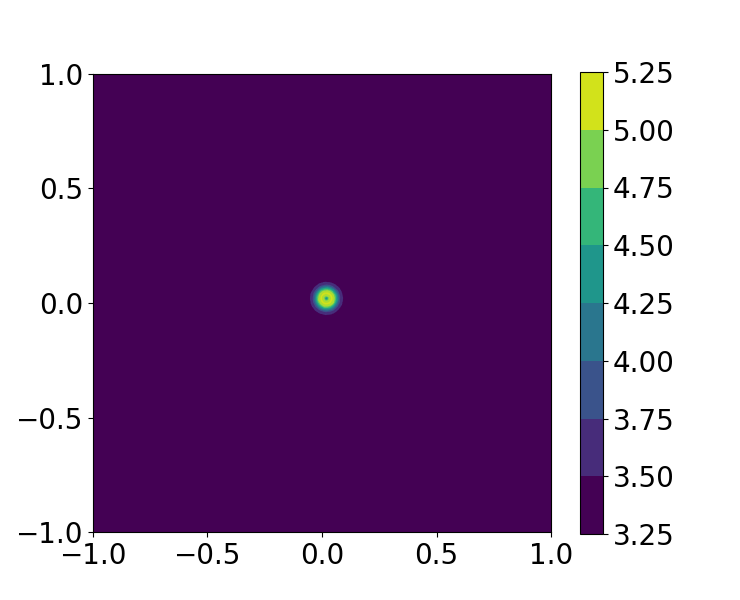}}
            \subfigure[$k = 4, \ t=1.0$ ]{		 
            \includegraphics[scale=0.23]{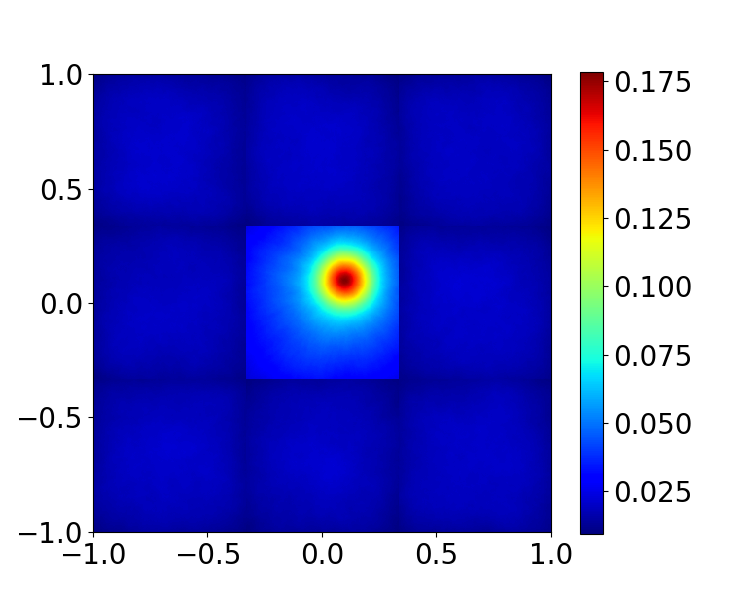}		
            \includegraphics[scale=0.23]{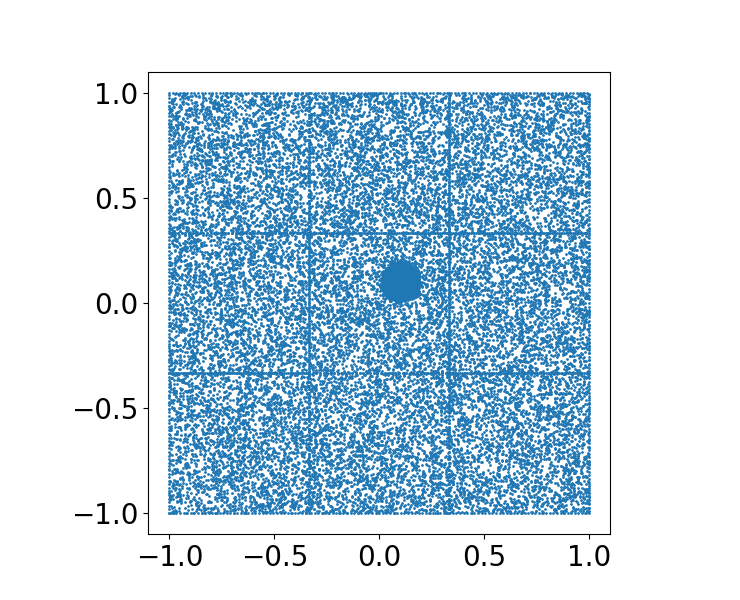}
            \includegraphics[scale=0.23]{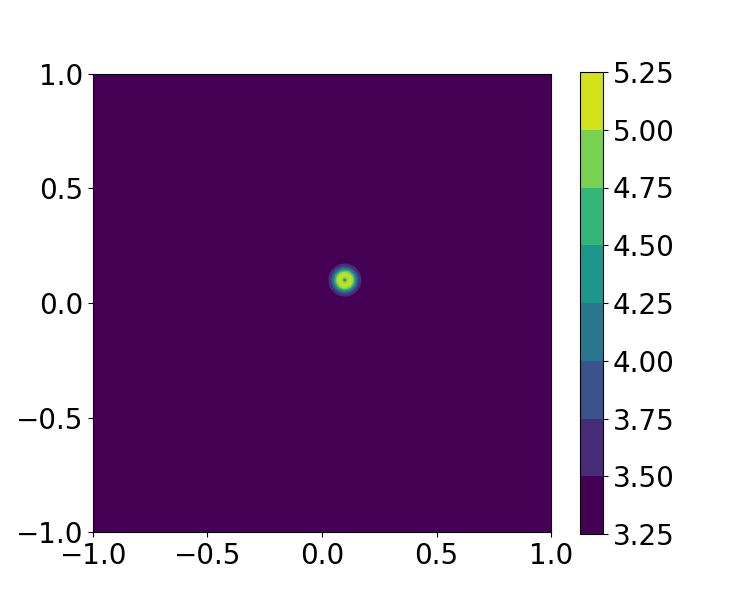}}
            
            \subfigure[$k = 4, \ t=2.0$ ]{		 
            \includegraphics[scale=0.23]{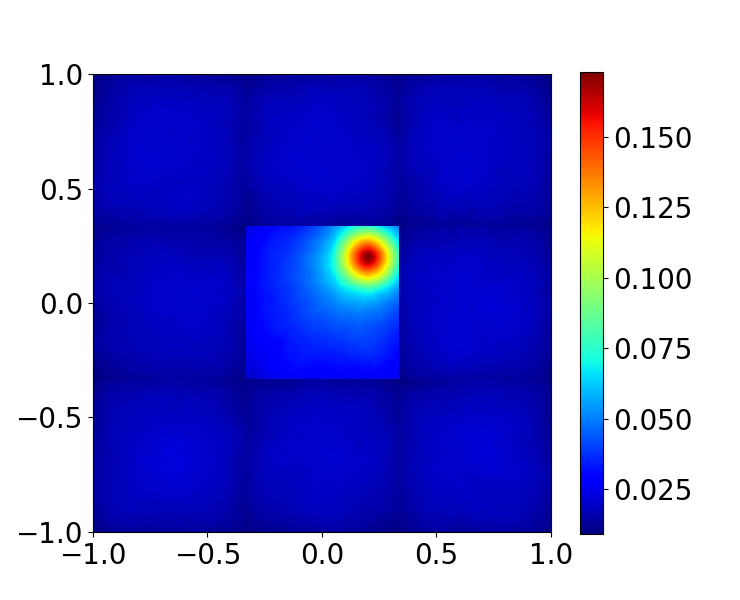}	 
            \includegraphics[scale=0.23]{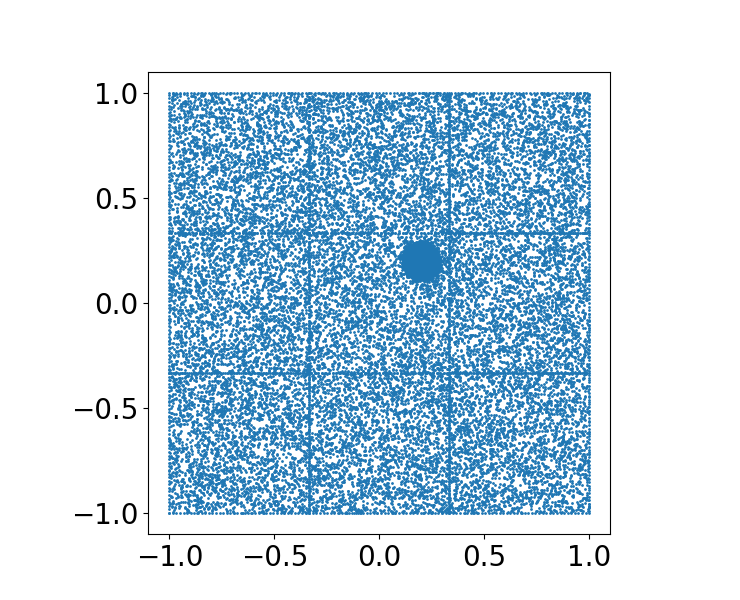}
            \includegraphics[scale=0.23]{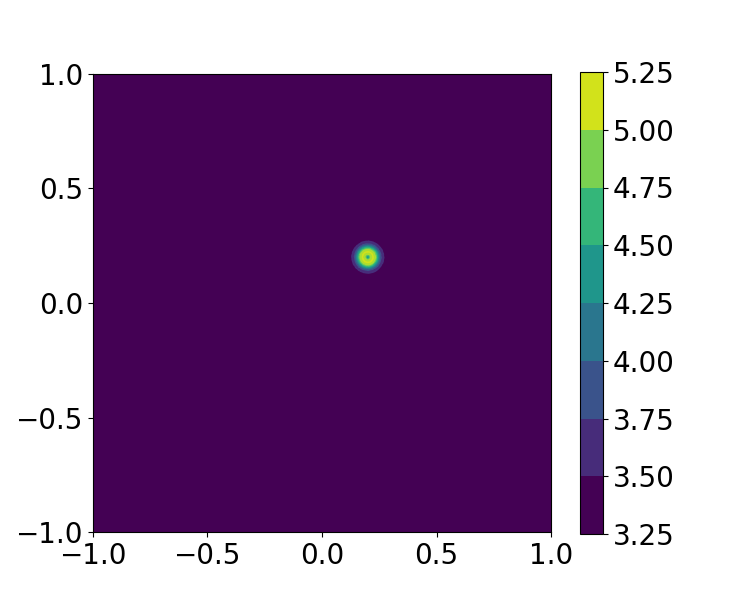}}
            }
	\caption{The partition hyperplane density (left), the collocation points (middle) and the shape parameters (right) for \eqref{eq:exactsolution6} with $K=4$ and $J_n=4000$ at the final $(k = K)$ iteration and $t =0.2, 1.0, 2.0$}
	\label{fig:fig13}
\end{figure}

The initial shape parameter $\gamma_n$  and the relative  errors for $\boldsymbol{\phi}(x,y,t)$ at the final $(k = K)$ iteration and $T=2.0$ with $K=0,4$, \deleted{and different $J_n$} \added{as well as the computational time and peak GPU memory consumption of the full adaptive procedure,} are summarized in Table \ref{Table:6} \added{for different $J_n$}. For $J_n=4000$, Figure \ref{fig:fig11} and Figure \ref{fig:fig12} display the exact solution $\boldsymbol{\phi}$, the approximate solutions $\tilde{\boldsymbol{\phi}}_{k}$, and the local surfaces of approximate solution $\tilde{\phi}_k$  on $[\frac{t}{10}-0.1,\frac{t}{10}+0.1]^2$ at the final $(k = K)$ iteration and $t =0.2, 1.0, 2.0$  with $K=0,4$. The partition hyperplane density, the collocation points and the shape parameters with $K=4$ and $J_n=4000$ at the final $(k = K)$ iteration and $t =0.2, 1.0, 2.0$ are shown in Figure \ref{fig:fig13}, which shows that  the positions of the concentration of the partition hyperplane and the collocation points, as well as the positions of the changes in the shape parameters, move with the shift of the peak position of the solution. \color{black} This indicates that AFCM can effectively integrate time discretization strategies (e.g., the Crank-Nicolson scheme) to achieve spatiotemporal adaptivity, making it suitable for PDEs with moving singularities or transient behaviors.

\color{black}

\begin{table}[htp]
	\begin{center} 
 \caption{\deleted{The initial shape parameter $\gamma_n$ and the relative errors for \eqref{eq:exactsolution6}  at the final $(k = K)$ iteration and $T=1.0$ with $K=0,4$ and different $J_n$}\added{The initial shape parameter $\gamma_n$, the relative errors at the final $(k = K)$ iteration and $T=1.0$ with $K=0,4$, the computational time, and the peak GPU memory consumption for \eqref{eq:exactsolution6} with different $J_n$}}
  \vspace{10pt}
		\begin{tabular}{|l|l|l|l|l|l|l|l|l|l|}
			\hline
            
		\multirow{2}{*}{$T$}&\multirow{2}{*}{$m$}& $J_n$ & 1500 & 2000  & 3000 & 4000         \\
        \cline{3-7}
        &   &   $\gamma_n$ & 2.0 & 2.6  & 2.8 & 3.4       \\
			\hline
	\multirow{6}{*}{$2.0$}&\multirow{6}{*}{$1.26 \times 10^6$}& $\vert\vert e_{K=0}\vert\vert_{L^{\infty}} $  & 8.22E-2 & 8.42E-2 & 2.44E-2 & 3.17E-2     \\
    \cline{3-7}
	&	&  $\vert\vert e_{K=0}\vert\vert_{L^2} $  & 1.46E-1 & 2.93E-1 & 7.78E-2 & 1.06E-1   \\
            \cline{3-7}
	&	&  $\vert\vert e_{K=4}\vert\vert_{L^\infty} $  & 2.67E-5 & 5.36E-7 & 4.34E-7 & 4.20E-7   \\
                \cline{3-7}
	&	&  $\vert\vert e_{K=4}\vert\vert_{L^2} $  & 4.45E-5 & 7.22E-7 & 4.65E-7 & 4.68E-7   \\
					\cline{3-7}
     & &  \added{time (s)}  & \added{500.4} & \added{665.7} & \added{1240.4} & \added{1488.6}   \\
                \cline{3-7}
      &  &  \added{Memory (GB)}  & \added{76.1} & \added{58.1} & \added{69.0} & \added{76.2}   \\
    \hline
		\end{tabular}
		\label{Table:6}
	\end{center}
\end{table}

\color{black}

\subsection{Two-Dimensional Poisson equation with corner near-singularity with L domain}\label{sec04:sec09}
 Consider the two-dimensional Poisson equation \eqref{eq:poissoneq} posed on the domain $\Omega=(-1,1)^2\backslash(-0.01,1) \times(-1,0.01)$, the exact solution is given by
\begin{equation}\label{eq:exactsolution9}
\begin{aligned}
\boldsymbol{\phi}=r^{\frac{2}{3}} \sin \frac{2\theta}{3}-\frac{1}{4}\left((x^2-1)+(y^2-1)\right).
\end{aligned}
\end{equation}

With the Dirichlet data $g(x, y)$ and source $f(x, y)$  derived from the analytical solution for consistency, we use $M_p = 6$ subdomains and $m = 7.5\times 10^4$ sampling points. Each subdomain uses $Q_n=Q_x\times Q_y = 91\times91$ collocation points, and Algorithm \ref{alg:Algorithm 1} is applied to compute the approximate solution $\boldsymbol{\phi}(x,y)$.

\begin{figure}[htp]
	\center
        {\subfigure[exact solution $\boldsymbol{\phi}$ ]{		 
            \includegraphics[scale=0.23]{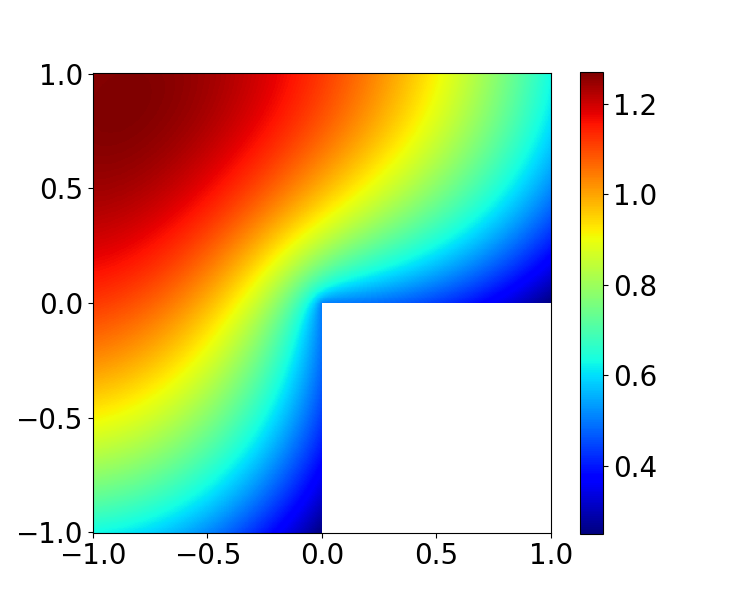}}
        \subfigure[approximte solution $\tilde{\boldsymbol{\phi}}_{0}  $ ]{		 
            \includegraphics[scale=0.23]{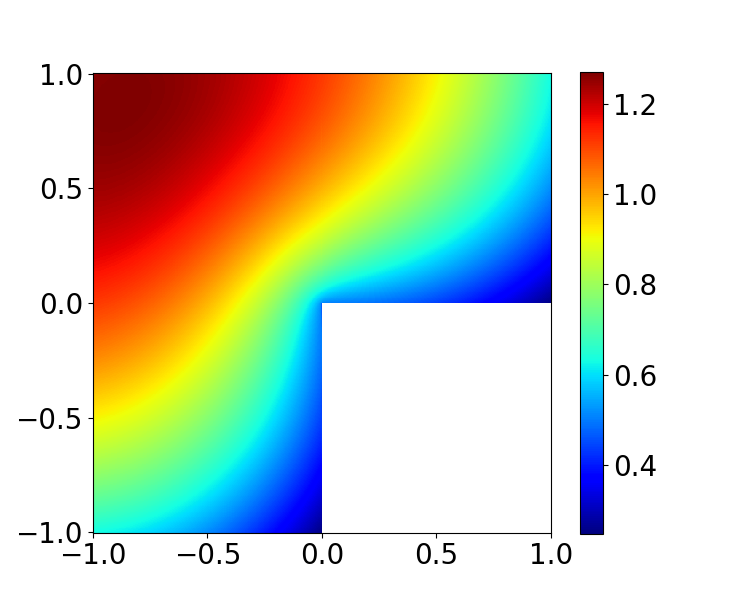}}
            \subfigure[local surface of  $\tilde{\phi}_0$ ]{		 
            \includegraphics[scale=0.32]{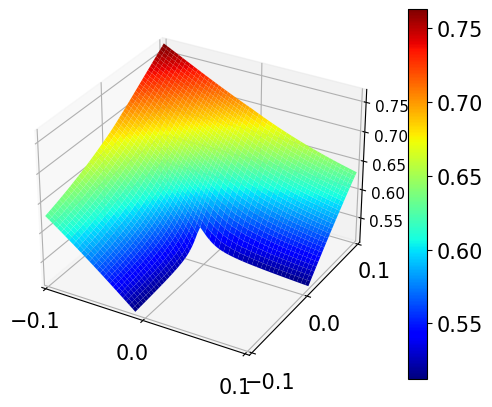}}
            \subfigure[Relative errors $\vert\vert e_{k}\vert\vert_{L^\infty}$ and $\vert\vert e_{k}\vert\vert_{L^2}$ ]{		
            \includegraphics[scale=0.191]{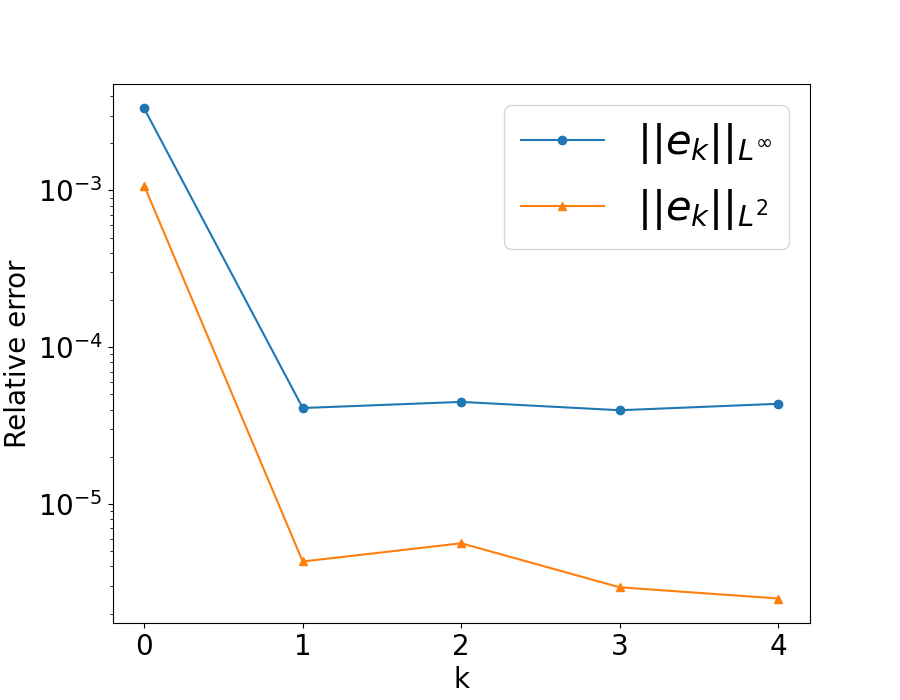}}
            \subfigure[approximte solution $\tilde{\boldsymbol{\phi}}_{K}  $  ]{		
            \includegraphics[scale=0.23]{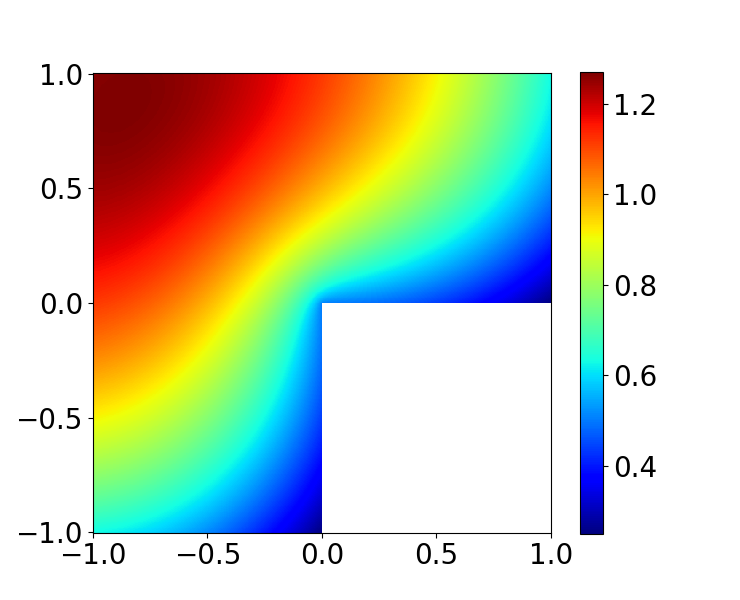}}
            \subfigure[local surface of  $\tilde{\phi}_K$ ]{		 
            \includegraphics[scale=0.32]{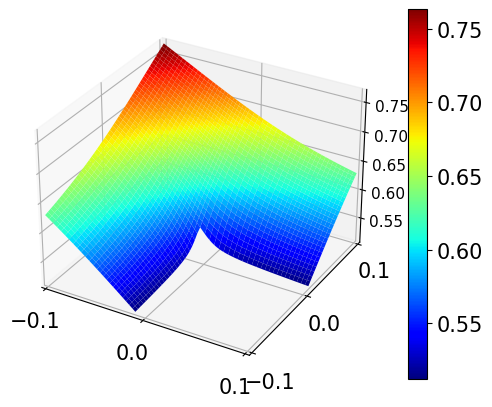}}
            }
	\caption{The exact solution and the numerical results for \eqref{eq:exactsolution9}  with $J_n=4000$ and $K=4$.}
	\label{fig:fig18}
\end{figure}

\begin{figure}[htp]
	\center
        {\subfigure[partition hyperplane density  $(k=0)$ ]{		 
            \includegraphics[scale=0.23]{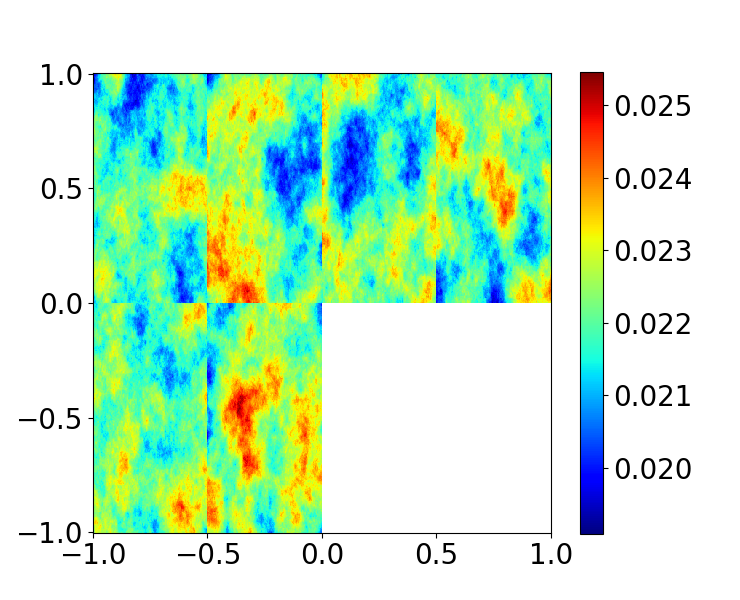}}
        \subfigure[collocation points $(k=0)$ ]{		 
            \includegraphics[scale=0.23]{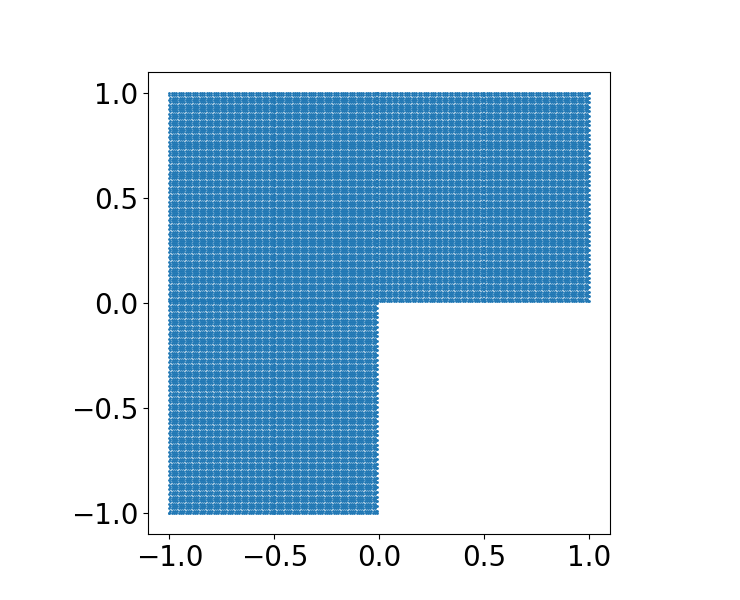}}
        \subfigure[shape parameter $(k=0)$ ]{		 
            \includegraphics[scale=0.23]{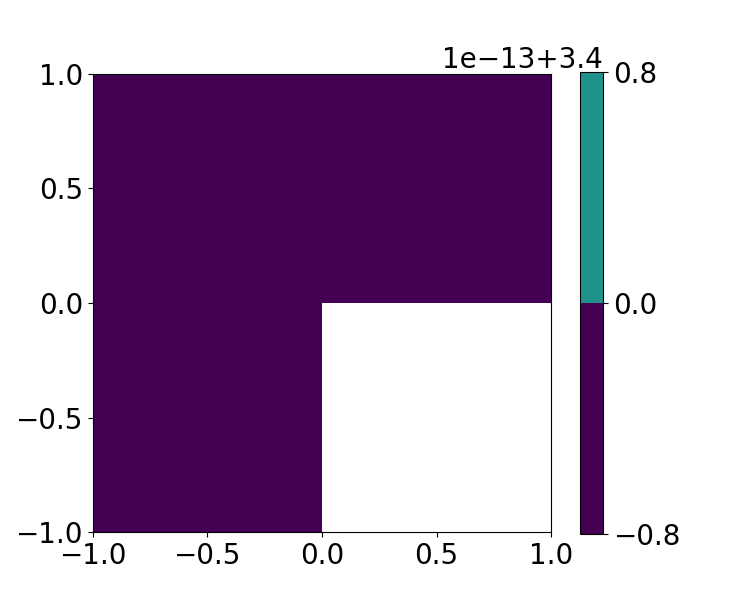}}
            \subfigure[partition hyperplane density  $(k=K)$  ]{		
            \includegraphics[scale=0.23]{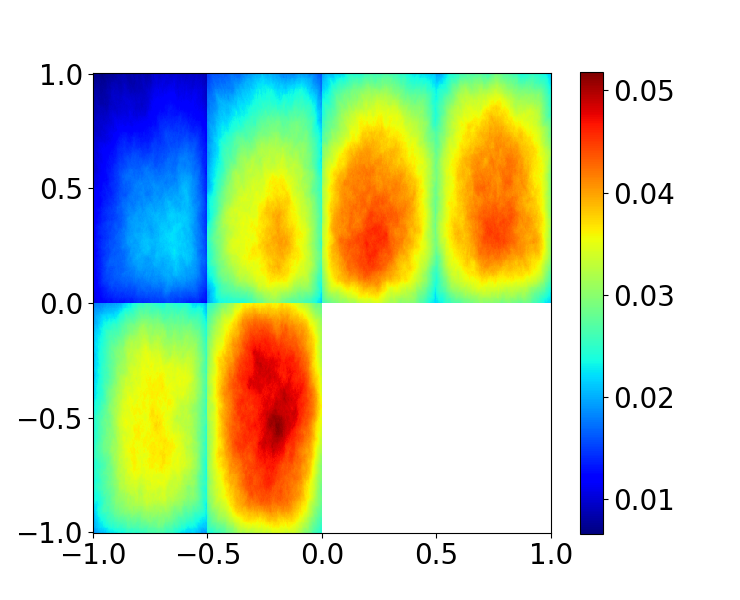}}
            \subfigure[collocation points $(k=K)$  ]{		
            \includegraphics[scale=0.23]{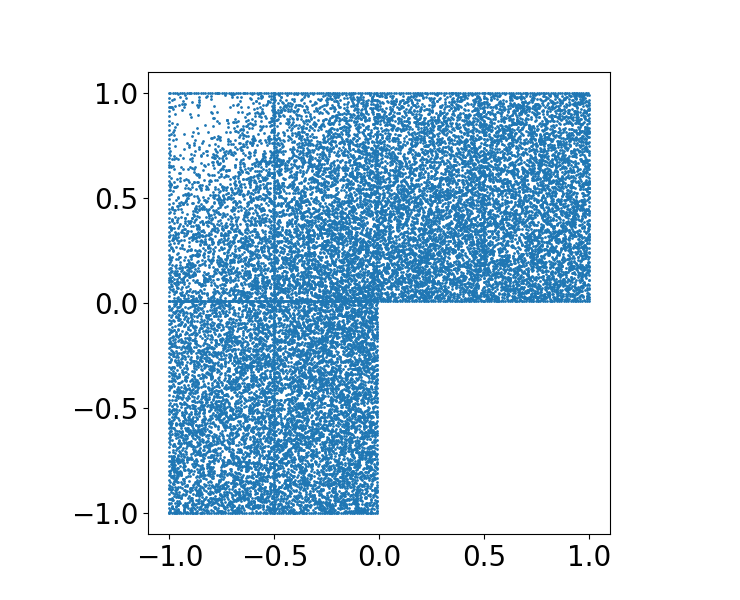}}
            \subfigure[shape parameter $(k=K)$ ]{		 
            \includegraphics[scale=0.23]{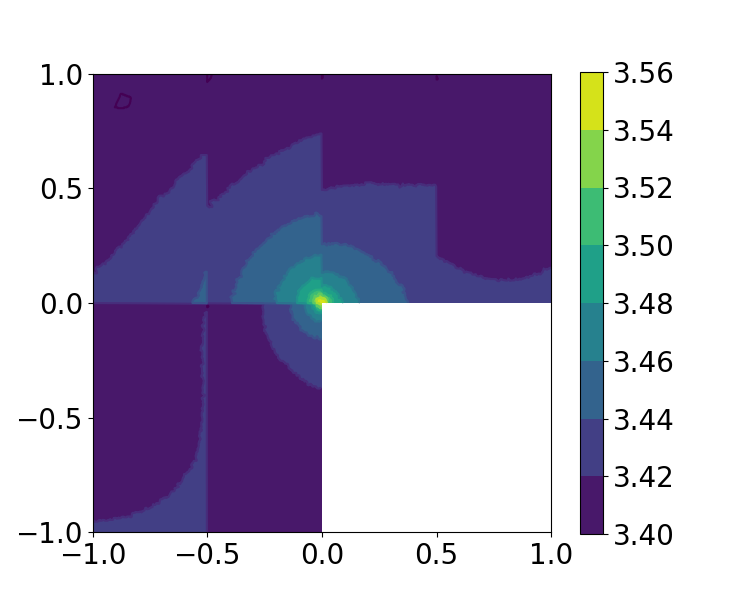}}
            }
	\caption{The partition hyperplane density, the collocation points and the shape parameters for \eqref{eq:exactsolution9} with $J_n=4000$ at  $k=0$ and $k=K=4$ iterations.}
	\label{fig:fig19}
\end{figure}

Table \ref{Table:9} summarizes the initial shape parameter $\gamma_n$ and relative errors of $\boldsymbol{\phi}(x,y)$ \deleted{for various $J_n$,} at both the initial iteration ($k=0$) and the final adapted iteration ($k=K=4$), \added{along with the computational time and peak GPU memory consumption of the full adaptive procedure for different $J_n$}. Figure \ref{fig:fig18} depicts the solution characteristics and adaptive progress by contrasting the exact solution against the initial and optimized approximations, along with their local surface profiles and the progression of error norms. \color{black}  Figure \ref{fig:fig19} shows that partition hyperplanes and collocation points achieve a certain degree of concentration in the corner region, and shape parameters are adjusted accordingly. The results indicate that AFCM can still improve accuracy through its adaptive mechanism under complex boundary conditions, though the extent of improvement is influenced by geometric complexity. This suggests that in regions with extreme singularities, local refinement or higher-order basis functions may be needed.

\color{black}

\begin{table}[htp]
	\centering 
	\caption{\deleted{The initial shape parameter $\gamma_n$ and the relative errors for \eqref{eq:exactsolution9}  with different $J_n$ at  $k=0$ and $k=K=4$ iterations}\added{The initial shape parameter $\gamma_n$, the relative errors at $k=0$ and $k=K=4$ iterations, the computational time, and the peak GPU memory consumption for \eqref{eq:exactsolution9} with different $J_n$}}
	\vspace{10pt}
	\small
	\begin{tabular*}{1.0\textwidth}{|l|l|l|l|*{4}{@{\extracolsep{\fill}}l|}}
		\hline
		\multirow{2}{*}{$M_p$}&\multirow{2}{*}{$Q_n$}&\multirow{2}{*}{$m$}& $J_n$   & 1500 & 2000  & 3000 & 4000      \\
		\cline{4-8}
		& & &  $\gamma_n$ &  2.0 & 2.6 & 2.8 & 3.4    \\
		\hline
		\multirow{6}{*}{$6$}&\multirow{6}{*}{$91\times 91$}&\multirow{6}{*}{$7.5\times 10^4$}& $\vert\vert e_{0}\vert\vert_{L^{\infty}} $   & 6.04E-4 & 3.17E-4 & 9.36E-4 & 3.33E-3     \\
		\cline{4-8}
		& & &  $\vert\vert e_{0}\vert\vert_{L^2} $   & 1.08E-5 & 2.06E-5 & 3.81E-4 & 1.07E-3 \\
		\cline{4-8}
		& & &  $\vert\vert e_{K=4}\vert\vert_{L^\infty} $   & 5.10E-4 & 3.01E-4 & 1.56E-4& 4.36E-5  \\
		\cline{4-8}
		& & &  $\vert\vert e_{K=4}\vert\vert_{L^2} $   & 1.06E-5 & 6.15E-6 & 4.91E-6 & 2.50E-6  \\
						\cline{4-8}
       & & &  \added{time (s)}  & \added{28.2} & \added{37.4} & \added{57.9} & \added{79.7}   \\
                \cline{4-8}
       & &  &  \added{Memory (GB)}  & \added{11.1} & \added{14.3} & \added{20.8} & \added{27.3}   \\
    \hline
	\end{tabular*}
	\label{Table:9}
\end{table}

\subsection{Three-Dimensional Poisson equation with a near singular solution with two \deleted{peak} \added{peaks}}\label{sec04:sec11}
Consider the three-dimensional Poisson equation in $\Omega=(0,0.5)^3$, the exact solution is given by
\begin{equation}\label{eq:exactsolution11}
\begin{aligned}
\boldsymbol{\phi}(x,y,z) =\ &e^{-1000\left((x-0.375)^2+(y-0.375)^2+(z-0.375)^2\right)}\\
&+e^{-1000\left((x-0.125)^2+(y-0.125)^2+(z-0.375)^2\right)}
\end{aligned}
\end{equation}

%The Dirichlet boundary condition $g(x, y, z)$ and source function $f(x, y, z)$ specified for the problem are derived directly from the provided analytical solution, ensuring consistency across all numerical experiments. In the present investigation, 
The computational domain is partitioned into $M_p = N_x\times N_y\times N_z=2\times 2\times 2$ subdomains, encompassing a total of $m = 4\times 10^5$ sampling points. Each subdomain is configured with $Q_n=Q_x\times Q_y \times Q_z= 20\times20\times 20$ collocation points, and Algorithm \ref{alg:Algorithm 1} is implemented to calculate the approximate solution 
$\boldsymbol{\phi}(x,y,z)$. For efficient GPU memory utilization,the batch-wise QR decomposition technique from \cite{pyRFM} is incorporated into the solution process.

\begin{figure}[htp]
	\center
        {\subfigure[exact solution $\boldsymbol{\phi}\ (z=0.375)$ ]{		 
            \includegraphics[scale=0.23]{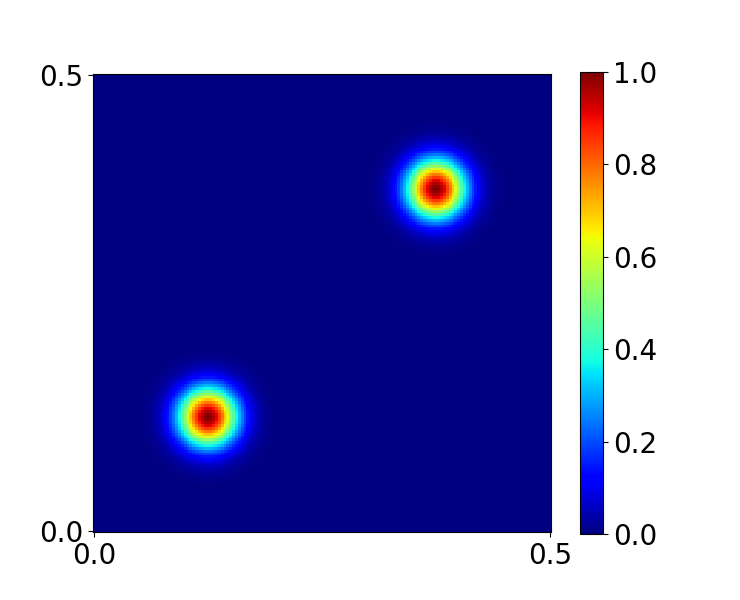}}
        \subfigure[approximte solution $\tilde{\boldsymbol{\phi}}_{0}\  (z=0.375)$ ]{		 
            \includegraphics[scale=0.23]{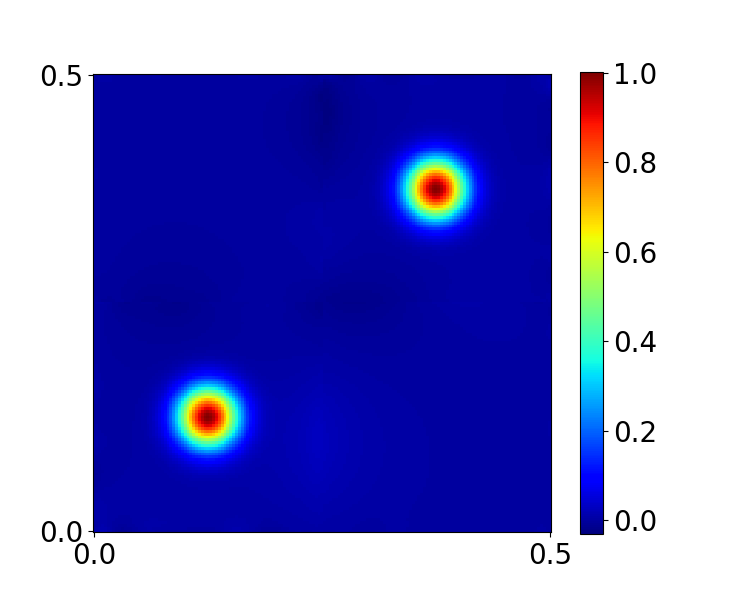}}
            \subfigure[local surface of  $\tilde{\phi}_0\ (z=0.375)$ ]{		 
            \includegraphics[scale=0.32]{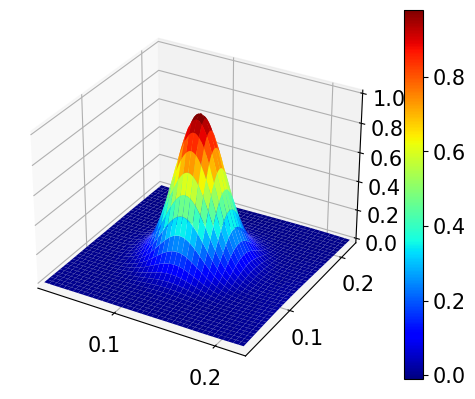}}
            \subfigure[Relative errors $\vert\vert e_{k}\vert\vert_{L^\infty}$ and $\vert\vert e_{k}\vert\vert_{L^2}$ ]{		
            \includegraphics[scale=0.191]{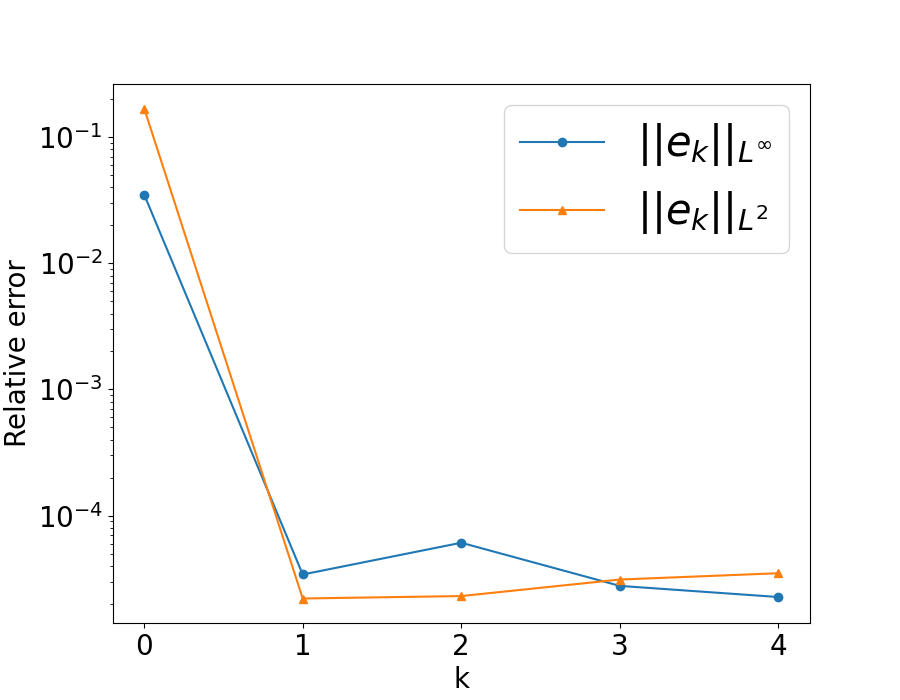}}
            \subfigure[approximte solution $\tilde{\boldsymbol{\phi}}_{K}\ (z=0.375) $  ]{		
            \includegraphics[scale=0.23]{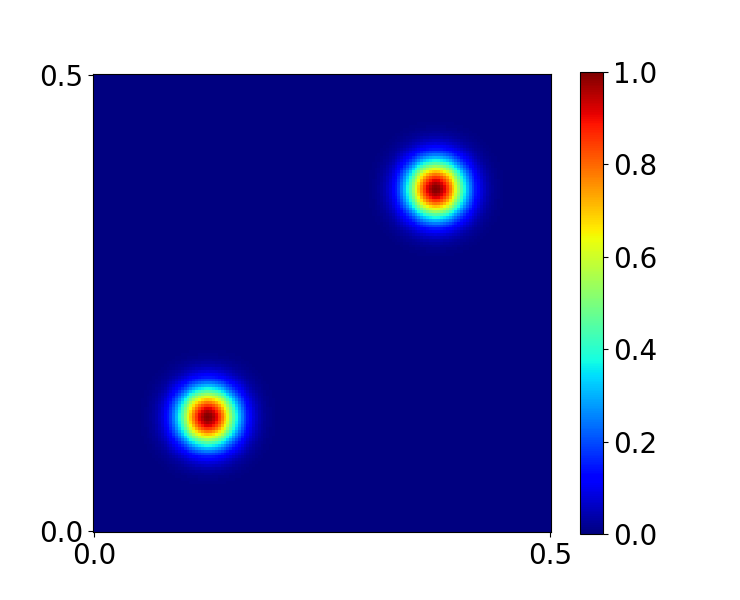}}
            \subfigure[local surface of  $\tilde{\phi}_K\ (z=0.375)$ ]{		 
            \includegraphics[scale=0.32]{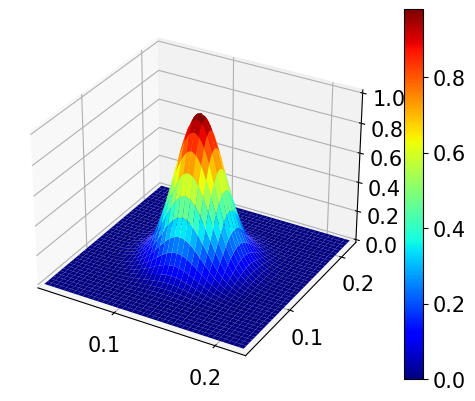}}
            }
	\caption{The exact solution and the numerical results for \eqref{eq:exactsolution11}  with $J_n=4000$ and $K=4$.}
	\label{fig:fig22}
\end{figure}

\begin{figure}[htp]
	\center
        {\subfigure[partition hyperplane density  $(k=0, z=0.375)$ ]{		 
            \includegraphics[scale=0.22]{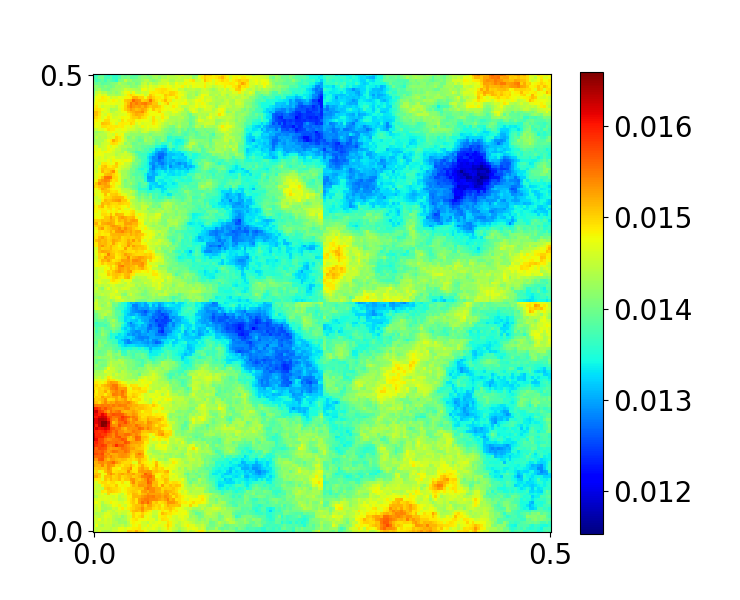}}
        \subfigure[collocation points $(k=0)$ ]{		 
            \includegraphics[scale=0.24]{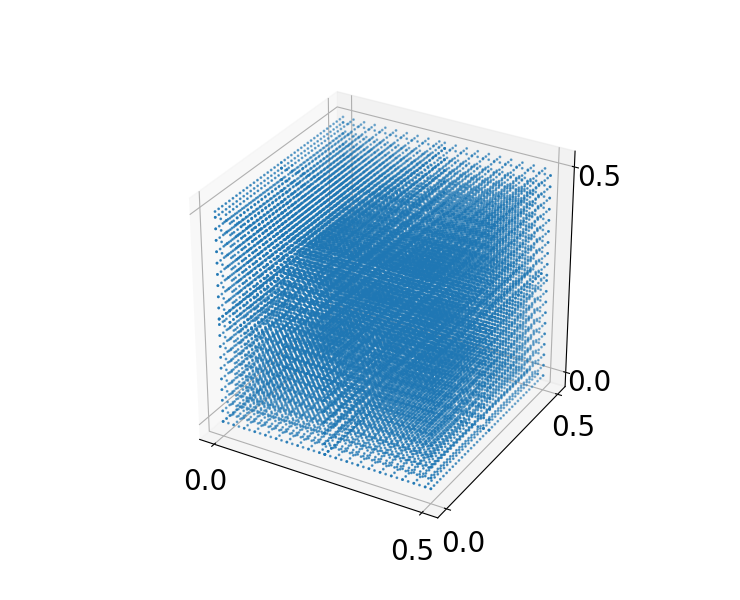}}
        \subfigure[shape parameter $(k=0, z=0.375)$ ]{		 
            \includegraphics[scale=0.22]{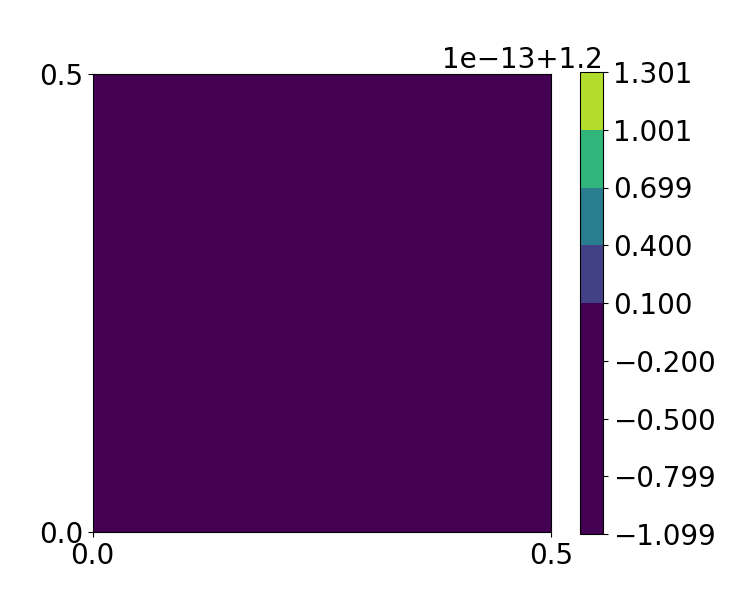}}
            \subfigure[partition hyperplane density  $(k=K, z=0.375)$  ]{		
            \includegraphics[scale=0.22]{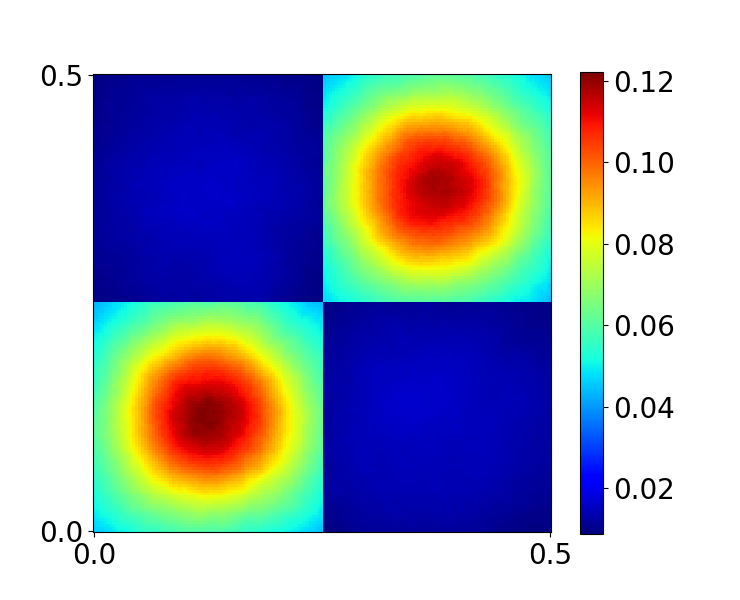}}
            \subfigure[collocation points $(k=K)$  ]{		
            \includegraphics[scale=0.24]{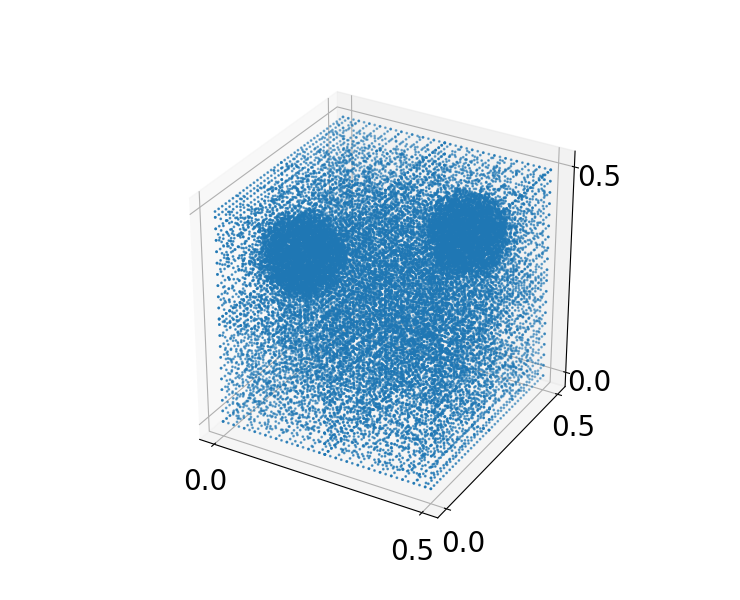}}
            \subfigure[shape parameter $(k=K, z=0.375)$ ]{		 
            \includegraphics[scale=0.22]{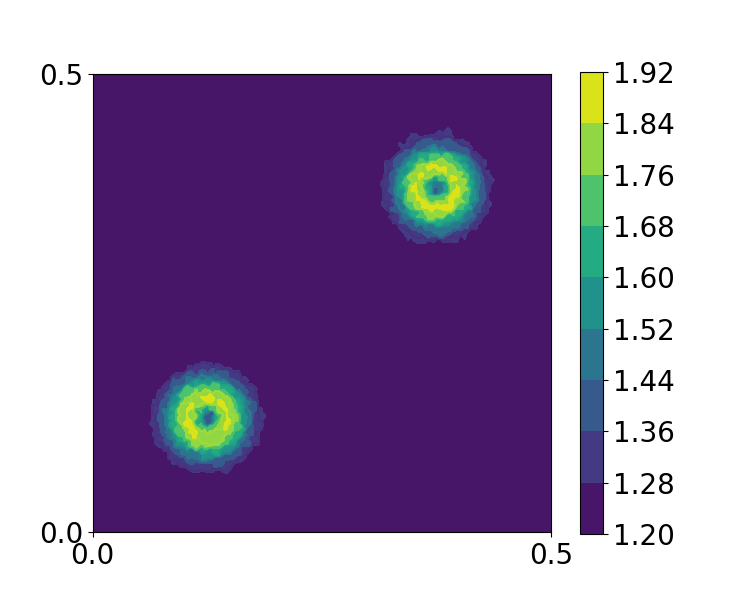}}
            }
	\caption{The partition hyperplane density, the collocation points and the shape parameters for \eqref{eq:exactsolution11} with $J_n=4000$ at  $k=0$ and $k=K=4$ iterations.}
	\label{fig:fig23}
\end{figure}

Table \ref{Table:11} reports the initial shape parameter $\gamma_n$ as well as the corresponding relative errors of $\boldsymbol{\phi}(x,y,z)$ \deleted{for different $J_n$ values}, at both the initial iteration ($k=0$) and the final adaptive step ($k=K=4$), \added{along with the computational time and peak GPU memory consumption of the full adaptive procedure for different $J_n$}. The solution characteristics and its adaptive process are illustrated in Figure \ref{fig:fig22}, which contrasts the exact solution, initial approximation, optimized outcome, their respective local surface profiles, and the iterative evolution of relative error norms. Furthermore, Figure \ref{fig:fig23} visualizes the hyperplane partition density, collocation point distribution, and shape parameter variation for $J_n=4000$ at the onset and conclusion of the adaptation process, thereby demonstrating the dynamic allocation of computational resources.

\begin{table}[htp]
	\begin{center} 
 \caption{\deleted{The initial shape parameter $\gamma_n$ and the relative errors for \eqref{eq:exactsolution11}  with different $J_n$ at  $k=0$ and $k=K=4$ iterations}\added{The initial shape parameter $\gamma_n$, the relative errors at $k=0$ and $k=K=4$ iterations, the computational time, and the peak GPU memory consumption for \eqref{eq:exactsolution11} with different $J_n$}}
  \vspace{10pt}
	\footnotesize 
	\begin{tabular*}{1.02\textwidth}{|l|l|l|l|*{4}{@{\extracolsep{\fill}}l|}} 
			\hline
            
		\multirow{2}{*}{$M_p$}&\multirow{2}{*}{$Q_n$}&\multirow{2}{*}{$m$}& $J_n$ & 1500 & 2000  & 3000 & 4000         \\
        \cline{4-8}
        &&&  $\gamma_n$ & 0.8 & 1.0  & 1.0 & 1.2       \\
			\hline
	\multirow{6}{*}{$2\times2\times2$}&\multirow{6}{*}{$20\times 20\times20$}&\multirow{6}{*}{$4\times 10^5$}& $\vert\vert e_{0}\vert\vert_{L^{\infty}} $  & 1.02E-0 & 3.48E-1 & 1.34E-1 & 3.51E-2     \\
    \cline{4-8}
		&&&  $\vert\vert e_{0}\vert\vert_{L^2} $  & 2.05E-0 & 7.25E-1 & 2.98E-1 & 1.69E-1   \\
            \cline{4-8}
		&&&  $\vert\vert e_{K=4}\vert\vert_{L^\infty} $  & 1.72E-2 & 1.14E-3 & 2.40E-4 & 2.26E-5   \\
                \cline{4-8}
		&&&  $\vert\vert e_{K=4}\vert\vert_{L^2} $  & 1.36E-2 & 7.50E-4 & 1.90E-4 & 3.49E-5   \\
									\cline{4-8}
       & & &  \added{time (s)}  & \added{788.7} & \added{1326.9} & \added{2146.6} & \added{4176.2}   \\
                \cline{4-8}
       & &  &  \added{Memory (GB)}  & \added{63.9} & \added{49.2} & \added{27.9} & \added{39.0}   \\
    \hline
		\end{tabular*}
		\label{Table:11}
	\end{center}
\end{table}

\color{black}

\section{Conclusions and remarks}
\label{sec05}
In this work, we propose the adaptive feature capture method (AFCM) based on RFM, a novel approach specifically developed to handle PDEs with near singular  solutions while maintaining high numerical accuracy without requiring additional computational resources and the prior information of the exact solutions. 
Building on TransNet's initialization framework, the AFCM employs the gradient norm of the approximate solutions of RFM as an indicator. This mechanism drives the partition hyperplanes of feature functions and collocation points to concentrate in regions characterized by larger solution gradients. Within these  regions, the method further enhances the steepness of feature functions' pre-activation values along the normal directions of the partition hyperplanes. 
Consequently, the AFCM achieves enhanced local expressive power in high-gradient regions, thereby yielding more accurate approximations. The AFCM repeats this adaptation process until no further improvement in the approximate solutions can be made. The proposed method delivers high accuracy in both space and time. 
As a mesh-free algorithm, it is inherently adaptable to complex geometric configurations and demonstrates efficacy in solving near-singular PDEs. 
We show a series of numerical experiments to validate our method. Results consistently confirm the stability and accuracy of the AFCM. Future work will explore extensions of this methodology to broader classes of PDEs and applications.

\section*{Acknowledgements}
\label{sec06}
X.-P. Wang acknowledges support from the National Natural Science Foundation of China (NSFC) (No. 12271461), the key project of NSFC (No. 12131010), Shenzhen Science and Technology Innovation Program (Grant: C10120230046, 10120240091), the Hetao Shenzhen-Hong Kong Science and Technology  Innovation Cooperation Zone Project (No. HZQSWS-KCCYB-2024016) and the University
Development Fund from The Chinese University of Hong Kong, Shenzhen (UDF01002028). The work of Q.L. He is partially supported by the National Natural Science Foundation of China (No. 12371434, No. U25A20200) and the National Key R \& D Program of China (No. 2022YFE03040002).

\bibliographystyle{elsarticle-num}
\bibliography{journal}

\end{document}